%% file: CategoricalArticle.tex
\documentclass{article}

\title{$\mathrm{LMod}_{R}(\V)$-enriched $\infty$-categories are left $R$-module objects of $\Cat^{\V}$ and $\Cat^{\V}$-enriched $\infty$-functors}
\author{Matteo Doni\footnote{Università degli Studi di Milano, Milan. Email address: matteo.doni@unimi.it}}
\input{CategoricalArticle.sty}
\begin{document}
%%% IMPOSTAZIONI TIPOGRAFICHE %%%
\newlength{\myindent} %Calcola e salva spazi di indentazione
\setlength{\myindent}{\parindent}
\parindent 0em %No indentazione pargrafi
\maketitle
\begin{abstract}
\subimport{./Chapters/}{Abstract}
\end{abstract}
\setcounter{tocdepth}{1}
\tableofcontents

\section{Introduction}
\subimport{./Chapters/}{Chap00-Preface2}

\section*{Acknowledgement}
\subimport{./Chapters/}{Chap00-Ringraziamenti}

\section*{Notation}
\subimport{./Chapters/}{Notation}
%%% CORPO %%%%
\subimport{./Chapters/}{RecallHigherCategory}

\subimport{./Chapters/}{RecallHigherAlgebra}

\subimport{./Chapters/}{RecallPresentableStable}

\subimport{./Chapters/}{EnrichedCategory}

\subimport{./Chapters/}{HigherCategoricalArticle}

\vskip 1cm

%\printbibliography 
%\bibliography{CategoricalArticle.bib}
%\bibliographystyle{alpha}
\bibliographystyle{unsrt}
\bibliography{CategoricalArticle}

\end{document}

%% file: Chapters/Abstract.tex
We establish the feasibility of investigating the theory of $\mathrm{LMod}_{R}(\V)$-enriched $\infty$-categories, where $R$ is an $\mathbb{E}_1$-ring of a presentable $\mathbb{E}_2$-monoidal $\infty$-category $\mathcal{V}$, through the framework of $\V$-enriched $\infty$-category theory. In particular, we prove that the $\infty$-category of $\mathrm{LMod}_{R}(\V)$-enriched $\infty$-categories $\mathcal{C}at_{\infty}^{\mathrm{LMod}_{R}(\V)}$ and the $\infty$-category of left $R$-module objects of the $\infty$-category of $\V$-enriched $\infty$-categories $\mathcal{C}at_{\infty}^{\V}$, $\mathrm{LMod}_{R}(\mathcal{C}at_{\infty}^{\V})$ are equivalent.
Moreover, if $R$ is an $\mathbb{E}_2$-ring of a presentable $\mathbb{E}_3$-monoidal $\infty$-category, we prove further that they are equivalent also to the $\infty$-category of $\mathcal{C}at_{\infty}^{\V}$-enriched $\infty$-functors $Fun^{\mathcal{C}at_{\infty}^{\V}}(B^{2}{R},\mathcal{C}at_{\infty}^{\V})$ are equivalent, where $B^{2}(-)$ is the \enquote{$2$-delooping} of a $(\infty,2)$-category with one object and one $1$-arrow.
The result has an interesting consequence.
When $R$ is an $\mathbb{E}_n$-ring of a presentable $\mathbb{E}_{n+1}$-monoidal $\infty$-category $\V$, it can be iterated to show that $(\infty,n)$-categories enriched in $\mathrm{LMod}_{R}(\V)$ are equivalent to $B^{n}R$-modules in $\V$-enriched $(\infty,n)$-categories, where $B^{n}$(-) is the \enquote{$n$-delooping}.

The specific case $\V=\Sp$ and $R=\Hk$, where $\Hk$ is the Eilenberg-Maclane Spectrum associated with a commutative and unitary ring $k$, holds significant value in derived algebraic geometry and spectral algebraic geometry. 
By specializing our results in this context we get that the $\infty$-categories $\mathcal{C}at_{\infty}^{\mathrm{Mod}_{\mathbb{H}\mathrm{k}}}$, $\mathrm{LMod}_{\mathbb{H}\mathrm{k}}(\mathcal{C}at_{\infty}^{\mathcal{S}p})$ and  $Fun^{\mathcal{C}at_{\infty}^{\mathcal{S}p}}(B^{2}\mathbb{H}\mathrm{k},\mathcal{C}at_{\infty}^{\mathcal{S}p})$ are equivalent.
Since $\Cat^{\LModk}$ is equivalent to the underlying $\infty$-category of the classical model category of dg-categories over $k$ with Dwyer-Kan equivalences, $\Cat^{\D(k)}$, this specific case provides us with two novel descriptions of the $\infty$-category of dg-categories over $k$; which is a fundamental object in derived algebraic geometry.

%% file: Chapters/Chap00-Preface2.tex
Categorification of classical homotopy theory led Quillen in the 1960s \cite{quillen2006homotopical} to define a more general homotopy theory. Quillen introduced model category theory, which allows us to formalize and work with various relationships. The basic object of study in this theory is a quadruple $(\C, \W, \mathscr{F}, \mathscr{C})$, where $\C$ is a category, $\W$ is a class of arrows in $\C$ called weak equivalences, and $\mathscr{F}$ and $\mathscr{C}$ are other classes of arrows providing additional structure to facilitate the theory's workings. However, these extra structures can make the theory cumbersome and complicated in proving certain results.

%\begin{example}[{}]
%\label{extop}
The category of topological spaces $\mathrm{Top}$ has two model structures whose classes of weak equivalences are composed of one of the canonical relations: weak homotopy equivalences or homotopy equivalence, \cite[Example 11.3.6]{riehl2014categorical}

Let $k$ be a commutative and unitary ring. The category of unbounded chain complexes $Ch(k)$ has two model structures whose classes of weak equivalences are composed of one of the canonical relations: quasi-isomorphisms or homotopy equivalence of chain complexes, \cite[Example 11.3.7]{riehl2014categorical}

The examples above show that some categories have two canonical interesting relations. The famous Whitehead theorem states that if we consider only CW-complexes, the two types of relations in $\mathrm{Top}$ coincide. Studying only CW-complexes instead of the entire category $\mathrm{Top}$ is not restrictive. A similar result holds for $Ch(k)$.

Despite the equivalence of the two relations in the above examples, where it matters, they have a different nature. The quasi-isomporphisms and weak homotopy equivalence are determined by an invariant, $H_{n}$, $\pi_{n}$, and have a homological essence. Instead, the homotopy equivalences have an enriched categorical essence. Indeed, in $\mathrm{Top}$ they are defined using the $\mathrm{Top}$-enrichment of $\mathrm{Top}$, and $Ch(k)$ using the $Ch(k)$-enrichment of $Ch(k)$, see \cite[\S Example 13.4.5]{riehl2014categorical}.

This dual nature of relations has many positive aspects; for example, it allows us to compute homotopy (co)limits using enriched category theory \cite[\S 7.7]{riehl2014categorical}. The need for a theory that accounts for these dual natures has led to the development of enriched model category theory.

In homotopy theory, the basic object of study is a pair $(\C, \W)$ with properties similar to the first two components mentioned above. Model category theory and enriched model category theory have the drawback of including two extra classes of arrows that can be an obstruction. To address this problem, some mathematicians (A. Joyal, C. Simpson, J. Lurie, J. Bergner, V. Hinich, R. Haugseng, D. Gepner, W. Dwyer, D. Kan, J. M. Boardman, R. M. Vogt, etc.) developed more flexible theories with the basic object being a pair $(\C, \W)$: the $\infty$-category theory and the enriched $\infty$-category theory.

The $\infty$-category theory is well-developed, but the enriched $\infty$-category theory is young and evolving.
The results of this paper are a small brick in this developing theory.

Considering hom-objects instead of hom-sets can seem daunting and might make one believe that working with them is more difficult. This complexity can sometimes simplify solving certain problems, which is why it is sometimes induced. This idea leads to defining stable homotopy theory \cite[\S 1.4]{HA}. The category $\mathrm{Top}$ is an $\mathrm{Top}$-enriched category, but its hom-$\mathrm{Top}$-objects do not have an algebraic structure. By inducing a topological algebraic structure on the hom-objects of $\mathrm{Top}$, we obtain the $\infty$-category of symmetric spectra $\Sp$. Just as for the category $\mathrm{Top}$, it is possible to define homotopy group invariants for $\Sp$, and since $\Sp$ has more structure, all functors land in the category of abelian groups $\pi_n: \Sp \to \Ab$. A naive way to think about the $\infty$-category $\Sp$ is that it bears the same relationship to the category of spaces as the ordinary category of abelian groups does to the ordinary category of sets.

Until now, we have been in the realm of algebraic topology, and to do algebraic geometry, we need to add further linearization. Let $k$ be a commutative and unitary ring. The Eilenberg-MacLane spectrum $\Hk$ is the internalization of $k$ inside the $\infty$-category $\Sp$. The $\infty$-category $\Sp$ forms a monoidal $\infty$-category with the famous smash product and spectrum sphere. So, it is possible to define the $\infty$-category of $\Hk$-modules inside $\Sp$, $\LModk$. Also, the homotopy group functors gain more linearity and land in the category of $k$-modules: $\pi_{n}: \LModk \to k$-$\mathrm{Mod}$.

In derived algebraic geometry, researchers work with categories whose hom-objects are defined up to some $k$-linear relation, so we are interested in studying the $\infty$-category of $\LModk$-enriched $\infty$-categories, $\Cat^{\LModk}$. The last statement can be made mathematically precise: the $\infty$-category of dg-categories over $k$, $\Cat^{\D(k)}$, and the $\infty$-category $\Cat^{\LModk}$ are equivalent, see \cite{DoniklinearMorita}.

The $k$-linearity in enrichment involves some difficulty. The $\Sp$-enrichment can be partially internalized, e.g., using stable $\infty$-category theory, or, as Heine proves in \cite[Theorem 7.21]{heine2023equivalence}, totally internalized using presentable stable $\infty$-category theory. However, the $\LModk$-enrichment does not have this property.

The aim of this paper is to prove that it is possible to study the $\infty$-category of $\LModk$-enriched $\infty$-categories, $\Cat^{\LModk}$, within the $\infty$-category of $\Sp$-enriched $\infty$-categories, $\Cat^{\Sp}$.

Using the canonical change of enrichment $U: \Cat^{\LModk} \to \Cat^{\Sp}$, it is possible to forget the action of $\Hk$ and consider the hom-objects of an $\LModk$-enriched $\infty$-category simply as hom-$\Sp$-objects.

To achieve our goal, we need to find a way to recover the forgotten action of $\Hk$ within the $\infty$-category $\Cat^{\Sp}$. In this paper, we will explore two ways to accomplish this.

The first method uses the higher algebra theory of modules within a symmetric monoidal $\infty$-category, \cite{HA}.

We start by internalizing the concept of commutative rings within $\Cat^{\Sp}$. Any ring $k$ can be considered as an $\Sp$-enriched $\infty$-category $\un{\Hk}$, which has only one object $\{*\}$, whose only non-trivial hom-object is $\Hom(*,*) = \Hk$ and whose composition is the multiplication in $k$, \Cref{defunderlineHigher}. Moreover, $\Cat^{\Sp}$ is a symmetric monoidal $\infty$-category, \Cref{notServePerIntro}.

The first result of this article is:

\begin{displayquote}
\textit{the $\infty$-category of $\LModk$-enriched categories $\Cat^{\LModk}$ and the $\infty$-category of left $\Hk$-module objects of the $\infty$-category of $\Sp$-enriched $\infty$-categories $\mathrm{LMod}_{\Hk}(\Cat^{\Sp})$ are equivalent.}
\end{displayquote}
Whose formalization is:

\begin{corollary}
\label{corIntroduzione3}
There exists an equivalence of $\infty$-categories 
\begin{equation}
\label{eqCatEnInLmodAreLmod}
\Cat^{\LModk}\simeq \mathrm{LMod}_{\Hk}(\Cat^{\Sp}).
\end{equation}
\end{corollary}

See \Cref{crCatEnInLmodAreLmodHKSp} for the rigorous statement.

Thus, we transition from a module structure in the enrichment of $\infty$-categories to $\infty$-categories endowed with a module structure themselves. 
This transition suggests that $\LModk$-enriched categories are $\Sp$-enriched functors. The vision of modules as enriched functors (and vice versa) is classical in representation theory and algebraic category theory. An incarnation of this analogy is the equivalence between the category of linear representations of $G$ on finite-dimensional $k$-vector spaces, denoted $Rep(G)$, and the category of comodules \cite{DeligneTanna, deligne2009ii}.

Furthermore, this analogy in the $\infty$-categorical setting has been extensively studied by V. Hinich \cite{HinichYon}. After rigorously defining all concepts, Hinich's work leads us to the second result, the proof of which is almost tautological.

Let $\2un{\Hk}$ denote the $\Sp$-enriched $\infty$-category with only one object $\{*\}$, where the only non-trivial hom-object is $\Hom(*,*)=\un{\Hk}$ (see \Cref{defDoubleFunctorHigherCateg}).

The second result can be articulated as follows:
\begin{displayquote}
\textit{the $\infty$-category of $\Cat^{\Sp}$-enriched $\infty$-functors $Fun^{\Cat^{\Sp}}(\2un{\Hk},\Cat^{\Sp})$ and the $\infty$-category of left $\Hk$-module objects in $\Cat^{\Sp}$ $\mathrm{LMod}_{\Hk}(\Cat^{\Sp})$ are equivalent.}
\end{displayquote}
Whose formalization is:

\begin{corollary}
\label{CorIntroduction2}
There exists a chain of equivalences of $\infty$-categories:
\begin{equation}
%\label{eqHkactionHKSp}
\mathrm{LMod}_{\Hk}(\Cat^{\Sp}) \simeq \mathrm{LMod}_{\un{\Hk}}(\Cat^{\Sp})\simeq  Fun^{\Cat^{\Sp}}(\underline{\un{\Hk}}, \Cat^{Sp}).
\end{equation}
\end{corollary}

See \Cref{crAction}.

Due to the underdevelopment of the theory of enriched $\infty$-categories, dealing with the $\infty$-category of enriched $\infty$-functors $Fun^{\Cat^{\Sp}}(\2un{\Hk}, \Cat^{\Sp})$ is not easy.
To make it more hands-on, we find a last description of $\Cat^{\LModk}$ by comparing it to a full $\infty$-subcategory of $Fun^{\Cat^{\Sp}}(\2un{\Hk}, \Cat^{\Sp})$.

Let $Alg$ be the $\infty$-category of associative symmetric spectra. This category is tautologically equivalent to the $\infty$-category of $\Sp$-enriched $\infty$-categories with only one object (actually, with contractible space of objects), \Cref{ChapEnrichCate}.

We take advantage of the particularity of the source in the second main result and we obtain the following:
\begin{displayquote}
\textit{the $\infty$-category of $Alg$-enriched $\infty$-functors from $\un{\Hk}$ to $\Cat^{\Sp}$ $Fun^{Alg}(\un{\Hk},\Cat^{\Sp})$ and the $\infty$-category of $\LModk$-enriched $\infty$-categories $\Cat^{\LModk}$ are equivalent.}
\end{displayquote}
Whose formalization is:

\begin{corollary}
\label{CorIntroductionAlg}
There exists a chain of equivalences of $\infty$-categories:
\begin{equation}
\Cat^{\LModk}\simeq \mathrm{LMod}_{\Hk}(\Cat^{\Sp})\simeq Fun^{Alg}(\un{\Hk},\Cat^{\Sp}).
\end{equation}
\end{corollary}
See \Cref{crCatLModkdescritionCalgtensored}.

%To summarize, these results can be stated as a single result as follows.
%\begin{theorem}
%\label{sumtheorem}
%Let $k$ be a commutative and unitary ring (also know as discrete $\E$-ring). There is a chain of equivalences of $\infty$-categories
%\[ \Cat^{\LModk}\simeq \mathrm{LMod}_{\Hk}(\Cat^{\Sp})\simeq Fun^{\Cat^{\Sp}}(\2un{\Hk},\Cat^{\Sp})\simeq Fun^{Alg}(\un{\Hk},\Cat^{\Sp}).\]
%\end{theorem}

The interest in \Cref{CorIntroductionAlg}, \Cref{CorIntroduction2} and \Cref{corIntroduzione3} come from derived algebraic geometry but the techniques of the proofs are completely higher algebraic categories. This leads the author to prove a more general results.

\begin{theorem}
\label{thCatEnInLmodAreLmodIntroduzione}
Let $(\V,-\otimes_{\V},\mathbb{I}_{\V}),$ be a presentable $\mathrm{E}_{2}$-monoidal $\infty$-category such that $\otimes_{\V}$ preserves small colimits componentwise. Let $R$ be an $\mathbb{E}_2$-algebra of $\V$,
then there exists an equivalence of $\infty$-categories 
\begin{equation}
\Cat^{\mathrm{LMod}_{R}(\V)}\simeq \mathrm{LMod}_{R}(\Cat^{\V}).
\end{equation}
\end{theorem}
See \Cref{thCatEnInLmodAreLmod}.

\begin{theorem}
Let $\V$ be a presentable $\mathbb{E}_3$-monoidal $\infty$-category and let $R$ be an $\mathbb{E}_2$-algebra of $\V$.

There exists a chain of equivalences of $\infty$-categories:
\begin{equation}
\mathrm{LMod}_{R}(\Cat^{\V}) \simeq \mathrm{LMod}_{\underline{R}}(\Cat^{\V})\simeq  Fun^{\Cat^{\V}}(\underline{\underline{R}}, \Cat^{\V}).
\end{equation}
\end{theorem}
see \Cref{thAction}; and  

\begin{proposition}
There exists a chain of equivalences of $\infty$-categories:
\begin{equation}
\Cat^{\mathrm{LMod}_{R}(\V)}\simeq \mathrm{LMod}_{R}(\Cat^{\V})\simeq Fun^{Alg(\V)}(\underline{R},\Cat^{\V}).
\end{equation}
\end{proposition}
See \Cref{eqCatLModkdescritionCalgtensored}.

In \cite[\S 6.3]{GepHauEnriched}, the authors build, for each $n\in\N_{\geq 0}\cup \{\infty\}$, the \enquote{n-delooping} $\infty$-functor. 
Let $\V$ be a $\mathbb{E}_{n+1}$-monoidal $\infty$-category, there exists a fully faithful and monoidal $\infty$-functors 
\begin{equation}
\label{eqBfunctor}
\mathcal{B}^{n}(-): Alg_{\mathbb{E}_n}(\V)\to \mathcal{C}at_{(\infty,n)}^{\V}: R\mapsto B^{n}.
\end{equation}

The case $n=1$ and $n=2$ are equivalent to the $\infty$-functors $\underline{(-)}$ and $\underline{\underline{(-)}}$ for the special case just mentioned.

The results \Cref{thCatEnInLmodAreLmodIntroduzione} can be iterated to obtain the following interesting corollary.

\begin{corollary}
Let $\V$ be a presentable $\mathbb{E}_{n+1}$-monoidal $\infty$-category with $n\in\N_{\geq 2}\cup \{\infty\}$ and let $R$ be $\mathbb{E}_{n+1}$-algebra of $\V$.
Then there exists an equivalence of $\infty$-category
\[   \mathcal{C}at_{(\infty,n)}^{\mathrm{LMod}_{R}(\V)}\simeq \mathrm{LMod}_{\mathcal{B}^{n}R}(\mathcal{C}at_{(\infty,n)}^{\V}).\] 
\end{corollary}
See \Cref{corIteratedResult}.

\section*{Related works}

%This work is part of a series of articles aimed at developing and comparing the Morita theory for $dg$-categories over $k$ (or $\D(k)$-enriched $\infty$-categories) and $k$-linear $\infty$-categories. 

This paper is the second of a series of articles, \cite{DoniCategorical, DoniklinearMorita}, aimed at developing and comparing the Morita theory for $dg$-categories over a commutative unitary ring $k$ (or $\D(k)$-enriched $\infty$-categories) and $k$-linear $\infty$-categories. 

We will not go as far in this paper as to define what Morita theory is, this happens in \cite{DoniklinearMorita}, we do not even deal with $dg$-categories.
The flavour of this paper is $\infty$-categorical and its target is to find the $\infty$-categorical analogue of the results in \cite{DoniCategorical}. 
Actually, we find the presentable $\infty$-categorial analogue of the result. Moreover, in $\infty$-categorical setting we need to deal with the fact that there are \enquote{infinite degree of associativity} between associativity and commutativity, so we need to deal with little cube operad, $\mathbb{E}_{n}$. 
I suspect it is possible to delete, or at least weaken the presentability hypothesis: maybe it is enough that $\V$ is an $\mathbb{E}_{2}$-monoidal (or $\mathbb{E}_3$-monoidal) $\infty$-category admitting geometric representation.  
$\\$

In my mind, this paper comes first and \cite{DoniCategorical} is an attempt to comprehend more deeply the results of this article.

At the end of this paper, a $k$-linearization will mean five things:
\begin{itemize}
    \item[(1)] a $\LModk$-enriched $\infty$-category;
    
    \item[(2)] or, a $\Cat^{\Sp}$-enriched $\infty$-functor from the $\Cat^{\Sp}$-enriched $\infty$-category $\2un{\Hk}$ to $\Cat^{\Sp}$;
    
    \item[(3)] or, a left $\un{\Hk}$-module object of $\Cat^{\Sp}$;

    \item[(4)] or, a $\Cat^{\Sp}$-enriched $\infty$-functor from the $Alg(\Sp)$-enriched $\infty$-category $\un{\Hk}$ to the $Alg(\Sp)$-left tensored $\infty$-category $\Cat^{\Sp}$;

    \item[(5)] or, a left $\Hk$-module object of $\Cat^{\Sp}$. 

\end{itemize}

$(1)$ represents the more classical derived algebro geometric view, \cite[\S 5]{DoniklinearMorita}.

%$(1)$ is the more classical derived algebraic geometry view, as we explain in \cite[\S 5]{DoniklinearMorita}.

$(4)$ is the starting point for a general notion of $k$-linearization. This description is important for us because it will allow us to give the right definition of $k$-linearization in the context of spectral algebraic geometry, \cite[\S 6]{DoniklinearMorita}.

\section*{Structure of the paper}
Initially, this article was intended to be self-contained but the technicalities of the $\infty$-category theory made this feat impossible.  For example, it lacks an introduction to the theory of $\infty$-operad.
So we opted for an article in which each concept is at least explained in words and its technical aspects are referred to in a few citations. 

It is split into two parts.
The first part is composed by \Cref{subsecHigherCategoryTheory,subsectionHigherAlgebraicStructure,subsecStablepresentable,ChapEnrichCate}. Here, we define, recall, and construct all mathematical concepts and structures needed to state and prove results. For a reader used to the work of Lurie, \cite{HTT,HA}, the \Cref{subsecHigherCategoryTheory,subsectionHigherAlgebraicStructure,subsecStablepresentable} can be skipped.

In \Cref{ChapEnrichCate}, we treat and present the young theory of $\infty$-enriched categories. This theory is far from being well-developed and with a commonly agreed standard notation. So in this part, we need to make some choices; for example, which of the three definitions of enriched $\infty$-categories to use.

The last part is composed only by \Cref{secLModEnrichedareModuleObjects}. Here, we enunciate and prove all the promised results in the introduction in the general setting and in \Cref{subsecSpecificCase} we specialized them to the specific case Of interest in derived algebraic geometry.

%% file: Chapters/Chap00-Ringraziamenti.tex
I extend my gratitude to Paolo Stellari for his guidance throughout the recent years, which has enabled me to consolidate my ideas.
I would like to thank Hinich and Haugseng for his assistance with my question, which helped me identify a mistake. In particular, the latter suggested me to generalize my former results, \Cref{crCatEnInLmodAreLmodHKSp,crAction}, into their final version that appears in this paper. I would also like to thank Devon Stockall for pointing out an error in Theorem 6.1 and, as a consequence, in Corollary 6.3.

%% file: Chapters/Notation.tex
\begin{itemize}
   
%    \item Following Lurie's convention, we use the term $\infty$-category to refer to a $(\infty,1)$-category. Among the various equivalent models of $\infty$-categories, we adopt the model of quasi-categories, as used by Lurie in \cite{HA}, \cite{HTT}, and \cite{SAG}. In certain instances, Riehl-Verity's $\infty$-cosmos theory will be more convenient, and we denote it with bold type, for example, $\mathcal{Q}cat$ for the (Lurie) $\infty$-category of quasi-categories and \textbf{Qcat} for the $\infty$-cosmos of quasi-categories. Note that in many cases, $\mathcal{Q}cat$ simplifies to $\Cat$, and it is an object in \textbf{Qcat}. In \cite{RiehlElements}, the authors introduce the \textit{model-independent} $\infty$-cosmos of $(\infty,1)$-categories \cite[Definition 1.3.10]{RiehlElements}. Although this generality can lead to complications we prefer the specific model \textbf{QCat}. However, we acknowledge that this work could have been formulated in the $\infty$-cosmos of $(\infty,1)$-categories;

 %   \item let $\mathbf{K}$ be an $\infty$-cosmos; we use the notation $\mathfrak{H}\mathbf{K}$ to denote its homotopy $2$-category;

    \item we adhere to Lurie's homological convention, as described in \cite{HA};

    \item in this work, an \textit{$\infty$-category} is an object in the $\infty$-cosmos of $(\infty,1)$-categories \cite[Definition 1.3.10]{RiehlElements};
    
    \item in \cite{RiehlElements,kerodon} the authors call \textit{elements of the $\infty$-category $\A$} a map $\mathbb{1}\to \A $ from the terminal object of a $\infty$-cosmos to an $\infty$-category $\A$. Instead, in \cite{HTT} the author calls \textit{object of the $\infty$-category $\A$} a map $\mathbb{1}\to \A $. Since both notations were appropriate, we decided that in this paper the terms \textit{element} and \textit{object} of an $\infty$-category are used interchangeably  (actually, element is more general);

    \item we use both $Fun(\C,\D)$ and $\D^{\C}$ to denote the $\infty$-category of $\infty$-functors between two $\infty$-categories $\C$ and $\D$;% For a pair of $\infty$-functors $F,G:\C\to\D$, we denote the $\infty$-category of natural transformations by $NT_{\D^{\C}}(F,G)$ or, when clear from context, simply $NT(F,G)$;

    \item  with a slight abuse of notation, we denote by $\pi_{0}:\mathcal{S}_{*}\to Set$ the composition of $\infty$-functors $U\pi_{0}$ where $\pi_{0}:\mathcal{S}\to Set$ is the component-of-space $\infty$-functor,\cite[Example 3.6.2]{RiehlCTIC}, and $U:\mathcal{S}_{*}\to\mathcal{S}$ is the canonical forgetful $\infty$-functor;

%    \item the symbol $:=$ signifies \enquote{defined as}. For instance, if $\C$ is a small $\infty$-category, writing $P(\C):=Fun(\C^{op},\mathcal{S})$ is equivalent to stating that $P(\C)$ is defined as the $\infty$-category of presheaves on $\C$;

    \item unless otherwise indicated by the symbol $\dashv$, the following conventions will be used: in a horizontally written adjunction, the left (right) adjoint has the arrowhead pointing to the right (left);

%    \item throughout this paper, $\simeq$ signifies \enquote{equivalence in the same setting}, while $\cong$ means \enquote{isomorphism in a category}. These conventions apply to Quillen equivalences of model categories, equivalences of ordinary categories, equivalences in $\infty$-categories, etc.
%     Inside an $\infty$-category the words isomorphism and the word equivalence have the same meaning, this is because different reference books have different conventions, \cite{RiehlElements},\cite{kerodon},\cite{HTT}. So, we decided that in this paper the terms \textit{isomorphism} and \textit{equivalence} are used interchangeably inside an $\infty$-category;

%    \item    The inclusion $\mathcal{S}\to \mathcal{C}at_\infty$ has both a left and a right adjoint. The right adjoint, denoted \[(-)^{eq}:\mathcal{C}at_\infty\to\mathcal{S},\] is referred to as \enquote{the maximal subspace} since it discards the non-invertible arrows in the $\infty$-category.
%    In \cite{HTT}, the author referees it by $(-)^{\simeq}$;

    \item We use the symbol $\mapsto$ to indicate where a functor sends objects or where an $\infty$-functor, up to equivalence, maps objects. For example, $-+1:\Z\to\Z:n\mapsto n+1$.

\end{itemize}

%% file: Chapters/RecallHigherCategory.tex
\section{Recalls: higher category theory}
\label{subsecHigherCategoryTheory}

The famous Barr-Beck-Lurie Theorem, \cite[Theorem 4.7.3.5]{HA}, is crucial in $\infty$-category theory because, nowadays, using it is the only strategy to prove that an $\infty$-category is equivalent to an $\infty$-category of algebra.
So, it is not surprising that the strategy of the proof of \Cref{thCatEnInLmodAreLmod} employs the Barr-Beck-Lurie Theorem.

This section is largely dedicated to recalling some $\infty$-categorical concepts in this famous theorem.
Let us begin by establishing some notations and recalling the useful comma $\infty$-category.

\begin{notation}
    $\Delta$ is the simplicial indexing category whose objects are the non-empty finite ordered sets $[n]:=\{0,1,2,...,n\}$ and whose morphisms are order-preserving functions between them. When we consider $[n]$ as an $\infty$-category (i.e., via the nerve), we write $\mathbb{n+1}$.
\end{notation}
 
\begin{notation}
We need to clarify our understanding of a diagram of $\infty$-categories $\A$. A diagram is in the sense of Riehl-Verity, \cite[Definition 2.3.1]{RiehlElements}: a morphism $F:D\to\A$ whose source is a simplicial set. In a cartesian closed $\infty$-cosmos, a diagram is equivalent to an $\infty$-functor between $\infty$-categories. Since the $\infty$-cosmos of quasi-categories is cartesian closed and we work with this model of $(\infty,1)$-categories, we use the terms \enquote{diagram} and \enquote{ $\infty$-functor} interchangeably.  
\end{notation}

Let $\C$ be an $\infty$-category, and let $j:\mathbb{1}\coprod \mathbb{1}\to\mathbb{2}$ denote the canonical inclusion of components. We define \textit{the $\infty$-category of arrows of $\C$} as the simplicial cotensor $\C^{\mathbb{2}} = Fun(\mathbb{2},\C)$. This construction is associated with the isofibration 
\[Fun(-,\C)(j):\C^{\mathbb{2}}\to \C^{\mathbb{1}\coprod\mathbb{1}}\simeq \C\times\C;\]
we denote this isofibration by $(s,t)$,  and we call the corestriction to the first component \textit{source} and to the second component \textit{target}.

Consider a cospan $f:\mathcal{A}\to \C\leftarrow \mathcal{B}:g$ in $\Cat$, or in the $\infty$-cosmos $\textbf{QCat}$. 

We recall that the comma $\infty$-category (\cite[Definition 3.4.1]{RiehlElements}) is defined as follows: 
\begin{equation*}
    \begin{split}
    \Hom_C (f,g):= &\;pullback(\mathcal{A}\times \mathcal{B}\xrightarrow{f\times g} \C\times \C\xleftarrow{(s,t)} \C^\mathbb{2}) \\ \simeq & \; \mathcal{A} \underset{Fun(\{0\},\C)}{\times} Fun(\mathbb{2}, \C) \underset{Fun(\{1\},\C)}{\times} \mathcal{B}.
    \end{split}
\end{equation*}

It is important to note that this is the same as Lurie's oriented fiber product $A\underset{C}{\tilde{\times}} B$, \cite[Definition 4.6.4.1]{kerodon}. 

Since comma $\infty$-categories are by definition pullbacks of isofibrations $\C^\mathbb{2}\to \C\times \C$, the comma $\infty$-categories are equipped with an isofibration $(p_1,p_0):\Hom_\C (f,g)\to \mathcal{B}\times \mathcal{A}$.

Let $\mathcal{A,B}$ be two $\infty$-categories and $\C\to \mathcal{A\times B}$ and $\D\to \mathcal{A\times B}$ two isofibrations (which is a module \cite[Example 7.2.8]{RiehlElements} from $\mathcal{B}$ to $\mathcal{A}$: cartesian on the left, $p_1$ cartesian, and cocartesian on the right, $p_0$ cocartesian). In short, there is a comma span  
\[ \mathcal{B}\xleftarrow[cocart]{p_1}\Hom_{\C}( f,g  )\xrightarrow[cart]{p_0}\mathcal{A}.\]

The comma $\infty$-category has more structure. It is associated with a natural transformation and it has a weak universal property:

% https://q.uiver.app/?q=WzAsNCxbMCwxLCJCIl0sWzIsMSwiQSJdLFsxLDIsIkMiXSxbMSwwLCJIb21fQyhmLGcpIl0sWzEsMiwiZiJdLFswLDIsImciLDJdLFszLDAsInBfMSwgY29jYXJ0IiwyXSxbMywxLCJwXzIsIGNhcnQiXSxbNCw2LCIiLDIseyJzaG9ydGVuIjp7InNvdXJjZSI6MjAsInRhcmdldCI6MjB9fV1d
\[\begin{tikzcd}
	& {\Hom_C(f,g)} \\
	B && A \\
	& {C.}
	\arrow[""{name=0, anchor=center, inner sep=0}, "f", from=2-3, to=3-2]
	\arrow["g"', from=2-1, to=3-2]
	\arrow[""{name=1, anchor=center, inner sep=0}, "{p_1, cocart}"', from=1-2, to=2-1]
	\arrow["{p_0, cart}", from=1-2, to=2-3]
	\arrow[shorten <=9pt, shorten >=9pt, Rightarrow, from=0, to=1]
\end{tikzcd}\]

\begin{example}
Let $\C$ be an $\infty$-category and let $c,d\in \C$ be two objects.
The comma $\infty$-category $\Hom_{\C}(c,d)$ associated to the cospan $\mathbb{1}\xrightarrow{c}\C\xleftarrow{d}\mathbb{1}$ is the mapping space $\text{Map}_{\C}(c,d)$.
So, we denote the mapping space by $\Hom_{\C}(c,d)$.
\end{example}

\begin{definition}[{\cite[Definition 3.2.7]{RiehlElements}}]
\label{defFiberEquivalence}
Let $m:\C\to \A\times\B$ and $n:\D\to \A\times \B$ be two $\infty$-functors in $\textbf{QCat}$. The notation $m\simeq_{\mathcal{A}\times \mathcal{B}}n$ (or $\C\simeq_{\mathcal{A}\times \mathcal{B}}\D$) means that they are equivalent in the $\infty$-cosmos of $\textbf{QCat}_{/\mathcal{A}\times \mathcal{B}}$. In this case, we will say that they are \textit{fibered equivalent}.   
\end{definition}

\begin{remark}[{\cite[Proposition 12.2.14, Proposition 7.4.6]{RiehlElements}}]
\label{rmkfiberequivalencesFiberwise}
 \Cref{defFiberEquivalence} is more general than what we need. We deal only with fibered equivalences between modules in the $\infty$-cosmos of ($\infty,1$)-categories (in this case also including the equivalences between comma $\infty$-categories). In this $\infty$-cosmos, for two modules being fibered equivalence is the same as being fiberwise equivalence:
 for each each objects $a\in\A$ and $b\in\B$:
the induced functor between the fibers $F_{a,b}:\C_{(a,b)}\to \D_{(a,b)}$ is an
equivalence of $\infty$-categories.
\end{remark}

\begin{notation}
\label{notEleComma}
Let $f:\mathcal{A}\to \C \leftarrow \mathcal{B}:g$ be a cospan of $\infty$-categories.
With a triple $(a,\alpha:f(a)\to g(b),b)$, we denote the element of the comma $\infty$-category $\mathbb{1}\to \Hom_{\C}(f,g)$ induced by the diagram:
% https://q.uiver.app/#q=WzAsNCxbMSwwLCJcXG1hdGhiYnsxfSJdLFswLDEsIlxcYnVsbGV0Il0sWzIsMSwiXFxidWxsZXQiXSxbMSwyLCJcXGJ1bGxldCJdLFswLDEsImMiLDJdLFswLDIsImEiXSxbMSwzLCJnIiwyXSxbMiwzLCJmIl0sWzcsNCwiYiIsMix7InNob3J0ZW4iOnsic291cmNlIjoyMCwidGFyZ2V0IjoyMH19XV0=
\[\begin{tikzcd}
	& {\mathbb{1}} \\
	B && A \\
	& {\C,}
	\arrow[""{name=0, anchor=center, inner sep=0}, "b"', from=1-2, to=2-1]
	\arrow["a", from=1-2, to=2-3]
	\arrow["g"', from=2-1, to=3-2]
	\arrow[""{name=1, anchor=center, inner sep=0}, "f", from=2-3, to=3-2]
	\arrow["\alpha"', shorten <=7pt, shorten >=7pt, Rightarrow, from=1, to=0]
\end{tikzcd}\]
by universal property. In other words, an element of $\Hom_{\C}(f,g)$ is a triple $(a,\alpha:f(a)\to g(b),b)$, where $a$ and $b$ are elements of $A$ and $B$ respectively, and $\alpha:f(a)\to g(b)$ is an arrow of $\C$.

In the case $B\simeq\mathbb{1}$, we are in the case of a slice $\infty$-category $\C_{/g}$ and an element is a triple $(a,\alpha:f(a)\to g,\mathbb{1})$ with the terminal $\infty$-category as the last component. We make an abuse of notation and denote this element by a pair consisting only of its first two components $(a,\alpha:f(a)\to g)$. The morphisms of a comma $\infty$-category also have a description as a triple. Moreover, for the coslice case, we use the same convention.
\end{notation}

The comma $\infty$-category can be used to describe some universal properties; for instance, see \Cref{defLimitColimit}.
The next definition explains how a comma $\infty$-category detects an universal property.

\begin{definition}[{\cite[Definition 3.5.2]{RiehlElements}}]
Given a cospan $C\xrightarrow{g} \A\xleftarrow{f}\B$, the comma $\infty$-category $\Hom_{\A}(f,g)\to \C\times \B$ is \textit{left
representable} if there exists a $\infty$-functor $l:\B\to\C$ so that
\[\Hom_{\A}(f,g)\simeq_{\C\times\B}\Hom_{\C}(l,\C);\]
and \textit{right representable} if there exists a $\infty$-functor $r:\C\to\B$ so that
\[\Hom_{\A}(f,g)\simeq_{\C\times\B}\Hom_{\B}(\B,r).\]
\end{definition}

One can use the comma construction to define the $\infty$-category of cones, \cite[Definition 4.2.1]{RiehlElements}. Let $\C$ be an $\infty$-category, $J$ be a simplicial set, $d:\mathbb{1}\to Fun(J,\C) $ be a diagram, and $\cDelta :\C\to Fun(K,\Delta )$ be the constant functor. 
\begin{notation}
The \textit{$\infty$-category of cones under }(\textit{over}) $d$ is $\Hom_{Fun(J,\C)}(d,\cDelta)$ ($\Hom_{Fun(J,\C)}(\cDelta ,d)$). 
\end{notation}
Just as in the classical case, we can characterize the (co)limits with a universal property.  

\begin{definition}[{\cite[Proposition 4.3.1]{RiehlElements}, co/limits represent cones}]
\label{DefLimCol}
A diagram $d:\mathbb{1}\to Fun(J,\C) $ admits a limit if and only if the $\infty$-category of cones $\Hom_{Fun(J,\C)}(d,\cDelta)$ over $d$ is right representable
\[\Hom_{Fun(J,\C)} (\cDelta, d) \simeq_{\mathbb{1}\times\C} \Hom_\C (\C, l),\]
in which case the representing functor $l:\mathbb{1}\to \C$ defines the limit functor. We call $l$ the \textit{apex} of the limit. Dually, a diagram $d:\mathbb{1}\to Fun(J,\C)$ admits a colimit if and only if the $\infty$-category of cones $\Hom_{Fun(J,\C)}(d,\cDelta)$ under $d$ is left representable
\[\Hom_{Fun(J,\C)} (d, \cDelta) \simeq_{\C\times \mathbb{1}} \Hom_\C(c, \C),\]
in which case the representing functor $c:\mathbb{1}\to\C$ defines the colimit functor. We call $c$ the \textit{nadir} of the colimit.
\end{definition}

We define limits and colimits cones in the $\infty$-categorical setting. The following proposition contains the desired definition simultaneously.

\begin{proposition}[{\cite[Proposition 4.3.2]{RiehlElements}}]
\label{defLimitColimit}
A diagram $d: \mathbb{1}\to Fun(J,\A)$ of shape $J$ in an $\infty$-category $\A$:
\begin{itemize}
    \item[(i)] admits a limit if and only if the $\infty$-category $\Hom_{Fun(J,\A)}(\cDelta,d )$ of cones over $d$ admits a terminal element,
    in which case the terminal element defines a \textit{limit cone}, and
    \item[(ii)] admits a colimit if and only if the $\infty$-category $\Hom_{Fun(J,\A)}(d,\cDelta)$ of cones under $d$ admits an initial element,
    in which case the initial element defines the \textit{colimit cone}.
\end{itemize}
\end{proposition}

\begin{notation}
Let $d:J\to \A$ be a diagram of shape $J$ in a $\infty$-category $\A$, let $\hat{d}:\cDelta\to d$ be a cone (over) under $d$ and let $F:\A\to \B$  be a $\infty$-functor between two $\infty$-categories.
We call the \textit{the image of cone $d$ via $F$} the cone (over) under $Fd$ which is
the image of the $\infty$-functor \[\hat{F}:\Hom_{Fun(J,\A)}(\cDelta, d)\to \Hom_{Fun(J,\B)}(\cDelta, Fd),\]
induced by the  universal property of pullbacks:
% https://q.uiver.app/#q=WzAsNixbMCwwLCJcXEhvbV97XFxBXntKfX0oXFxjRGVsdGEsZCkiXSxbMSwxLCJcXEhvbV97XFxCXntKfX0oXFxjRGVsdGEsZCkiXSxbMywxLCJcXEJee0pcXHRpbWVzXFxtYXRoYmJ7Mn19Il0sWzEsMywiXFxjRGVsdGFcXHRpbWVzIEoiXSxbMywzLCJcXEJcXHRpbWVzXFxCIl0sWzIsMCwiXFxBXntKXFx0aW1lc1xcbWF0aGJiezJ9fSJdLFswLDEsIiIsMCx7InN0eWxlIjp7ImJvZHkiOnsibmFtZSI6ImRhc2hlZCJ9fX1dLFsxLDJdLFsxLDNdLFszLDRdLFsyLDRdLFswLDMsIihwX3swfSxwX3sxfSkiLDAseyJjdXJ2ZSI6Mn1dLFsxLDQsIiIsMSx7InN0eWxlIjp7Im5hbWUiOiJjb3JuZXIifX1dLFswLDUsIiIsMSx7ImN1cnZlIjotMn1dLFs1LDIsIkZee0pcXHRpbWVzXFxtYXRoYmJ7Mn19Il1d
\[\begin{tikzcd}[ampersand replacement=\&]
	{\Hom_{\A^{J}}(\cDelta,d)} \&\& {\A^{J\times\mathbb{2}}} \\
	\& {\Hom_{\B^{J}}(\cDelta,d)} \&\& {\B^{J\times\mathbb{2}}} \\
	\\
	\& {\cDelta\times J} \&\& {\B\times\B.}
	\arrow[dashed, from=1-1, to=2-2]
	\arrow[from=2-2, to=2-4]
	\arrow[from=2-2, to=4-2]
	\arrow[from=4-2, to=4-4]
	\arrow[from=2-4, to=4-4]
	\arrow["{(p_{0},p_{1})}", curve={height=12pt}, from=1-1, to=4-2]
	\arrow["\lrcorner"{anchor=center, pos=0.125}, draw=none, from=2-2, to=4-4]
	\arrow[curve={height=-12pt}, from=1-1, to=1-3]
	\arrow["{F^{J\times\mathbb{2}}}", from=1-3, to=2-4]
\end{tikzcd}\]
\end{notation}

There are some useful behaviours that an $\infty$-functor with respect to the arrows, diagram cones, and (co)limits.

\begin{definition}
[{\cite[ Definition 4.4.2.7.]{kerodon}}] Let $\C$ and $\D$ be $\infty$-categories. We say that an $\infty$-functor $F : \C \to \D$ is \textit{conservative} if it satisfies the following condition:
\begin{itemize}
\item Let $u : X \to Y$ be a morphism in $\C$. If $F(u) : F(X) \to F(Y )$ is an isomorphism in the $\infty$-category $\D$, then $u$ is an isomorphism.
\end{itemize}
\end{definition}

\begin{definition}
\label{defcreates}
For any class of diagrams $K=\{K_{\gamma}: J_{\gamma} \to C\}_{\gamma}$ valued in $\C$, an $\infty$-functor $F : \C \to \D$:
 \begin{itemize}
 \item \textit{preserves} those (co)limits if for any diagram $K_{\gamma}: J_{\gamma} \to C$ and (co)limit cone (over) under $K_{\gamma}$, the image of this cone defines a (co)limit cone (over) under the composite diagram $FK_{\gamma}: J_{\gamma} \to D$;
\item \textit{reflects} those (co)limits if any cone (over) under a diagram $K_{\gamma}: J_{\gamma} \to \C$, whose image upon applying $F$ is a (co)limit cone for the diagram $FK_{\gamma}: J_{\gamma} \to \D$, is a (co)limit cone (over) under $K_{\gamma}$; 
\item \textit{creates} those (co)limits if whenever $FK_{\gamma}: J_{\gamma} \to \D$ has a (co)limit in $\D$, there is some (co)limit cone (over) under $FK_{\gamma}$ that can be lifted to a (co)limit cone (over) under $K_{\gamma}$, and moreover $F$ reflects the (co)limits in the class of diagrams.
\end{itemize}
\end{definition}

As previously mentioned, we will use the Lurie-Barr-Beck Theorem in the proof of \Cref{thCatEnInLmodAreLmod}. Therefore, we need to prove that an $\infty$-functor $U_{!}$ creates colimits of $U_{!}$-split simplicial sets. Despite this we do not define $U_{!}$-split simplicial sets explicitly. We avoid using them in the proof thanks to the following fact: every colimit of $U_{!}$-split simplicial sets is a sifted colimit (\cite[Lemma 5.5.8.4]{HTT}). We will prove that our $\infty$-functor $U_{!}$ preserves sifted colimits.

For a detailed presentation and definition of $U_{!}$-split simplicial objects, see \cite[\textsection 4.7.2]{HA}, or refer to \cite[Definition 2.3.13, Proposition 2.3.15]{RiehlElements} for a more categorical exposition.

We define the important concept of sifted simplicial sets.

\begin{definition}
A simplicial set $K$ is \textit{sifted} if it satisfies the following conditions:
\begin{itemize}
    \item [(1)] $K$ is nonempty;
    \item[(2)] The diagonal map $K\to  K \times K$ is a initial $\infty$-functor.
\end{itemize}
\end{definition}

\begin{notation}
We say that a diagram $D: A\to \C$ is \textit{sifted} if its source $A$ is a sifted simplicial set. We say that a colimit of a diagram $D: A\to \C$ is \textit{sifted} if the diagram $D$ is sifted.
\end{notation}

\begin{remark}
The idea is that sifted colimits are the colimits preserved by the forgetful functor from the category of groups $U: Gr\to Set$.
Sifted colimits play a crucial role in the study of additive $\infty$-categories see \cite[\S C.1.5]{SAG}. %which we hope will become a topic of future work.
\end{remark}

%%%%%%%%%%%%%%%%%%%%%%%%%%%%%%%%%%%%%%%%%%%%%%%%%%%%%%%%%%%%%%%%%%%%%%%%%%%%%%%%%%%%%%%%%%%%%%%%%%%%%%%%%%%%%%%%%%%%%%%%%%%%%%%%%%%%%%%%%%%%%%%%%%%%%%%%%%%%%%%%%%%%%%%%%%%%%%%%%%%%%%%%%%%%%%%%%%%%%%%%%%%%%%%%%%%%%%%%%%%%%%%

%% file: Chapters/RecallHigherAlgebra.tex
\section{Higher algebraic structure}
\label{subsectionHigherAlgebraicStructure}

In this section, we recall some Higher Algebraic structures. This presentation is not exhaustive, and we recommend interested readers to consult \cite[\S 4.2]{HA}, \cite{heine2023equivalence}, or \cite[\S 2.7]{HinichYon} for a comprehensive treatment.

In \cite{DoniCategorical}, we provide categorical algebraic definitions of these structures. Beginners may find it useful to read before it. While the structures are conceptually the same, the style of the definitions differs significantly. In ordinary algebraic category theory, definitions consist of mathematical objects with a list of axioms. In the higher world, such definitions are not feasible due to the endless list of axioms.

In higher category theory, the only viable method to define higher algebraic structures is using $\infty$-operad theory (also known as categorical pattern, fibrous theory). Initially, we intended to include a section on this theory in the paper, but unfortunately, we did not realize this plan. For an introduction to this topic, see \cite[\S 2]{HA}, \cite{HinichYon}, or \cite{GepHauEnriched}. Note that $\infty$-operad theory is the $\infty$-categorical analog of fc-multicategories theory or colored operad theory.

Five ordinary categories serve as the \enquote{shapes} of the higher categorical algebraic structure.

\begin{definition}
\label{defBMLMASS}
We define five categories: $Ass=\Delta^{op}$, $BM:=(\Delta_{/[1]})^{op}$, $LM\subset BM$ the full subcategory spanned by $\sigma : [n] \to [1]$ satisfying $\sigma(0) = 0$ and having at most one value equal to $1$, and $RM:= LM^{op}$. Let $Fin_{*}$ denote the category of finite pointed sets; we define $Comm:=Fin_{*}$.
\end{definition}

\begin{remark}[{Technical point}]
The categories $LM,RM,BM,Ass$ in \Cref{defBMLMASS} are defined in \cite{HinichYon} and they are not $\infty$-operads. They are a strong approximation of $\mathcal{A}ss^{\otimes}$, $\mathcal{LM}^{\otimes}$, $\mathcal{LM}^{\otimes}$, or $\mathcal{BM}^{\otimes}$ defined in \cite{HA}, as per \cite[\S 2.9.3]{HinichYon}. The existence of a strong approximation between an $\infty$-category and an $\infty$-operad, implying that the $\infty$-category \enquote{shapes the same higher algebraic structure as the $\infty$-operad} (\cite[Proposition 2.7.2]{HinichYon}). In this article, we do not distinguish between an $\infty$-operad and its strong approximation (e.g., in \Cref{defMonoidalCategory}, $\mathcal{A}ss^{\otimes}$, $\mathcal{LM}^{\otimes}$, $\mathcal{LM}^{\otimes}$, or $\mathcal{BM}^{\otimes}$ should be replaced with $Ass$, $LM$, $BM$, and $RM$, respectively). $Comm$ is defined in \cite[Notation 2.0.0.2]{HA}.
\end{remark}

\begin{notation}
Sometimes we use the notation $\E$ and $\mathbb{E}_{1}$ to denote $Comm$ and $Ass$ respectively. This is abuse because they are not the same mathematical objects but they \enquote{shapes the same higher algebraic structure as the $\infty$-operad}. For a deep study of $\E$ and $\mathbb{E}_{1}$ see \cite[\S 5]{HA}.
\end{notation}

$Comm$ and $Ass$ have only one color $\{a\}$, $BM$ has three colors $\{a,m,b\}$, $LM$ has two colors $\{a,m\}$, and $RM$ has two colors $\{m,b\}$.

The idea is that $Ass$ shapes monoidal $\infty$-categories, $Comm$ shapes symmetric monoidal $\infty$-category, $BM$ shapes bimodules objects, $LM$ left-module objects and $RM$ right-module objects.

Let us begin with the definition of a monoidal $\infty$-category. We provide a more general definition.

\begin{definition}[{\cite[Definition 2.1.2.13]{HA}}]
\label{defMonoidalCategory}
Let $O\in\{Comm, Ass,BM,RM,LM\}$ be a category. We say that an $\infty$-functor $p : \C^{\otimes}\to O$ is a cocartesian fibration of $\infty$-operads if it satisfies the hypotheses of \cite[Proposition 2.1.2.12]{HA}. In this case, we also say that $p$ exhibits $\C^{\otimes}$ as an $O$-monoidal $\infty$-category.
\end{definition}

\begin{example}[Definition of monoidal $\infty$-category]
\label{ExMonoidalHigherCategory}
A monoidal $\infty$-category is an $Ass$-monoidal $\infty$-category $\mathcal{F}:\W^{\otimes}\to Ass=\Delta^{op}$.

A symmetric monoidal $\infty$-category is a $Comm$-monoidal $\infty$-category $\mathcal{F}:\W^{\otimes}\to Comm$.
\end{example}

Roughly speaking, an $Ass$($Comm$)-monoidal $\infty$-category $\W^{\otimes}\to Ass=\Delta^{op}$ defines an $\infty$-category $\W$ and an $\infty$-functor $-\otimes_{W}-:\W\times \W\to \W$ which plays the role of a (commutative) operation with a unit $\I_{\W}$. As is customary, we will make an abuse.
The symbol $\otimes$ is typically employed to emphasize that $\mathcal{W}^{\otimes}$ is an $\infty$-category, along with additional higher algebraic structures. When the context is clear, we may simply use $\mathcal{W}$ or $(\mathcal{W}, \otimes_{\mathcal{W}}, \I_{\mathcal{W}})$ to denote a (symmetric) monoidal category.

\begin{example}
    The triple ($\Cat,\times, \mathbb{1}$) whose second component is the product and the third component is the terminal $\infty$-category forms a symmetric monoidal $\infty$-category.
\end{example}

\begin{example}[Definition of $\V$-left tensored $\infty$-category]
\label{exVLeftTensoredCategory}
Let $\V$ be a monoidal $\infty$-category. A $\V$-left(-right) tensored $\infty$-category is an $LM$($RM$)-monoidal $\infty$-functor $\mathcal{F}:\W^{\otimes}\to LM$ such that there is a monoidal equivalence of monoidal $\infty$-categories $\V\simeq \W\times_{RM}\{a\}$ (see \Cref{defMorphismLaxMonoidal}).
\end{example}

%%%%%%%%%%%%%%%%%%%%%%%%%%%%%%%%%%%%%%%%%%%%%%%%%%%%%%%%%%%%%%%%%%%%%%%%%%%%%%%%%%%

Roughly speaking, a $LM$-monoidal $\infty$-category $\W^{\otimes}\to LM$ defines an $\infty$-category $\W$ and an $\infty$-functor $-\;^{\V}\otimes_{\W}-:\V\times \W\to \W$ which respect the higher categorical version of the axioms in \cite[Definition 2.5]{DoniCategorical} and which plays the role of an action on $\W$. As is customary, we will make abuse and we will write only $\W$ or $(\W,^{\V}\otimes_{\W},\I_{\W})$ to denote a $LM$ $(RM)$-monoidal category.

An $RM$-monoidal $\infty$-category is similar to above only that the operation is $-\;^{\V}\otimes_{W}-:\W\times \V\to \W$.

\begin{example}[Definition of $\V$-$\A$ bitensored $\infty$-category]
\label{exVABITensoredHigherCategory}
Let $\V$ and $\A$ be two monoidal $\infty$-categories, a $\V$-$\A$ bitensored $\infty$-category is $BM$-monoidal $\infty$-functor $\mathcal{F}:\W^{\otimes}\to BM$ such that there are two monoidal equivalence of monoidal $\infty$-categories $\V\simeq \W\times_{LM}\{a\}$ and $\A\simeq \W\times_{LM}\{b\}$ (see \Cref{defMorphismLaxMonoidal}).
\end{example}

Roughly speaking, a $BM$-monoidal $\infty$-category $\W^{\otimes}\to BM$ defines an $\infty$-category $\W$ and two $\infty$-functors $-\;^{\V}\otimes_{W}-:\V\times \W\to \W$ and $-\;^{\A}\otimes_{W}-:\W\times \A\to \W$ which respect the same axioms of a $LM$-monoidal $\infty$-category and $RM$-monoidal $\infty$-category plus same axioms that permit to intertwine the right and the left action on $\W$, see \cite[\S 4.3.3]{HA}. As is customary, we will make an abuse and we will write only $\W$ or $(\V,\W,\A)$ to denote a $BM$-monoidal category.

Now, we define how to compare two $O$-monoidal $\infty$-categories. The absence of an introduction to $\infty$-operads theory makes the fundamental definition of lax monoidal and monoidal $\infty$-functor not completely clear. For the sake of completeness, we write it anyway.

\begin{definition}
\label{defMorphismLaxMonoidal}
Let $\W$ and $\V$ be two (symmetric) monoidal $\infty$-categories.
We call \textit{lax monoidal $\infty$-functor} with source $\W^{\otimes}$ and target $\V^{\otimes}$ an element $F:\W^{\otimes}\to \V^{\otimes}$ in $\Hom_{Cat_{\infty/Ass}}(\W^{\otimes},\V^{\otimes})$ such that
\begin{itemize}
    \item[(*)] $F$ carries inert morphisms in $\W^{\otimes}$ to inert morphisms in $\V^{\otimes}$.
\end{itemize}

We call \textit{monoidal $\infty$-functor} with source $\W^{\otimes}$ and target $\V^{\otimes}$ an element in $F:\W^{\otimes}\to \V^{\otimes}$ in $\Hom_{Cat_{\infty/Ass}}(\W^{\otimes},\V^{\otimes})$ such that
\begin{itemize}
    \item[(**)] $F$ carries cocartesian morphisms in $\W$ to cocartesian morphisms in $\V$.
\end{itemize}
\end{definition}

A lax monoidal $\infty$-functor $F:\W^{\otimes}\to \V^{\otimes}$ is an $\infty$-functor $f:\W\to\V$, called \textit{underlying $\infty$-functor of $F$}, plus, for each pair of objects $x,y\in\W$, a morphism in $\V$
\[\theta_{x,y}:f(x)\otimes_{\V} f(y)\to f(x\otimes_{\W}y)\]
and a morphism between the unit objects
\[\mu:\I_{\V}\to f(\I_{W}).\]
Instead, an $\infty$-functor $F:\W^{\otimes}\to\V^{\otimes}$ is monoidal if $\theta_{x,y}$ and $\mu$ are isomorphism.

In this article, we make the abuse to denote a (lax) monoidal $\infty$-functor by its underlying $\infty$-functor.

\begin{remark}
\label{rmkINERTMORPHISMMANNAGGIA}
We have never defined what an inert morphism is, and we will not. We only spend a few words about them in this remark. In $\infty$-operads theory, there are two important types of morphisms: inert and active. 
The collection of cocartesian morphism includes both.
The idea is that inert morphisms \enquote{create the structure}, so the property $(*)$ means that $s$ \enquote{respects the structure}. Instead active morphisms \enquote{define the operation} of the (symmetric) monoidal $\infty$-category, so the property $(**)$ means that $s$ also preserves \enquote{the operation}.
The concept of inert arrow also appears in \Cref{notUltimaCredodellarticolo}.
\end{remark}

%%%%%%%%%%%%%%%%%%%%%%%%%%%%%%%%%%%%%%%%%%%%%%%%%%%%%%%%%%%%%%%%%%%%%%%%%%%%%%%%%%%%%%%%%%%%%%%%%%%%%%%%%%%%%%%%%%%%%%%%%%%%%%%

Until now, the higher algebraic mathematical structures are composed of $\infty$-categories and $\infty$-functors, but in this series of articles, we need more general algebraic structures and we need to add properties to these $\infty$-categories (e.g., presentable $\infty$-categories) and these $\infty$-functors (e.g., preserving-small-colimits $\infty$-functors).

To obtain this, we need the internal module objects theory in an $O$-monoidal $\infty$-category for $O\in\{Comm, Ass, BM, RM, LM\}$. The absence of an introduction to $\infty$-operads theory makes the fundamental definition of algebra in a $O$-monoidal $\infty$-category (\Cref{defAlgebra}) not completely clear. For the sake of completeness, we write it anyway.

\begin{definition}[{\cite[ Definition 2.1.2.7]{HA}}]
\label{defAlgebra}
Let $O\in\{Comm, Ass, BM, RM, LM\}$ be a category and $\W$ be an $O$-monoidal. An \textit{$O$-algebra objects (or, $O$-algebra)} is a global section $s:O\to \W$ with the property:
\begin{itemize}
    \item[(*)] that $s$ carries inert morphisms in $O$ to inert morphisms in $\W$.
\end{itemize}
By $Alg_{O}(\W)\subseteq Fun_{Cat_{\infty/O}}(O,\W)$ we denote the full $\infty$-subcategory spanned of $O$-algebras and we call it the \textit{$\infty$-category of $O$-algebras}. 
\end{definition}

\begin{example}
\label{defAlgebraHigherAlgebra}
Let $\V\to Ass(Comm)$ be a (symmetric) monoidal $\infty$-category. We let $Alg(\V)$ $(CAlg(\V))$ denote the $\infty$-category of (commutative and) associative algebra objects in $\V$, i.e. $Alg_{Ass}(\V)$ $(Alg_{Comm}(\V))$.
\end{example}

Roughly speaking, a (commutative and) associative algebra object in $\V$ $s: Ass (Comm)\to \V$ chooses an object $v\in\V$ and a morphism $m_v:v\otimes_{\V}v\to v$ in $\V$ which plays the role of a (commutative and) associative operation with a unit which satisfies the categorical version of the axioms in \cite[Definition 2.6]{DoniCategorical}. As customary, we will make abuse and write only $v$ or $(v,m_v)$ to denote a (commutative and) associative algebra object in $\V$.

\begin{notation}[{\cite[Definition 2.3.1]{GepHauEnriched}}]
\label{notUltimaCredodellarticolo} Let $(\Delta^{op})^{int}$ denote the wide subcategory of $\Delta^{op}\simeq Ass$ where the morphisms are the inert morphisms in $\Delta^{op}$.

Let $\V$ and $\W$ two monoidal $\infty$-category and let $F$ be an $\infty$-functor between $\infty$-categories of algebra objects:
\begin{equation*}
\label{eqFTilda}
F:Alg_{Ass}(\V)\to Alg_{Ass}(\W).  
\end{equation*}
The precomposition with the inclusion $(\Delta^{op})^{int}\to Ass$ induces an $\infty$-functor 
\begin{equation}
\label{eqCiao}
F^{triv} := \V \simeq Alg_{(\Delta^{op})^{int}}(\V)\to Alg_{(\Delta^{op})^{int}}(\W) \simeq \W;
\end{equation}
where the two equivalences above are proved in \cite[\S 3.4]{GepHauEnriched}. We denote by  $F^{triv}$ the $\infty$-functor in \eqref{eqCiao}. 

\end{notation}
\begin{remark}
   In the definition of $(\Delta^{op})^{int}$, we encounter a similar issue as mentioned in Remark \ref{rmkINERTMORPHISMMANNAGGIA}. However, this issue does not pose a problem for us. The crucial aspect of this notation lies in understanding that from a given $\infty$-functor between the $\infty$-category of algebra objects \eqref{eqFTilda}, we are able to derive an $\infty$-functor between $\infty$-categories, as depicted in \eqref{eqCiao}.
\end{remark}

\begin{example}
\label{defLModHigherAlgebra}
 Let $\V\to Ass$ be a monoidal $\infty$-category and let $q : \mathcal{M}^{\otimes}\to LM$ $(RM)$ be a $\V$-left (right) tensored $\infty$-category such that $\mathcal{M}^{\otimes}\times_{LM}m\simeq \V$ $(\mathcal{M}^{\otimes}\times_{RM}m\simeq \V)$. We let $\mathrm{LMod}(\mathcal{M})$ $(\mathrm{RMod}(\mathcal{M}))$ denote the $\infty$-category $Alg_{LM}(\mathcal{M}^{\otimes})$ $(Alg_{RM}(\mathcal{M}^{\otimes}))$. We will refer to $\mathrm{LMod}(\mathcal{M})$ as \textit{the $\infty$-category of $\V$-left (right) module objects of $\mathcal{M}$}.
 Composition with the inclusion $Ass\to LM$ determines a categorical fibration
 
\[\mathrm{LMod}(\mathcal{M}) = Alg_{LM}(\mathcal{M}^{\otimes}) \to Alg(\V)= Alg_{Ass}(\V^{\otimes}).\]
\[(\mathrm{RMod}(\mathcal{M}) = Alg_{RM}(\mathcal{M}^{\otimes}) \to Alg(\V)= Alg_{Ass}(\V^{\otimes})).\]

If $A$ is an algebra object of $\V$, we let $\mathrm{LMod}_A(\mathcal{M})$ ($\mathrm{RMod}_A(\mathcal{M})$) denote the fiber $\mathrm{LMod}(\mathcal{M}) \times_{Alg(\V)} \{A\}$ ($\mathrm{RMod}(\mathcal{M}) \times_{Alg(\V)} \{A\}$); we will refer to $\mathrm{LMod}_{A}(\mathcal{M})$ ($\mathrm{RMod}_A(\mathcal{M})$) as the \textit{$\infty$-category of left (right) $A$-module objects of $\mathcal{M}$}.
\end{example}

Roughly speaking, a left $A$-module object of $\mathcal{M}$ $s: LM\to \mathcal{M}^{\otimes}$ chooses $m\in \mathcal{M}$ an object of $\mathcal{M}$ and $l: A\;^{\V}\otimes_{\mathcal{M}}m\to m$ a morphism of $\mathcal{M}$ which respect the higher categorical version of the axioms in \cite[Definizione 2.8]{DoniCategorical} and which plays the role of an action on $m$. As customary, we will make abuse and write only $m$, $(A,m)$ or $(m,l)$ to denote a left $A$-module objects of $\mathcal{M}$.

A right $A$-module object of $\mathcal{M}$ has a similar description.

\begin{example}
\label{defBOMDHigherAlgebra}
  Let $\V\to Ass$ and $\A\to Ass$ be monoidal $\infty$-categories and let $q : \mathcal{M}^{\otimes}\to BM$ be a $\V$-tensored $\infty$-category such that $\mathcal{M}^{\otimes}\times_{BM}m\simeq \mathcal{M}$. We let $\mathrm{BMod}(\mathcal{M})$ denote the $\infty$-category $Alg_{BM}(\mathcal{M})$. We will refer to $\mathrm{BMod}(\mathcal{M})$ as the $\infty$-category of $\V$-$\A$- bimodule objects of $\mathcal{M}$.

 Composition with the inclusion $Ass\times Ass\to BM$ determines a categorical fibration
 
\[\mathrm{BMod}(\mathcal{M}) = Alg_{BM}(\mathcal{M}^{\otimes}) \to Alg(\V)\times Alg(\A)= Alg_{Ass}(\V^{\otimes})\times Alg_{Ass}(\A^{\otimes}).\]

If $C$ is an algebra object of $\V$ and $B$ is an algebra object of $\A$, we let $_{C}\mathrm{Mod}_B(\mathcal{M})$ or $_{C}\mathrm{BMod}_{B}(\mathcal{M})$ denote the fiber $\mathrm{BMod}(\mathcal{M}) \times_{Alg(\V)\times Alg(\A)} \{(C,B)\}$; we will refer to $_{C}\mathrm{Mod}_B(\mathcal{M})$ as the $\infty$-category of \textit{$C$-$B$-bimodule objects of $\mathcal{M}$}.
\end{example}

Roughly speaking, a $C$-$B$-bimodule object of $\mathcal{M}$ denoted by $s:BM\to \mathcal{M}^{\otimes}$ selects $m\in \mathcal{M}$ an object of $\mathcal{M}$ and two morphisms in $\mathcal{M}$: $l:C\;^{\V}\otimes_{\mathcal{M}}m\to m$ and $r:m\;^{\A}\otimes_{\mathcal{M}}B\to m$. These morphisms respect the same axioms as a left $C$-module object and right $B$-module object, along with axioms that permit the intertwining of the actions $r$ and $l$ on $m$, see \cite[\S 4.3.3]{HA}. As customary, we use the notation $(C,m,A)$ or $m$ to denote a $C$-$B$-bimodule object of $\mathcal{M}$.
Furthermore, if we want to emphasize that an element $m$ of $_C \mathrm{BMod}_B(\mathcal{M})$ has a left action of $C$ and a right action of $B$, we will write the element $m$ with two subscripts $_C m_B$.

\begin{remark}
In \Cref{defAlgebraHigherAlgebra}, \Cref{defLModHigherAlgebra}, and \Cref{defBOMDHigherAlgebra}, if one substitutes the $\infty$-categories $\V$, $\A$, and $\mathcal{M}$ with the $O$-monoidal $\infty$-category $(\Cat,\times, \mathbb{1})$, where $O$ is an element of the set $\{Comm, Ass, BM, RM, LM\}$, after some computation, one finds the \Cref{ExMonoidalHigherCategory}, \Cref{exVLeftTensoredCategory}, and \Cref{exVABITensoredHigherCategory}, respectively.
\end{remark}

%% file: Chapters/RecallPresentableStable.tex
\section{Stable $\infty$-categories and presentable $\infty$-categories}
\label{subsecStablepresentable}
\subsection{Stable and presentable $\infty$-category}
In this paper the notion of stability is marginal but in \cite{DoniklinearMorita}, it is a fundamental notion. So we choose the definition of stable $\infty$-category the most similar to the $dg$-categorical notion of pretriangulated $dg$-category.  

\begin{definition}[{\cite[Definitions 1.1.1.9, Corollary 1.4.2.27]{HA}}]
\label{defStable}
An $\infty$-category $\C$ is stable if it satisfies the following conditions:
\begin{itemize}
    \item There exists a zero object $\mathbf{0}$;
    \item Every morphism in $\C$ admits cofiber;
    \item The suspension $\infty$-functor
    \begin{equation}
    \label{eqsuspension}    
    \Sigma:\C\to \C:X\mapsto \operatorname{pushout}(\mathbf{0}\leftarrow X\to \mathbf{0})
    \end{equation}
    is an equivalence of $\infty$-categories, which inverse is the loop $\infty$-functor
    \begin{equation}
    \label{eqloop}
          \Omega:\C\to\C:X\mapsto \operatorname{pullback}(\mathbf{0}\to X \leftarrow \mathbf{0}).      
    \end{equation}

\end{itemize}
\end{definition}

In \cite{HA}, Lurie has proved that the above definition is equivalent to the following one; this second one allows defining the $\infty$-category of (small) stable $\infty$-categories.

\begin{definition}[{\cite[Definitions 1.1.1.9, Proposition 1.1.3.4]{HA}}]
\label{defPref2}
An $\infty$-category $\C$ is stable if it is pointed (i.e., it has zero objects $\mathbf{0}$), it admits finite limits, it admits finite colimits, and each square, $\Delta^1\times\Delta^1\to \C$, in $\C$ 
is a pushout if and only if it is a pullback.
\end{definition}

%\begin{notation}
%We will denote by $\Cat^{Ex}\subseteq \Cat$ the $\infty$-subcategory whose objects are the stable $\infty$-categories and whose functors are exact, i.e., they preserve finite limits and finite colimits.
%\end{notation}

\begin{remark}
Directly from the definition follows that the $\infty$-category $A$ is stable if and only if its dual category $A^{op}$ is stable.
\end{remark}

Thanks to the work of Groth, Lurie \cite{HA}, there are many equivalent definitions of stable $\infty$-category, see \cite[theorem 4.4.12]{RiehlElements}.

Before listing alternative definitions, we recall some facts.

If $\C$ is an $\infty$-category with terminal elements $t:\mathbb{1}\to \C$, we will call the comma $\infty$-category $\Hom_{\C}(t,\C)$ the \textit{$\infty$-category of pointed elements of $\C$} and we denote it by $\C_{*}$. This $\infty$-category is pointed and it is the canonical way to pointed an $\infty$-category with terminal elements, \cite[Lemma 4.4.2. ]{RiehlElements}.

Moreover, if $\C$ is also finite complete, there exists an adjunction
\begin{equation}
\label{eqAdPointedForgetful}
\begin{tikzcd}[ampersand replacement=\&]
	\C \& {\C_{*};}
	\arrow[""{name=0, anchor=center, inner sep=0}, "{-\coprod t}", shift left=2, from=1-1, to=1-2]
	\arrow[""{name=1, anchor=center, inner sep=0}, "U", shift left=2, from=1-2, to=1-1]
	\arrow["\dashv"{anchor=center, rotate=-90}, draw=none, from=0, to=1]
\end{tikzcd}
\end{equation}
where the right adjoint is the cocartesian projection of the comma $\infty$-category $\C_{*}$ and the left adjoint is the $\infty$-functor: \[-\coprod t:a\mapsto \{(a\to a\coprod t,t\to a\coprod t)\}.\]    

Indeed, using the dual result of \cite[Corollary 12.2.7]{RiehlElements}, to find the left adjoint of $U$  it is sufficient to find for each element of $\C$ $a:\mathbb{1}\to \C$ the initial element of the comma $\infty$-category $\Hom_{\C_{*}}(a, U)$. It is easy to see that the initial element is the pair composed by the legs of the coprod itself with the terminal elements, $-\coprod t:a\mapsto \{(a\to a\coprod t,t\to a\coprod t)\}$; it is possible to prove this last sentence using the universal property of coproduct.

It is obvious that if $\C$ is a pointed $\infty$-category, then the above adjunction is an equivalence of $\infty$-categories.

  Now, we recall the suspension-loop adjunction in the case that $\C$ is not necessarily stable.
 Let $\C$ be a pointed $\infty$-category with finite limits and finite colimits, then the $\infty$-functors \Cref{eqloop} and \Cref{eqsuspension} define an adjunction:
 \begin{equation}
 \label{eqAdSigmaOmega}
 \begin{tikzcd}[ampersand replacement=\&]
	\C \& {\C.}
	\arrow[""{name=0, anchor=center, inner sep=0}, "\Sigma", shift left=2, from=1-1, to=1-2]
	\arrow[""{name=1, anchor=center, inner sep=0}, "\Omega", shift left=2, from=1-2, to=1-1]
	\arrow["\dashv"{anchor=center, rotate=-90}, draw=none, from=0, to=1]
\end{tikzcd}
 \end{equation}
 
If $\C$ is not pointed but only finite complete and finite cocomplete, the loop-suspension adjunction is defined as the composition of two adjunctions (\ref{eqAdSigmaOmega}) and (\ref{eqAdPointedForgetful}): 
\begin{equation}
 \label{eqAdSigmaOmegaplus}
 \Sigma:=\Sigma(-\coprod t):\C\leftrightarrows\C_{*}\leftrightarrows \C_{*}:U\Omega=:\Omega.
 \end{equation}

 With a small abuse of notation, we denote the $\infty$-functors in the adjunction (\ref{eqAdSigmaOmegaplus}) by $\Sigma$ and $\Omega$, just as in (\ref{eqAdSigmaOmega}).      
\\ \\
We need to recall another tool from \cite{HA} before giving other definitions.

\begin{definition}    
Let $F:\C\to \D$ be an $\infty$-functor between $\infty$-categories,
\begin{itemize}
    \item If $\C$ admits pushouts, then we will say that $F$ is \textit{excisive} if $F$ carries pushout square in $\C$ to pullback square in $\D$;
    \item if $\C$ admits a terminal element $t:\mathbb{1}\to \C$, we say that $F$ is \textit{reduced} if $F(t)$ is a terminal element of $\D$; 
    \item let $S^{fin}$ be the $\infty$-category generated by a single element under finite colimits, we denote by $Sp(\C)$ the full subcategory of $Fin(S^{fin}_*,\C)$ spanned by reduced, excisive $\infty$-functors. The objects of $Sp(\C)$ are called \textit{$\C$-spectrum objects};
    
    \item with $\Omega^{\infty}:Sp(\C)\to \C:F\mapsto F(\mathcal{S}^{0})$, we denote the evaluation $\infty$-functor in the pointed $0$-sphere.
\end{itemize}
\end{definition}

\begin{example}
\label{exOmegaInftyCambioArricchimento}
In the case that $\C\simeq\mathcal{S}$, $\Sp(\Top)$ is denote by $\Sp$ and it is called \textit{the $\infty$-category of spectra}.
We have a nice description for $\Omega^{\infty}:\Sp\to \mathcal{S}_{*}$ and its left adjoint.
The $\infty$-category $\Sp$ is symmetric monoidal with the famous smash product as a tensor and the sphere spectrum $\mathbb{S}$ as unit. There are many ways to define this symmetric monoidal structure and one of them is as the inherited structure by the following localization:
% https://q.uiver.app/#q=WzAsMixbMCwwLCJGdW4oXFxtYXRoY2Fse1N9XntmaW59XyosXFxtYXRoY2Fse1N9KSJdLFsyLDAsIlxcU3AiXSxbMCwxLCIiLDAseyJvZmZzZXQiOi0yfV0sWzEsMCwiIiwwLHsib2Zmc2V0IjotMywic3R5bGUiOnsidGFpbCI6eyJuYW1lIjoiaG9vayIsInNpZGUiOiJib3R0b20ifX19XV0=
\[\begin{tikzcd}
	{Fun(\mathcal{S}^{fin}_*,\mathcal{S})} && \Sp;
	\arrow["L",shift left=2, from=1-1, to=1-3]
	\arrow["i",shift left=3, hook', from=1-3, to=1-1]
\end{tikzcd}\]
where the functors $\infty$-category $Fun(\mathcal{S}^{fin}_*,\mathcal{S})$ is symmetric monoidal with the Day Convolution structures.

In this model, the sphere spectrum is the image through $L$ of the covariant Yoneda; that is $\mathbb{S}\cong L\Yo^{S^{0}}$.

Moreover, thanks to the Yoneda Lemma, the evaluation $\Omega^{\infty}$ is the corepresented $\infty$-functor \[\Hom_{\Sp}(\mathbb{S},-)\cong\Hom_{\Sp}(L\Yo^{S^{0}},-)\cong \Hom_{Fun(\Top^{fin}_{*},\Top)}(\Yo^{S^{0}}, i(-))\cong ev_{0}(-)= \Omega^{\infty}(-);\] and its left adjoint is the tensor with the sphere spectrum $-\otimes \mathbb{S}$,
\[-\otimes \mathbb{S}:\mathcal{S}\rightleftarrows \Sp:\Hom_{\Sp}(\mathbb{S} ,-)\cong \Omega^{\infty}.\]
We are in the same situation as \Cref{exFreePrecUnderSeg}.
\end{example}

\begin{remark}[{\cite[Remark 1.4.2.18]{HA}}] Let $\C$ be an $\infty$-category which admits finite limits. Then the forgetful functor $U: \C_{*} \to \C$ in \Cref{eqAdPointedForgetful} induces an equivalence of
$\infty$-categories $\Sp(\C_*) \to \Sp(\C)$. 
In particular, this means that our definition of $\Sp$ is equivalent to the definition of Lurie \cite[Definition 1.4.3.1]{HA}.
\end{remark}

%%%%%%%%%%%%%%%%%%%%Fino A Qua RIVISTO%%%%%%%%%%%%%%%%%%

Finally, we list some reformulations of the definition of a stable $\infty$-category:

\begin{itemize}
     
    \item[(i)] let $\C$ be an $\infty$-category with finite limits. Then, $\C$ is stable if and only if the evaluation $\infty$-functor $\Omega^\infty$ is an equivalence;

    \item[(ii)] let $\C$ be an $\infty$-category with finite limits. $\C$ is stable if and only if the following equivalences hold: 
    \[\C\cong \lim \left\{ \dots \C_*\xrightarrow{\Omega} \C_*\xrightarrow{\Omega}\C_*\xrightarrow{\Omega}\dots \right\}\cong \lim \left\{ \dots \xrightarrow{\Omega} \C_* \xrightarrow{\Omega} \C_*. \right\};\]

     \item[(iii)] let $\C$ be an $\infty$-category with finite colimits. $\C$ is stable if and only if the following equivalences hold: 
     \[\C\cong \operatorname{colim} \left\{ \dots \C_*\xleftarrow{\Sigma} \C_*\xleftarrow{\Sigma}\C_*\xleftarrow{\Sigma}\dots \right\}\cong \operatorname{colim} \left\{ \dots \xleftarrow{\Sigma} \C_* \xleftarrow{\Sigma} \C_* \xleftarrow{\Sigma} \C_*. \right\};\]

     \item[(iv)] $\C$ is stable if $\C$ is a pointed $\infty$-category and there exists an adjunction equivalence \[\Sigma:\C\rightleftarrows\C:\Omega;\] where $\Sigma$ is the suspension $\infty$-functor, and $\Omega$ is the loop $\infty$-functor.

    \end{itemize}

In \cite[\S 1.4]{DoniPhDThesis}, we  explain why $(ii)$ and $(iii)$ are dual.

In particular, $(ii)$ allow us to represent an element of a stable $\infty$-category $\C$ as:

\begin{itemize}
    
    \item a sequence of elements in $\C$ $(X_{n})_{n\in\Z}$ such that $\Omega X_{n+1}\cong X_{n}$. Equivalently, by adjoint equivalence $\Sigma\dashv \Omega$, as a sequence $ X_{n+1}\cong \Sigma X_{n}$;  

    \item a numeral sequence of elements in $\C$ $(X_{n})_{n\in\N}$ such that $\Omega X_{n+1}\cong X_{n}$. Equivalently, by the adjoint equivalence $\Sigma\dashv \Omega$, as a sequence $ X_{n+1}\cong \Sigma X_{n}$. 
   
\end{itemize} 

We might have used the dual description, i.e. $(iii)$, which would have led to a dual representation, but we have chosen $(ii)$ to keep the same convention as in \cite{SchwedeSymSp} and \cite{HA}, see in \Cref{ExOurCase}.

%%%%%%%%%%%%%% Fino A QUA Rivisto %%%%%%%%%%%%%%%%%%%%

Having multiple definitions is very helpful (although it may be confusing at first) because it allows for different points of view.

\begin{notation}
Let $\C$ be an $\infty$-category. By $h\C$, we denote the ordinary category whose objects are the same as those of $\C$ and whose hom-set between two objects $x, y$ is the set $\pi_0(\C(x, y))$. We call $h\C$ the \textit{homotopy category of $\C$}.
For a comprehensive presentation of $h\C$, we recommend reading \cite[\S 1.1]{RiehlElements}.
\end{notation}

\begin{notation}[{\cite[\textsection 1.1]{HA}, \cite[\textsection 4.4]{RiehlElements}}]
\label{rmkhCtriang}
 In a stable $\infty$-category $\C$, the homotopy category $h\C$ is a triangulated category. The idea is that an exact triangle is a cofiber (which inside a stable $\infty$-category is equivalent to being a fiber) diagram in $\C$:
\[
\begin{tikzcd}
	X:=\operatorname{fiber}(Y\to Z) & Y \\
	{\mathbf{0}} & Z:=\operatorname{cofiber}{X\to Y}.
	\arrow[from=1-1, to=2-1]
	\arrow[from=2-1, to=2-2]
	\arrow[from=1-2, to=2-2]
	\arrow[from=1-1, to=1-2]
\end{tikzcd}
\]
Let $\C$ be a stable $\infty$-category. In order to maintain consistent notation with the triangulated category $h\C$, we denote, as everyone does,
by $[n]:\C\to\C$ the $\infty$-functor that sends $X\in\C$ to $\Sigma^{n}X$ if $n\geq 0$ and to  $\Omega^{n}X$ if $n\leq 0$.
%Furthermore, it follows from the adjunction equivalence $\Sigma\dashv \Omega$ that a square in $\C$ is a pullback if and only if it is a pushout.   
\end{notation}

\begin{notation}    
\label{notationEinfinito}
In \Cref{exOmegaInftyCambioArricchimento}, we recalled that $(\Sp,\otimes_{\Sp},\mathbb{S})$ forms a monoidal $\infty$-category.
If $\V=\Sp$ in \Cref{defAlgebraHigherAlgebra}, as customary, we use the convention that we write $Alg$ $(CAlg)$ instead of $Alg(\Sp)$ ($CAlg(\Sp)$). 
An \textit{$\E$-ring} is an object of $CAlg$ and an $\mathbb{E}_1$-ring is an object of $Alg$.
If $\W=\Sp$ in \Cref{defLModHigherAlgebra}, as customary, we use the convention that sometimes we write $\mathrm{LMod}_{A}$ ($\mathrm{RMod}_{A}$) instead of $\mathrm{LMod}_{A}(\Sp)$ ($\mathrm{RMod}_{A}(\Sp)$). 
If $\mathcal{M}=\Sp$ in \Cref{defBOMDHigherAlgebra}, as customary, we use the convention that we write $\mathrm{Mod}$ instead of $\mathrm{BMod}(\Sp)$.

For each $\E$-ring $A$, there is a chain of equivalences:
\begin{equation}
\label{eqLModBModRMod}
_\mathbb{S}\mathrm{Mod}_A \simeq \mathrm{LMod}_A \simeq \; _A \mathrm{Mod}_A \simeq \mathrm{RMod}_A,
\end{equation}
see \cite[\S 4.3.2]{HA}.

\end{notation}

\begin{example}[{\cite[Example 6.27]{SchwedeSymSp}, \cite[\textsection 0]{SS}}]
\label{ExOurCase}
Let $k$ be a commutative unitary ring.
The \textit{Eilenberg-Maclane spectrum} $\Hk=({K(k,n)}_{n\in\Z})\in \Sp$ is the spectrum whose $n$-component is the \textit{Eilenberg-Maclane space} of type $n$. It is an $\mathbb{E}_\infty$-$ring$. 
  See \Cref{exEilenbergMaclaneSpectrumDiscrete} for more about Eilenberg-Maclane spectra.
\end{example}

\begin{example}
\label{exStable}
The $\infty$-category $\Sp\simeq Sp(\mathcal{S})= Ext_{*}(\mathcal{S}_{*},\mathcal{S})$ is stable by construction.  
\end{example}

We need to know what a $t$-structure on a stable $\infty$-category is. Now we will recall its definition and explain some useful examples for us; for a complete introduction, see \cite[\textsection 1.2.1]{HA}.
Roughly speaking, giving a $t$-structure to a stable $\infty$-category $\C$ means adding additional information to $\C$ that allows $\C$ to be broken up into homologous parts, and by studying these parts it is possible to derive information on the whole of $\C$.

\begin{definition}[{\cite[Definition 1.2.1.1 Definition 1.2.1.4 ]{HA}}]
    Let $\C$ be a stable $\infty$-category. A triple $(\C,\C_{\geq 0},\C_{\leq 0})$ where $\C_{\geq 0}$ and $\C_{\leq 0}$ are two subcategories of $\C$ is a \textit{$t$-structure} for $\C$ if the triple $(h\C,h\C_{\geq 0},h\C_{\leq 0})$ is a $t$-structure for the triangulated category $h\C$.
  
\end{definition}

\begin{notation}
Let $(\C,\C_{\geq 0},\C_{\leq 0})$ be a stable $\infty$-category equipped with a $t$-structure, let $\C_{\geq n}$ and $\C_{\leq n}$ denote the full subcategories of $\C$ spanned by those objects that belong to $h\C_{\geq n}$ and $h\C_{\leq n}$, respectively.
\end{notation}

 Given a stable $\infty$-category with a $t$-structure $\C$, for each $n\in \Z$, there exist two adjunctions 
 \[\begin{tikzcd}[ampersand replacement=\&]
	\C \& {\C_{\leq n},}
	\arrow[""{name=0, anchor=center, inner sep=0}, "{\tau_{\leq n}}", shift left=2, from=1-1, to=1-2]
	\arrow[""{name=1, anchor=center, inner sep=0}, "i", shift left=2, hook', from=1-2, to=1-1]
	\arrow["\dashv"{anchor=center, rotate=-90}, draw=none, from=0, to=1]
\end{tikzcd}\]
\[\begin{tikzcd}[ampersand replacement=\&]
	{\C_{\geq n}} \& {\C.}
	\arrow[""{name=0, anchor=center, inner sep=0}, "{\tau_{\geq n}}", shift left=2, from=1-2, to=1-1]
	\arrow[""{name=1, anchor=center, inner sep=0}, "i", shift left=2, hook', from=1-1, to=1-2]
	\arrow["\dashv"{anchor=center, rotate=-90}, draw=none, from=1, to=0]
\end{tikzcd}\]
The former is a localization, and the latter is a colocalization.

\begin{definition}[{\cite[Definition 1.2.1.11]{HA}}]
    Let $\C$ be a stable $\infty$-category equipped with a $t$-structure. \textit{The heart} $\C^{\heartsuit}$ of $\C$ is the full subcategory $\C_{\leq 0}\cap\C_{\geq 0}\subseteq\C $. For each $n\in\Z$, we let $\pi_{0}:\C\to \C^{\heartsuit}$ denote the $\infty$-functor $\tau_{\geq 0}\circ \tau_{\leq 0} \simeq \tau_{\leq 0} \circ\tau_{\geq 0}$, and let $\pi_{n}:\C\to\C^{\heartsuit}$ denote the composition of $\pi_{0}$ with the shift $\infty$-functor $X\mapsto X[-n]$. 
\end{definition}

\begin{remark}[{\cite[Remark 1.2.1.12]{HA}}]
    Let $\C$ be a stable $\infty$-category; its heart $\C^{\heartsuit}$ is an abelian category. 
\end{remark}    

%%%%%%%%%%%%%%%%%%%%%%FINO A QUA CORRETTO %%%%%%%%%%%%

In the next part, we list the $t$-structures that are important to us.

\begin{example}[{\cite[Proposition 1.4.3.6]{HA}}]
\label{extstructureSp}
  Let $\Sp_{\leq 0}$ be the full $\infty$-subcategory of $\Sp$ spanned by spectra $X=(X_n)_{n\in\Z}$ whose $n$-component is an $n$-connected space. 
Let $\Sp_{\geq 0}$ be the full $\infty$-subcategory of $\Sp$ spanned by spectra $X=(X_n)_{n\in\Z}$ whose $n$-component is an $n$-connective space.
The triple $(\Sp,\Sp_{\geq 0}, \Sp_{\leq 0})$ forms a $t$-structure, sometimes called the Postnikov $t$-structure.
    In particular, $\Sp^{\heartsuit}\simeq \Ab$, and for each $n\in\Z$, we have the homotopy groups $\infty$-functors
    \[ \pi_{n}: \Sp\to \Ab. \]
    We recall that these $\infty$-functors factor as follows:
    \[
    \pi_{n}: \Sp\xrightarrow{\Omega^{\infty}\simeq\Hom_{\Sp}(\mathbb{S},-)}\mathcal{S}_{*}\xrightarrow{\pi_{n}} \Ab, n\geq 2;
    \]
    \[
    U\pi_{1}: \Sp\xrightarrow{\Omega^{\infty}\simeq\Hom_{\Sp}(\mathbb{S},-)}\mathcal{S}_{*}\xrightarrow{\pi_{1}} Grp;
    \]
    \[
    U\pi_{0}: \Sp\xrightarrow{\Omega^{\infty}\simeq\Hom_{\Sp}(\mathbb{S},-)}\mathcal{S}\xrightarrow{\pi_{0}} Set,
    \]
    where $U$'s are the canonical forgetfuls.
\end{example}
    
%\begin{example}
 %\label{extstructureD(k)-}
 
%Let $D(k)^-_{\leq 0}$ be the full $\infty$-subcategory of $\D(k)^-$ spanned by chain complexes with trivial positive homology groups, i.e., $H_{n}(C_*)\simeq 0$ if $ n>0 $.
%Let $D(k)^-_{\geq 0}$ be the full $\infty$-subcategory of $\D(k)^-$ spanned by chain complexes with trivial negative homology groups, i.e., $H_{n}(C_*)\simeq 0$ if $n<0$.
% The triple $(D(k)^-,D(k)^-_{\leq 0},D(k)^-_{\geq 0})$ forms a $t$-structure.  
% In particular $\D(k)^{\heartsuit}\simeq k\text{-}Mod$.
%\end{example}

%\begin{example}[{\cite[Proposition 1.3.5.21]{HA}}]
%\label{extstructureD(k)} Let $\D(k)_{\leq 0}$  be the full $\infty$-subcategory of $\D(k)$ spanned by chain complexes with trivial positive homology groups, i.e., $H_{n}(C_*)\simeq 0$ if $n>0 $. 
%Let and $\D(k)_{\geq 0}$ be the full $\infty$-subcategory of $\D(k)$ spanned by chain complexes with trivial negative homology groups, i.e., $H_{n}(C_*)\simeq 0$ if $n<0 $.
% The triple $(D(k),D(k)_{\leq 0},D(k)_{\geq 0}) $ forms a right complete $t$-structure.  
% In particular $\D(k)^{\heartsuit}\simeq k\text{-}Mod$.
%\end{example}

For us, there is another important stable $\infty$-category with a $t$-structure, but before we can consider it, we need a definition.

\begin{notation}[{\cite[Remark 7.1.0.3]{HA}}]
\label{notDisceteRings}
We say that a $\E$-ring $A$ is \textit{discrete} if for each $n\in\Z\text{-}\{0\}$, its $n^{th}$ homotopy group is trivial, i.e. $\pi_{n}A\simeq 0$.
We denote the full $\infty$-subcategory of $\Sp$ spanned by discrete $\E$-rings by $CAlg^{\heartsuit}$: this notation is not accidental and its meaning is imaginable but we do not explain it. 
\end{notation}

\begin{notation}
Let $A\in CAlg $ be an $\E$-ring; we call \textit{the underlying spectrum of $A$} its image via the canonical forgetful $\infty$-functor $U:CAlg\to \Sp$. 
Let $R\in CAlg $ be an $\E$-ring and $M\in \mathrm{Mod}_{R}$ be a $R$-module; we call \textit{underlying spectrum of $M$} its image via the forgetful $\infty$-functor $U:\mathrm{Mod}_R\to \Sp$. 
\end{notation}

\begin{remark}
\label{rmkHomotopyGroupRMod}
Let $R$ be an $\E$-ring; we can define its homotopy groups as the homotopy groups of its underlying spectrum; that is, for each $n\in\Z$, $\pi_n(M):=\pi_n(UM)$.
The action map $R\otimes R\to R$ endows the sums $\sum_{n\in\Z}\pi_*(R)$ with the structure of a graded ring.    

Let $M$ be a left $R$-module spectrum; we can define its homotopy groups as the homotopy groups of its underlying spectrum; for each $n\in\Z$, $\pi_n(M):=\pi_n(UM)$.
The action map $R\otimes M\to M$ endows the sums $\sum_{n\in\Z}\pi_*(M)$ with the structure of a graded left module over $\pi_{*}R$.

In particular, if $R$ is a discrete $\E$-ring, every homotopy group $\pi_{n}M$ has the structure of an (ordinary) $\pi_{0}R$-module; in other words, we have a collection of $\infty$-functors, for each $n\in\Z$, 
\[
\pi_n:\mathrm{Mod}_{R}\to \pi_{0}R\text{-}\mathrm{Mod}.
\]
 \end{remark}

Now that we have defined all the necessary tools, we can easily deduce that the Eilenberg-Maclane Spectrum $\Hk$ is discrete.

\begin{example}[{\cite[Remark 7.1.0.3]{HA}, \cite[ Example 6.27]{SchwedeSymSp}, \cite[\textsection 0]{SS}}]
\label{exEilenbergMaclaneSpectrumDiscrete}
Let $k$ be a commutative unitary ring.
By \Cref{rmkHomotopyGroupRMod}, the homotopy groups of an $\E$-ring $\chi$ are the homotopy groups of its underlying spectrum. In \Cref{extstructureSp}, we have seen that they are the homotopy groups of the space $\Omega^{\infty}(\chi)$.
Hence, the homotopy groups of the Eilenberg-Maclane Spectrum $\Hk=({K(k,n)}_{n\in\Z})$, whose $n$-component is the \textit{Eilenberg-Maclane space} of type $n$, are computed as follows:
\begin{equation}
\pi_{n}\Hk \cong \pi_{n}\Omega^{\infty}\Hk \cong\pi_{n}K(k,0)\cong \begin{cases}
    0 & \text{if } n\neq 0 \\
    k & \text{if } n =0.
\end{cases}
\end{equation}

So, $\Hk$ is discrete (\Cref{notDisceteRings}).

There is an equivalence of categories between the category of discrete $\mathbb{E}_{\infty}$-rings $CAlg^{\heartsuit}$ and the ordinary category of commutative rings:  
\[  CAlg^{\heartsuit}\simeq CRing  \]
\[   X\mapsto \pi_{0}X\]
\[ \Hk\mapsfrom k \]
that at every commutative ring $k$ it associates the \textit{Eilenberg-Maclane Spectrum} $\Hk=({K(k,n)}_{n\in\Z})$ whose $n$-component is the Eilenberg-Maclane space of type $n$.  

Since the above equivalence exists we will make an abuse of notation, which is common, and we will say \textit{$k$-linear $\infty$-categories} 
instead of \textit{$\Hk$-linear $\infty$-categories}.      \end{example}

%\begin{example}[{\cite[Proposition 7.1.1.5]{HA}}]
%\label{exStableModHK}

%\end{example}

\begin{example}[{\cite[Proposition 7.1.1.13,Proposition 7.1.1.5]{HA}}]
\label{exSperiamoUltimoDaCitare}
The $\infty$-category $\LModk$ is stable.  
Let $ \LModk^{\leq 0}$ be the full $\infty$-subcategory of $\LModk$ spanned by modules with trivial positive homotopy groups, i.e., $\pi_{n}(M)\simeq 0$ if $n>0 $.
Let $\LModk^{\geq 0}$ be the full $\infty$-subcategory of $\LModk$ spanned by modules with trivial negative homotopy groups, i.e., $\pi_{n}(M)\simeq 0$ if $n<0 $.
The triple $(\LModk,\LModk^{\leq 0},\LModk^{\geq 0}) $ forms a $t$-structure. 
In particular, $\LModk^{\heartsuit}\simeq k\text{-}Mod$.   
\end{example}

\begin{remark}
    Note that we could have also defined the $t$-structure before and then defined $\pi_n:\mathrm{Mod}_{R}\to \pi_{0}R\text{-}Mod$ using the $t$-structure, but our choice seems clearer to us.
Obviously, the two definitions of $\pi_n:\mathrm{Mod}_{R}\to \pi_{0}R\text{-}Mod$ are equivalent.
\end{remark}

\begin{notation}
\label{notRelativeTensorProduct}
Let $C$, $D$, $B$ be three $\E$-rings. We denote by $_B\otimes_{A}$ the relative tensor product between bimodules:
\begin{equation}
\begin{split}
-_B\otimes_{B}- : \; _C \mathrm{Mod}_B & \times \; _B \mathrm{Mod}_D \to \; _C \mathrm{Mod}_D: \\
(M & ,N)\mapsto Bar_{B}(M,N).
\end{split}
\end{equation}

Via the equivalences in \eqref{eqLModBModRMod}, for each $\E$-ring $A$, $\mathrm{RMod}_{A}$ and $\mathrm{LMod}_{A}$ inherit a monoidal structure. For example, $\mathrm{LMod}_{A}$ has as tensor the $\infty$-functor:
\[
-_{A}\otimes_{A} - : \mathrm{LMod}_{A} \times \mathrm{LMod}_{A} \simeq \; _{A}\mathrm{Mod}_{A} \times \; _{A}\mathrm{Mod}_{A} \xrightarrow{-\;_{A}\otimes_{A}-} \; _{A}\mathrm{Mod}_{A} \simeq \mathrm{LMod}_{A},
\]
and as unit $A$ seen as an object of $\mathrm{LMod}_{A}$.
Since the tensor of $\mathrm{LMod}_{A}$, $\mathrm{RMod}_{A}$, and $_{A}\mathrm{Mod}_{A}$ are all relative tensor products, we use the same notation for all of them.

For a detailed introduction about the relative tensor product, see \cite[\S 4.4]{HA}. In this paper, it is enough to consider the situation in \cite[Corollary 4.4.2.15]{HA}.
\end{notation}

\begin{example}
\label{exsmashproduct}
The smash product of $\Sp$, which we have already mentioned in \Cref{exOmegaInftyCambioArricchimento}, can be described as a relative tensor product.
Indeed, there is the equivalence $\Sp\simeq \mathrm{Mod}_{\mathbb{S}}$ and $\mathbb{S}$ is canonically an $\E$-ring.
\end{example}

\begin{example}
\label{exRelativeProductHK}
Let $k$ be a commutative unitary ring and $\Hk$ its associated Eilenberg-Maclane spectrum.
Specializing \Cref{notRelativeTensorProduct} in the situation $A=\Hk$, we obtain that the following triple $(\mathrm{Mod}_{\Hk},\; _{\Hk}\otimes_{\Hk},\;\Hk)$ forms a symmetric monoidal $\infty$-category.
\end{example}

In the last part of this subsection, we recall the presentable $\infty$-category theory, \cite[\S 5.5]{HTT}.

\begin{definition}[{\cite[Definition 5.4.2.1]{HTT}}] 
\label{defAccessible}
Let $\gamma$ be a regular cardinal. An $\infty$-category $\C$ is $\gamma$-accessible if there exists a small $\infty$-category $\C_0$ and an equivalence $Ind_{\gamma}(\C_0) \simeq \C$.
We will say that $\C$ is accessible if it is $\gamma$-accessible for some regular cardinal $\gamma$.
\end{definition}

\begin{definition}
\label{defPresentable}
An $\infty$-category $\C$ is presentable if it admits small colimits and is accessible.
We say that an $\infty$-category $\C$ is \textit{compactly generated} if it is presentable and $\omega$-accessible.
\end{definition}

\begin{remark}
In \cite[5.5.1]{HTT}, it is shown that there is a list of alternative definitions of presentable $\infty$-categories. Most of them are found by Simpson.     
\end{remark}

\begin{notation}
We denote by $\widehat{\Cat}$ the $\infty$-category of non-necessarily small $\infty$-categories.
\end{notation}

\begin{notation}
We denote by $P_r^L\subseteq \widehat{Cat}_\infty$ the $\infty$-subcategory spanned by presentable $\infty$-categories and whose morphisms are small-colimits-preserving functors.
%We denote by $Pr^{L}_{cg}\subseteq Pr^{L}$ the full $\infty$-subcategory spanned by compactly generated $\infty$-categories.
\end{notation}

\begin{example}[{\cite[Proposition 1.4.3.7,  Proposition 7.2.4.2.]{HA}}]
\label{exPresentable}
The $\infty$-categories $\LModk$, and $\Sp$ are presentable. Actually, they are compactly generated. 
\end{example}

In \cite[Proposition 4.8.1.15]{HA}, the author proves that $Pr^{L}$ inherits a closed symmetric monoidal structure such that the inclusion $\iota: Pr^{L}\to \widehat{\Cat}$ is lax monoidal. 
The closed monoidal structure in $Pr^{L}$ has as unit the $\infty$-category $\Top$, as hom-$\Cat$-object between two presentable $\infty$-categories $\A,\B$ the $\infty$-category $LFun(\A,\B)$ and the tensor is defined as follows:
\begin{equation}
\label{eqTensorPrl}
\begin{split}    
-\otimes^{L}-:Pr^{L} & \times Pr^{L}\to Pr^{L}\\ (\A & ,\B)\mapsto \A\otimes^{L}\B:=RFun(\A^{op},\B),
\end{split}
\end{equation}
and it preserves small colimits componentwise.  

The algebra objects in the symmetric monoidal $\infty$-category $Pr^{L}$ are important, so they deserve a specific name.

\begin{definition}
\label{defpresentablymonoida}
 A \textit{presentably monoidal $\infty$-category} is an associative algebra object of $Pr^{L}$, see \Cref{defAlgebraHigherAlgebra}.    
\end{definition}

\begin{remark}
In \cite{GepHauEnriched}, the authors give an equivalent definition of presentably monoidal $\infty$-category, which is clearer and so we report here: let $\C$ be a monoidal $\infty$-category, we say that it is \textit{presentably monoidal} if the tensor is compatible with small colimits and $\C$ is presentable.
For the definition of compatible with small colimits see \cite[Denition 3.1.1.18 and Variant 3.1.1.19]{HA}.
\end{remark}

\begin{example}[{\cite[Corollary 4.8.2.19,Corollary 4.4.2.15]{HA}
}]
\label{exPresentableMonoidal}
The $\infty$-categories $\LModk$, and $\Sp$ are presentably monoidal $\infty$-category with the monoidal structure in \Cref{exRelativeProductHK} and \Cref{exsmashproduct} (or \Cref{exOmegaInftyCambioArricchimento}).
\end{example}

\begin{example}[{\cite[Corollary 3.2.3.5, Example 3.2.4.4 and Proposition 3.2.3.1]{HA}}]
\label{exMonoidalStructureAlgSpCAlgSp}
The $\infty$-categories $CAlg$ and $Alg$ are presentably monoidal $\infty$-category, and the canonical forgetful $\infty$-functors: 
\[Alg\to \Sp,\]
and
\[CAlg\to \Sp,\]
are monoidal.
\end{example}

\begin{definition}[{\cite[Definition 3.2.10 and 3.2.9]{GepHauEnriched}}]
We denote by $Mon^{Lax}$ the $\infty$-category of monoidal $\infty$-categories and lax monoidal $\infty$-functors.

We denote by $PMon^{Lax}\subseteq Mon^{Lax}$ the full $\infty$-subcategory spanned by presentably monoidal $\infty$-categories.

\end{definition}

\begin{definition}[{\cite[Definition 3.1.24]{GepHauEnriched}}]
Let $Mon^{Pr}$ be the $\infty$-category $Alg_{Ass}(P^{L}_r)$ of associative algebra objects in
$P^L_r$ equipped with the tensor product in \Cref{eqTensorPrl}. Thus $Mon^{Pr}$
 is the $\infty$-category
of monoidal $\infty$-categories $\C^\otimes$ compatible with small colimits such that $\C$ is presentable, with
morphisms monoidal $\infty$-functors that preserve colimits. .
\end{definition}
 The objects of $Mon^{Pr}$ are
presentably monoidal $\infty$-categories and it is canonically a (non-full) $\infty$-subcategory $Mon^{Pr}\subseteq PMon^{Lax}$.

%%%%%%%%%%%%%%%%%%%%%%%%%%%%%%%%%%%%%%%%%%%%%%%%%%%%%%%%%%%%%%%%%%%%%%%%%%%%%%%%%%%%%%%%%%%%%%%%%%%%%%%%%%%%%%%%%%%%%%%%%%%%%%%%%%%%%%%%%%%%%%%%%%%%%%%%%%%%%%%%%%%%%%%%%%%%%%%%%%%%%%%%%%%%%%%%%%%%%%%%%%%%%%%%%%%%%%%

\subsection{Technical non-sense}
\label{subsecTecnialInftyCase}

 In this subsection, we study the categorical version of the extension-restriction and restriction-coextension adjunction.

\begin{lemma}
\label{LemmaHighLeftRightAdjointforgMod}
Let $\C$ be a monoidal $\infty$-category, and let $\mathcal{M}$ be a $\C$-left tensored $\infty$-category. Assume that $\C$ and $\mathcal{M}$ admit small colimits and the tensor $\infty$-bifunctors:
\[-\otimes_{\C}-:\C\times\C\to \C;  \]
\[-\;^{\C}\otimes_{\mathcal{M}}-: \C\times \mathcal{M}\to \mathcal{M},\]
preserve small colimits componentwise. Then, for every map $f: B \to A$ in $Alg(\C)$, the forgetful functor $\mathrm{LMod}_A(\mathcal{M}) \to \mathrm{LMod}_B(\mathcal{M})$ admits a left adjoint, given by the relative tensor product construction
\[\mathrm{LMod}_B(\mathcal{M}) \to \mathrm{LMod}_A(\mathcal{M})      :\]
\[m\mapsto \;_{A}B_{B}\;_{B}\otimes_{B} m.\]
\end{lemma}
\begin{proof}
  This lemma is a subcase of \cite[Proposition 4.6.2.17]{HTT}.  
\end{proof}

The following example is useful for us.

\begin{example}
\label{exTensorHkForgetful}
$\Sp$ is a presentable $\infty$-category, and the smash product preserves small colimits componentwise \cite[Corollary 4.8.2.19]{HA}. Canonically, $\Sp$ is a $\Sp$-left tensored $\infty$-category. We are in the hypothesis of \Cref{LemmaHighLeftRightAdjointforgMod}. We apply the lemma with $f$ being the unit morphism of the Eilenberg-Maclane spectrum $\Hk$ $u:\mathbb{S}\to \Hk$, which is a morphism in $Alg$. We obtain the restriction-extension adjunction:

\begin{equation}
\label{eqTensorHkForgetful}
-\otimes \Hk:\Sp\simeq \mathrm{LMod}_{\mathbb{S}} \rightleftarrows \mathrm{LMod}_{\Hk}: U 
\end{equation}
This adjunction is equivalent to finding the left adjoint to a forgetful functor.
\end{example}

For the adjunction (\ref{eqTensorHkForgetful}), we need the right adjoint to the forgetful, and our next result explicitly finds it.

\begin{remark}
We can say more about the right adjoint above; indeed, it can be described in two ways:
\begin{itemize}
\label{rmkForgetisRepresented}
    \item 
    Let $R\in Alg$ be an associative algebra (for us, $R$ will almost always be $\Hk$). Then the forgetful functor $U:\mathrm{LMod}_{R}\to \Sp$ is the $\Sp$-enriched representable $\infty$-functor 
    \[ \mathrm{Mod}_{R}(R,-):\mathrm{Mod}_{R}\to \Sp \] 
    where, for each $M $ element of $\mathrm{LMod}_{R}$ and for each $n\in\N$,
    \[\mathrm{LMod}_{R}(R,M)_{-n}= \mathrm{Hom}_{\mathrm{LMod}_{R}}(R,[n](M)).\]
    This result is expected because every presentable stable category is an $\Sp$-enriched category with $\mathrm{LMod}_{R}(-,-)$ as hom-$\Sp$-object (see \cite[Theorem 1.2]{heine2023equivalence}). 
Let $-\otimes R\dashv U:\Sp\to \mathrm{LMod}_{R}$ be the tensor-forgetful adjunction (\ref{eqTensorHkForgetful}), the $\infty$-functor $\Omega_{\infty}U$ is equivalent to the representable $\infty$-functor $\Hom_{\mathrm{LMod}_{R}}(R,-)$. Indeed, for each element $M$ in $\mathrm{LMod}_{R}$, we have the following chain of equivalence:
\[\Omega_{\infty}U(M)\simeq \Hom_{\Sp}(\mathbb{S},UM)\simeq\Hom_{\mathrm{LMod}_{R}}(\mathbb{S}\otimes R,M)\simeq\Hom_{\mathrm{LMod}_{R}}(R,M),\]
where the first holds by definition, the second holds by the universal property of the adjunction $-\otimes R\dashv U$, and the last because $\mathbb{S}$ is the unit of the monoidal structure of $\Sp$.

Now by definition of the various $\pi_0$'s, $\pi_0U\simeq \pi_0$, \Cref{extstructureSp}, the following isomorphism $\pi_0\Omega_{\infty}\cong \pi_0$ holds, which implies that $U$ and $\mathrm{Mod}_{R}(R,-)$ coincide.
%    In the proof of \Cref{propDgCatePreadditive}, we will prove this fact;
    \item The forgetful $\infty$-functor \[U:\mathrm{LMod}_{R}\to \Sp\] is equivalent to the $\infty$-functor giving the tensor of bimodules:  
    \[ R\;_R\otimes_{R}-: \mathrm{LMod}_{R}\simeq \mathbb{1}\times \mathrm{LMod}_{R} \xrightarrow{R\times id}\;_\mathbb{S}\mathrm{Mod}_{R} \times _{R}\mathrm{Mod}_\mathbb{S}\xrightarrow{ \;_{R}\otimes_{R}}\Sp.  \]
    This is the typical situation because $R$ is a dualizable element of $\mathrm{LMod}_{R}$. Hence the right adjoint of the tensor product $\infty$-functor $-\otimes_{\mathbb{S}}R $ (which is the internal-hom) may be described as the tensor product with the dual of $R$, \cite[Proposition 4.6.2.1, Proposition 7.2.4.4]{HA} or \cite{hoyois2017higher}, and in this case, the dual is $R$ itself but considered as a right $R$-module.   
\end{itemize}
\end{remark}

\begin{lemma}
\label{LemmaRightaadjointforgetful}
Let $R\in Alg$ be an associative algebra object. The forgetful $\infty$-functor \[U:\mathrm{LMod}_{R}\to \Sp;\] has a right adjoint.
\end{lemma}
\begin{proof}
The $\infty$-functor $U$ is equivalent to the $\infty$-functor $R \;_R\otimes_{R}-:\;_R \mathrm{Mod}\to \Sp$ where $R$ is considered as a right $R$-module, see \Cref{rmkForgetisRepresented}.
$\mathrm{Mod}_{R}$ and $\Sp$ are both compactly generated; then, if we show that $U$ preserves small colimits, then using the famous Adjoint Functor Theorem, \cite[Corollary 5.5.2.9]{HTT}, we get the thesis.
The smash product in $\Sp$ preserves small colimits, \cite[Corollary 4.8.2.19]{HA}, then the functor $R\;_R\otimes_{R}-$ preserves small colimits, \cite[Corollary 4.4.2.15]{HA}, so the thesis.
\end{proof}

%%%%%%%%%%%%%%%%%%%%%%%%%%%%%%%%%%%%%%%%%%%%%%%%%%%%%%%%%%%%%%%%%%%%%%%%%%%%%%%%%%%%%%%%%%%%%%%%%%%%%%%%%%%%%%%%%%%%%%%%%%%%%%%%%%%%%%%%%%%%%%%%%%%%%%%%%%%%%%%%%%%%%%%%%%%%%%%%%%%%%%%%%%%%%%%%%%%%%%%%%%%%%

%% file: Chapters/EnrichedCategory.tex
\section{Enriched $\infty$-categories}
\label{ChapEnrichCate}

Throughout this introduction $\mathcal{V}$ is a monoidal symmetric category.

Nowadays, there are three different definitions of enriched $\infty$-categories: V. Hinich's \cite{HinichYon}, D. Gepner's and R. Haugseng's \cite{GepHauEnriched}, and J. Lurie's \cite{HA}.

The idea behind these definitions is also present in ordinary category theory. Let $\mathcal{V}$ be a symmetric monoidal category, the first two are the vision of a $\mathcal{V}$-enriched category, with $S$ as the set of objects, as an $\text{Ass}^{\mathcal{V}}$-algebra in the category $\mathcal{V}$-graphs, with $S$ as the set of objects, $Graph_{S}(\mathcal{V})$; $\text{Ass}^{\mathcal{V}}$ is the (non-symmetric) associative operad in $\mathcal{V}$ (see \cite{MuroHomOp}). Instead, the idea behind the last is that enriched categories are tensored $\mathcal{V}$-enriched categories \enquote{without tensor}. At first glance, this can be strange but in the theory of $\infty$-operads, this concept has a precise meaning: we lift only the arrows (inert maps) that shape the structure that may contain the tensor but we do not lift the arrows (active maps) that define the tensor, see \cite[\S 2]{HA} or \cite[\S 2.6]{HinichYon}.

In this paper, we take advantage of this double vision to write an $\LModk$-enriched $\infty$-category, with some conditions, as a left $\LModk$-module object of $\Cat^{\Sp}$. Hence we need a way to compare the different definitions.

In recent years, A. Macpherson in \cite{Macpherson} proved that Hinich's and Haugseng-Gepner's definitions are equivalent and H. Heine in \cite{Heine} proves that Hinich's and Lurie's definitions are equivalent. Hence, as one would expect, they are mathematically the same concept.

%We are particularly interested in Heine's work. His comparison is the starting point in the proof of our Main Result.

For us, unless otherwise specified, an $\V$-enriched $\infty$-category is an algebra in the (non-symmetric) $\infty$-operad of graphs $Quiv_X(\V)$ with spaces of objects the $\infty$-groupoid $X$, \Cref{defHigherenrichedCategory}.

%In classical categorical homotopy theory, one can use the Enriched Yoneda embedding to embed a (small) dg-category over $k$, $\mathcal{C}$, in the smallest pretriangulated dg-category which contains itself. 
%The above construction is functorial and defines a localization from the category of dg-categories into the category of pretriangulated dg-categories. 

Any generalization of the ordinal category theory, to be well-formed, needs its Yoneda embedding. Today there are two possible definitions of enriched categorical Yoneda Embedding, Hinich's and Heine's, and, as far as we know, it has not yet been proved that they are equivalent.
But, it is folklore, that the two Yoneda embeddings are the same and that the difference is that H. Heine defines the Yoneda embedding between two enriched $\infty$-categories, instead Hinich defines it as a \enquote{$\infty$-functor with source a $\mathcal{V}$-enriched $\infty$-category and as a target $\mathcal{V}$-left tensored $\infty$-category}, see \Cref{subsecEnrichedCatFun}.
In \cite{DoniklinearMorita}, which belongs to this series of articles, we work with presentable stable $\infty$-categories which are left $\Sp$-module objects of $Pr^{L}$, then Hinich's definition is convenient for us.
Thanks to this choice the proof of \Cref{thAction} is almost tautological.

%%%%%%%%%%%%%%%%%%%%%%%%%%%%%%%%%%%%%%%%%%%%%%%%%%%%%%%%%%%%%%%%%%%%%%%%%%%%%%%%%%%%%%%%%%%%%%%%%%%%%%%%%%%%%%%%%%%%%%%%%%%%%%%%%%%%%%%%%%%%%%%%%%%%%%%%%%%%%%%%%%%%%%%%%%%%%%%%%%%%%%%%%%%%%%%%%%%%%%%%%%%%%%%%%%%%%%%%%%%%%%%%%%%%%%%%%%%%%%%%%%

%%%%%%%%%%%%%%%%%%%%%%%%%%%%%%%%%%%%%%%%%%%%%%%%%%%%%%%%%%%%%%%%%%%%%%%%%%%%%%%%%%%%%%%%%%%%%%%%%%%%%%%%%%%%%%%%%%%%%%%%%%%%%%%%%%%%%%%%%%%%%%%%%%%%%%%%%%%%%%%%%%%%%%%%%%%%%%%%%%%%%%%%%%%%%%%%%%%%%

\subsection{Enriched $\infty$-categories}
\label{subsecEnrichedCatFun}

\begin{definition}[{\cite[Definition 3.2.10 and 3.2.9]{GepHauEnriched}}]
We denote by $Mon^{Lax}$ the $\infty$-category of monoidal $\infty$-categories and lax monoidal $\infty$-functors.

We denote by $PMon^{Lax}\subseteq Mon^{Lax}$ the full $\infty$-subcategory spanned by presentably monoidal $\infty$-categories.

\end{definition}

\begin{definition}[{\cite[Definition 3.1.24]{GepHauEnriched}}]
Let $Mon^{Pr}$ be the $\infty$-category $Alg_{Ass}(P^{L}_r)$ of associative algebra objects in
$P^L_r$ equipped with the tensor product in \Cref{eqTensorPrl}. Thus $Mon^{Pr}$
 is the $\infty$-category
of monoidal $\infty$-categories $\C^\otimes$ compatible with small colimits such that $\C$ is presentable, with
morphisms monoidal $\infty$-functors that preserve colimits. .
\end{definition}
The objects of $Mon^{Pr}$ are
presentably monoidal $\infty$-categories and it is canonically a (non-full) $\infty$-subcategory $Mon^{Pr}\subseteq PMon^{Lax}$.
$\\$

In this paper, we work with $\infty$-categories enriched in presentably monoidal $\infty$-category, \Cref{defAlgebraHigherAlgebra}, this type of enrichment is widely studied in \cite{GepHauEnriched}.

In this section, we fix a presentably monoidal $\infty$-category $\V$;
 we work in a slightly more general setting than what is needed in the results in \Cref{secLModEnrichedareModuleObjects}, so the reader, if he is only interested in the results of this paper, can replace all $\V$ by $\LModk$ or $\mathrm{LMod}_{\mathbb{S}}(\Sp)\simeq \Sp$.

%Our definition of an enriched $\infty$-category is Hinich's definition, see \Cref{defHigherenrichedCategory}. We chose it because Hinich, in \cite{HinichYon}, defines the $\V$-enriched \categorical Yoneda embedding; but it would have been possible to use any other definition. 

There is a weak condition with respect to the enriched $\infty$-categories: \precategories. We do not use this generality, but we need this concept to define what an enriched $\infty$-category is. The difference between enriched $\infty$-categories and \precategories is similar to the difference between Segal space and complete Segal space (see \cite{Rezk} for this more classical case). One must add the complete condition to have the correct notion of $\V$-equivalence between $\V$-enriched $\infty$-categories; i.e. $\V$-enriched $\infty$-functor which are homotopically essentially surjective and homotopically fully-faithful, see \cite[\textsection 5.3]{GepHauEnriched}. In other words, they are the \categorical analogous of Dwyer-Kan model categories, \cite[\S 1]{MuroDK} or \cite[Remark 3.1.20]{DoniPhDThesis}, as Haugseng has proved in \cite[Theorem 5.8]{HaugsengRect}. 

Hinich's definition is based on a $BM$-monoidal $\infty$-category (or, $BM$-operad) $Quiv_X(M)$, which we recall here: $BM$-monoidal $\infty$-category are in particular $\infty$-categories over $BM$, that why Hinich uses the method described in \cite[\S 2]{DoniPhDThesis} to define $Quiv_X(M)$.

Our introduction to enriched $\infty$-categories is far from being technical and detailed, if the reader wants to know them, he has to read \cite{HinichYon}.
 $\\$
 
Since $Pr^{L}$ is a monoidal $\infty$-category, canonically, it is possible to consider it as $BM$-monoidal $\infty$-category.

Let $M=(M_{a}, M_{m}, M_{b})$ be an bimodules object of $Pr^{L}$. We denote by $\otimes_{a,m}^{M}$ the left action and by $\otimes_{m,b}^{M}$  the right action; if it is not confusing we will delete subscripts or superscripts. 

In \cite[Theorem 4.4.8]{HinichYon}, the author proves that, for every space $X$, there is a $Quiv_{X}(M_{a})$-$M_{b}$-bimodule object of $Pr^{L}$ denoted  \[Quiv_{X}^{BM}(M)=(Quiv_X(M_{a})= Fun(X^{op}\times X,M_{a}), Fun(X,M_{m}), M_{b}),\]  where all the actions have a nice description using coend.% and left Kan extension theories

\begin{itemize}
    \item[(i)] The first component is an associative algebra of $Pr^{L}$ $(Fun(X^{op}\times X,M_{a}),\otimes_{a})$, usually denoted $Quiv_{X}(M_{a})$, where the associative operation behaves as follows:    
    \begin{equation}
     \label{eqActionFun(XxX,M)xFun(XxX,M)} 
     \begin{split}
    -\otimes_{a}-:Fun( X^{op}\times X,M_{a}  )\times Fun(X^{op}\times & X,M_{a})\to Fun(X^{op},X,M_{a}): \\ (A_1 & ,A_2 )\mapsto A_1\otimes_{a}A_2,
    \end{split}
    \end{equation}
     for each pair of objects $(x_1,x_2)$,
\[ A_1\otimes_{a}A_2(x_1,x_2):= \oint^{y\in X} A_1(x_1,y)\otimes^{M_a}_{a} A_2(y,x_2);  \]
%In other word, it is the left Kan extension drawn below:
% https://q.uiver.app/#q=WzAsNSxbMCwxLCJYXFx0aW1lcyBYXFx0aW1lcyBYIl0sWzEsMCwiWFxcdGltZXMgWCJdLFsyLDEsIihYXFx0aW1lcyBYKVxcdGltZXMgKFhcXHRpbWVzIFgpIl0sWzMsMSwiTV97YX1cXHRpbWVzIE1fe2F9Il0sWzQsMSwiTV97YX0iXSxbMCwxLCJcXGdhbW1hIl0sWzAsMiwiKFxcYWxwaGEsXFxiZXRhKSIsMl0sWzIsMywiKEFfMSxBXzIpIiwyXSxbMyw0LCJcXG90aW1lc197YX0iLDJdLFsxLDQsIkFfMVxcb3RpbWVzX3thfUFfMSJdLFs2LDEsIiIsMCx7InNob3J0ZW4iOnsic291cmNlIjoyMH19XV0=
%\[\begin{tikzcd}[sep=small]
%	& {X\times X} \\
%	{X\times X\times X} && {(X\times X)\times (X\times X)} & {M_{a}\times M_{a}} & {M_{a},}
%	\arrow["\gamma", from=2-1, to=1-2]
%	\arrow[""{name=0, anchor=center, inner sep=0}, "{(\alpha,\beta)}"', from=2-1, to=2-3]
%	\arrow["{(A_1,A_2)}"', from=2-3, to=2-4]
%	\arrow["{-\otimes_{a}-}"', from=2-4, to=2-5]
%	\arrow["{A_1\otimes_{a}A_1}", from=1-2, to=2-5]
%	\arrow[shorten <=3pt, Rightarrow, from=0, to=1-2]
%\end{tikzcd}\]

%here, the $\infty$-functor of spaces $(\alpha,\beta)$ acts as follows $(\alpha,\beta):(x_1,x_2,x_3)\mapsto (x_1,x_2,x_2,x_3)$, and $\gamma:(x_1,x_2,x_3)\mapsto (x_1,x_3) $ is the projection;

\item[(ii)] the chosen objects in $Pr^{L}$ is  $Fun(X,M_{m})$;

\item[(iii)] the third component is the associative algebra of $Pr^{L}$ $M_{b}$ with the starting operation;

\item[(iv)] the left action of $Quiv_{X}(M_{a})$ over $Fun(X,M_m)$ 
\begin{equation}
\label{eqLeftActionFun(xm)}
\begin{split}
-\otimes_{a,m}-:Fun(X^{op}\times X,M_{a} & )\times Fun(X,M_{m})\to Fun(X,M_{m}):\\ (A & ,F)\mapsto A\otimes_{a,m}F 
\end{split}
\end{equation}
can be described very well using the coend.
For each $x\in X$, \[A\otimes_{a,m}F(x):=\oint^{c\in X}A(c,x)\otimes_{a,m}^{M} F(c);  \]

%Another time this is a left Kan extension: 
% https://q.uiver.app/#q=WzAsNSxbMCwxLCJYXFx0aW1lcyBYIl0sWzEsMCwiIFgiXSxbMiwxLCIoWFxcdGltZXMgWClcXHRpbWVzICggWCkiXSxbMywxLCJNX3thfVxcdGltZXMgTV97bX0iXSxbNCwxLCJNX3ttfSJdLFswLDEsIlxcZ2FtbWEiXSxbMCwyLCIoaWQsXFxiZXRhKSIsMl0sWzIsMywiKEEsRikiLDJdLFszLDQsIlxcb3RpbWVzX3thLG19IiwyXSxbMSw0LCJBXFxvdGltZXNfe2EsbX1GIl0sWzYsMSwiIiwwLHsic2hvcnRlbiI6eyJzb3VyY2UiOjIwfX1dXQ==
%\[\begin{tikzcd}
%	& { X} \\
%	{X\times X} && {(X\times X)\times ( X)} & {M_{a}\times M_{m}} & {M_{m};}
%	\arrow["\gamma", from=2-1, to=1-2]
%	\arrow[""{name=0, anchor=center, inner sep=0}, "{(id,\beta)}"', from=2-1, to=2-3]
%	\arrow["{(A,F)}"', from=2-3, to=2-4]
%	\arrow["{-\otimes_{a,m}-}"', from=2-4, to=2-5]
%	\arrow["{A\otimes_{a,m}F}", from=1-2, to=2-5]
%	\arrow[shorten <=3pt, Rightarrow, from=0, to=1-2]
%\end{tikzcd}\]
\item[(v)] instead the right action is the diagonal action
\begin{equation*}
\begin{split}
-\otimes_{m,b}-: Fun(X,M_{m})\times & M_{b}\to Fun(X,M_{m}): \\ ( & F,B)\mapsto F\otimes_{m,b}B; 
\end{split}
\end{equation*}for each  $x\in X$, 
\begin{equation}
\label{eqRightAction}
    F\otimes_{m,b}B(x)=F(x)\otimes_{a,m}^{M} M.
\end{equation}

%This also has its description as left Kan extension: 
% https://q.uiver.app/#q=WzAsNCxbMCwxLCJYXFxzaW1lcSBYXFx0aW1lcyBcXG1hdGhiYnsxfSJdLFsxLDAsIiBYIl0sWzMsMSwiTV97bX1cXHRpbWVzIE1fe2J9Il0sWzQsMSwiTV97bX0iXSxbMCwxLCJpZCIsMCx7InN0eWxlIjp7ImhlYWQiOnsibmFtZSI6Im5vbmUifX19XSxbMiwzLCJcXG90aW1lc197bSxifSIsMl0sWzEsMywiQVxcb3RpbWVzX3thLG19RiJdLFswLDIsIihGLE0pIiwyXSxbNywxLCIiLDIseyJzaG9ydGVuIjp7InNvdXJjZSI6MjB9fV1d
%\[\begin{tikzcd}
%	& { X} \\
%	{X\simeq X\times \mathbb{1}} &&& {M_{m}\times M_{b}} & {M_{m}.}
%	\arrow["id", no head, from=2-1, to=1-2]
%	\arrow["{-\otimes_{m,b}-}"', from=2-4, to=2-5]
%	\arrow["{F\otimes_{a,m}M}", from=1-2, to=2-5]
%	\arrow[""{name=0, anchor=center, inner sep=0}, "{(F,M)}"', from=2-1, to=2-4]
%	\arrow[shorten <=4pt, Rightarrow, from=0, to=1-2]
%\end{tikzcd}\]

\end{itemize}

The above structure is fundamental for this series of articles.
$(i)$ allows to define $\LModk$-enriched $\infty$-categories, $(iv)$ permits to define the $\infty$-category of $\LModk$-enriched $\infty$-functors $Fun^{\LModk}(\C,\D)$. In particular, the $\infty$-category of $\LModk$-enriched presheaves on $\C$, $P_{\LModk}(\C)$, and $(iii)$ permits to define a left action of $\LModk$ on $P_{\LModk}(\C)$: but this paper is not concerned with the $\LModk$-enriched presheaves on $\C$ and its associated structures do not concern this paper.

\begin{remark}
\label{rmkHigherAlgebraicStructureInQuivBM}
Again, we do not prove the following statements, but we only list them. For a proof, see \cite[Corollary 4.4.9]{HinichYon}.
Let $M=(M_{a}, M_{m}, M_{b})$ be an bimodules object in $Pr^{L}$. The structure of $Quiv_{X}^{BM}(\V)$ contains the following higher algebraic structures:
\begin{itemize}
    \item[(0)] $(Fun(X^{op}\times X,M_{a}),\otimes_{a})$ is an associative algebra in $Pr^{L}$. Usually, it is denoted $Quiv_{X}(M)$. In particular, since $Pr^{L}\to \Cat$ is lax monoidal, $Quiv_{X}(M)$ is also a monoidal $\infty$-category;
    
    \item[(1)] $Fun(X,M_{m})$ is a left $Quiv_{X}(M_{a})$-module object in $Pr^{L}$. This left $Quiv_{X}(M_{a})$-module object is usually denoted $Quiv^{LM}_{X}(M)$. In particular, since $Pr^{L}\to \Cat$ is lax monoidal, $Quiv^{LM}_{X}(M)$ is also a $Quiv_{X}(M_{a})$-left tensored $\infty$-category;
    
    \item[(2)] $Fun(X,M_{m})$ is a right $M_{b}$-module object in $Pr^{L}$. This right $M_{b}$-module object is usually denoted $Quiv^{RM}_{X}(M)$. In particular, since $Pr^L\to \Cat$ is lax monoidal, $Quiv^{RM}_{X}(M)$ is a is also a $M_{b}$-right tensored $\infty$-category;
    
    \item[(3)] $Fun(X,M_{m})$ is a $Fun(X^{op}\times X,M_{a})$-$M_{b}$-bimodule object of $Pr^{L}$. This $Fun(X^{op}\times X,M_{a})$-$M_{b}$-bimodule object is usually denoted $Quiv^{BM}_{X}(M)$. In particular, since $Pr^{L}\to \Cat$ is lax monoidal, $Quiv^{RM}_{X}(M)$ is also a $Fun(X^{op}\times X,M_{a})$-$M_{b}$-bitensored $\infty$-category.
\end{itemize}

Note that to construct $Quiv^{LM}_{X}(M)$, $Quiv^{RM}_{X}(M)$, or $Quiv_{X}(M)$, we do not need bimodule object of $Pr$ $M$ but it is enough to have an $LM$-module object, $RM$-module object, or an associative algebra.
\end{remark}

The construction of $Quiv^{*}_{X}(M)$, with $*\in \{BM,LM,\emptyset\}$, is functorial in $X$ and $M$. Indeed, it defines an $\infty$-functor,

\begin{equation}
\label{eqQuivFunctorialMonLax}
Quiv^{*}: \mathcal{S}^{op} \times Mon^{Lax} \to Mon^{Lax}.
\end{equation}

Thanks to \cite[Proposition 4.3.5, Corollary 4.3.16]{GepHauEnriched}, we can restrict the above $\infty$-functor as follows:

\begin{equation}
\label{eqQuivFunctorialPMonLax}
Quiv^{*}_{(-)}(-): \mathcal{S}^{op} \times PMon^{Lax} \to PMon^{Lax};
\end{equation}
and 
\begin{equation}
\label{eqQuivFunctorialMonPr}
Quiv^{*}_{(-)}(-): \mathcal{S}^{op} \times Mon^{Pr} \to Mon^{Pr}.
\end{equation}

\begin{notation}
In this paper, we use the same notation for \eqref{eqQuivFunctorialMonLax}, \eqref{eqQuivFunctorialPMonLax}, and \eqref{eqQuivFunctorialMonPr}.
\end{notation}

To define what is an enriched $\infty$-category, we need only the case where $*=\emptyset$. As mentioned earlier, we define precategories before defining enriched $\infty$-categories.

\begin{definition}
\label{defPreca}
Let $\V$ be a presentably monoidal $\infty$-category, and let $X$ be a space. An $\V$-enriched precategory with the space of objects $X$ is an associative algebra in $Quiv_X (\V)$.
\end{definition}

\begin{remark}
In \Cref{rmkHigherAlgebraicStructureInQuivBM}, we have seen that $Quiv_{X}(\V)$ is an associative object in $Pr^{L}$, and we have never defined what an associative algebra is for such algebraic structure.

But since the inclusion $Pr^{L}\to \Cat$ is lax monoidal, an associative object in $Pr^{L}$ is an associative object in $\Cat$, which is a monoidal $\infty$-category. Then in \Cref{defPreca}, we mean an associative algebra in the monoidal $\infty$-category $Quiv_{X}(\V)$.
\end{remark}

We define the $\infty$-functor $PreC(-)$ as the composition of $Quiv_{(-)}(-)$ with the $\infty$-functor:
\[Alg_{Ass}(-):=Fun^{coCart}_{Ass}(id_{Ass},-): Mon^{Lax}\to \Cat\]  
or
\[Alg_{Ass}(-):=Fun^{coCart}_{Ass}(id_{Ass},-):PMon^{Lax}\to \Cat\]  
or 
\[Alg_{Ass}(-):=Fun^{coCart}_{Ass}(id_{Ass},-):Mon^{Pr}\to {Pr^L};\]
where for each presentably monoidal $\infty$-category $\W\to Ass$, $Fun_{Ass}^{coCart}(id_{Ass},\W)$ is the full $\infty$-subcategory of $Fun_{\mathcal{C}at_{\infty/Ass}}(id_{Ass},\W)$ spanned by the $\infty$-functors which sends cocartesian morphisms in cocartesian morphisms.

And we obtain the $\infty$-functors:

\begin{equation}
\label{eqVaccaLogia0}
PreC_{(-)}(-):=Alg_{Ass}(Quiv_{(-)}(-)): \mathcal{S}^{op} \times Mon^{Lax}\to \Cat; 
\end{equation}

and

\begin{equation}
\label{eqVaccaLogia}
PreC_{(-)}(-):=Alg_{Ass}(Quiv_{(-)}(-)): \mathcal{S}^{op} \times PMon^{Lax}\to \Cat;
\end{equation}

and 

\begin{equation}
\label{eqVaccaLogia1}    
PreC_{(-)}(-):=Alg_{Ass}(Quiv_{(-)}(-)): \mathcal{S}^{op} \times Mon^{Pr}\to {Pr^{L}}.  
\end{equation}

\begin{notation}
We denote the $\infty$-functors \eqref{eqVaccaLogia0}, \eqref{eqVaccaLogia} and \eqref{eqVaccaLogia1} in the same way.  
\end{notation}

Now, we define the $\infty$-category of precategories.
Let $(\V,\otimes_{\V},I_{\V})$ be a presentably monoidal $\infty$-category. Using the straightening-unstraightening, we obtain the cartesian fibration \[\int PreC(\V)\to \mathcal{S}.\]

\begin{definition}
Let $\V$ be a presentably monoidal $\infty$-category. We call the total spaces $\int PreC(\V)$ the \textit{$\infty$-category of $\V$-precategories} and we denote it by $PreC(\V)$.
\end{definition}

\begin{remark}
Haugseng and Gepner in \cite{GepHauEnriched} call $PreC(\V)$ the \textit{$\infty$-categories of categorical algebras} and denote it by $Alg_{cat}(\V)$. Actually, $Alg_{cat}(\V)$ and $Quiv(\V)$ are not the same $\infty$-category, but A.W. Macpherson in \cite[Theorem 1.1]{Macpherson} proves that they are equivalent.
\end{remark}

In \cite[Lemma 4.3.9 and Corollary 4.3.16]{GepHauEnriched}, the authors prove that the construction of $PreC(\V)$ is functorial in $\V$. Indeed, it defines two $\infty$-functors:

\begin{equation}
\label{eqVaccaLogiaPreC}    
PreC(-): PMon^{Lax}\to \Cat;
\end{equation}

and

\begin{equation}
\label{eqVaccaLogia1PreCat}    
PreC(-): Mon^{Pr}\to Pr^{L}.  
\end{equation}

\begin{notation}
\label{notprec2}
    In this paper, we use the same notation for both the $\infty$-functors (\ref{eqVaccaLogiaPreC}) and (\ref{eqVaccaLogia1PreCat}).
    Let $f:\A^{\otimes}\to \B^{\otimes}$ be an arrow of $PMon^{Lax}$ or $Mon^{pr}$, we denote by $f_{!}$ the 
\[f_{!}:=PreC(f):PreC(\A)\to PreC(\B). \]

\end{notation}

\begin{notation}[{\cite[Lemma 4.3.19]{GepHauEnriched}}]
\label{notPreC} 
 Let $f:\A^{\otimes}\to \B^{\otimes}$ be a monoidal $\infty$-functor between presentably monoidal $\infty$-categories $\A^{\otimes}$ and $\B^{\otimes}$, such that its underlying $\infty$-functor $f:\A\to \B$ preserves small colimits and has a right adjoint $f\dashv g:\B\to \A$ (this implies that $g$ is lax monoidal). Then there exists an adjunction 
\[f_{!}\dashv g_{!}:Prec(\B)\to Prec(\A) \] in the slice $\infty$-cosmos $\mathbf{QCat}_{/\mathcal{S}}$. Note that the right $g_{!}$ is not necessarily an arrow in $Pr^{L}$.

In particular, for each space $X$, we have a specific case of the above: 
\[f_!:=PreC(id_{X},f):PreC_{X}(\A)\to PreC_{X}(\B)\]
and an adjunction 
\[f_{!}\dashv g_{!}:Prec_X(\B)\to Prec_X(\A). \]
We use the same notation for the specific and the general case.
\end{notation}

\begin{remark}
In \cite{GepHauEnriched}, the authors adopt the notation $f_*$ instead of our notation $f_{!}$.
\end{remark}

\begin{example}
\label{exFreePrecUnderSeg}
Let $(\mathcal{V},\otimes,\I_{\mathcal{V}})$ be a presentably monoidal $\infty$-category, and let $\mathcal{S}^\times:=(\mathcal{S},\times, \mathbb{1})$ be the cartesian monoidal structure of the $\infty$-category of spaces.

Let $\I_{\mathcal{V}}:\mathbb{1}\to \mathcal{V}$ be the unit object of $\mathcal{V}$. Taking the left Kan extension of $\I_{\mathcal{V}}$ along the Yoneda embedding $\Yo:\mathbb{1}\to P(\mathbb{1})\simeq \mathcal{S}$ we obtain a small-colimits-preserving $\infty$-functor, that we denote by $-\otimes \I_{\mathcal{V}}:\mathcal{S}\to \mathcal{V}$. The notation $-\otimes \I_{\mathcal{V}}$ is not casual but is because every presentable $\infty$-category $\mathcal{V}$ is $\mathcal{S}$-left tensored $\infty$-category (\cite[Remark 5.5.1.7]{HTT}) and this $\infty$-functor is the tensor with $\I_{\mathcal{V}}$. This functor has as a right adjoint the corepresented $\infty$-functor $\Hom_{\mathcal{V}}(\I_{\mathcal{V}},-):\mathcal{V}\to \mathcal{S}$. Now we prove that $-\otimes\I_{\mathcal{V}}$ is a  monoidal $\infty$-functor.

Let $X$ be a space. Using the formula of left Kan extension, we obtain the equivalence

\begin{equation}
\label{eqMariano}   
-\otimes \I_{\mathcal{V}}(X)\simeq \mathrm{colim}(d:= \Hom_{\mathcal{S}}(\Yo, X)\xrightarrow{p_0} \mathbb{1}\to \mathcal{V}).
\end{equation}

Remarking that $\Hom_{\mathcal{S}}(\Yo, X)\simeq \Hom_{\mathcal{S}}(\mathbb{1}, X)\simeq X$, we obtain that the above diagram $d$ is the constant $\infty$-functor to $\I_{\mathcal{V}}$ with source $X$. So we can change the right-hand part in \eqref{eqMariano} to \[-\otimes \I_{\mathcal{V}}(X)\simeq \coprod_{\pi_{0}X}\I_{\mathcal{V}}. \]

Let $Y$ be another space. Since $-\otimes^{L}-$ preserves small colimits componentwise \cite[Remark 4.8.1.18]{HA}, $\I_{\mathcal{V}}$ is the unit of $-\otimes_{\V}-$ and that component-of-space $\infty$-functor $\pi_0$ preserves finite products, the following chain of equivalence holds: 

\begin{equation*}
\begin{split}
(-\otimes \I_{\V}(X))\otimes (-\otimes \I_{\V}(Y))\simeq & \coprod_{\pi_{0}X}\I_{\V}\otimes \coprod_{\pi_{0}Y}\I_{\V} \\ \simeq & \coprod_{\pi_{0}Y}(\coprod_{\pi_{0}X}\I_{\V}\otimes\I_{\V})
\simeq  \coprod_{\pi_{0}Y}\coprod_{\pi_{0}X}(\I_{\V}\otimes\I_{\V}) \\ \simeq & \coprod_{\pi_{0}Y\times\pi_{0}X}\I_{\V}\otimes\I_{\V}\simeq \coprod_{\pi_{0}(Y\times X)}\I_{\V} \\ \simeq & -\otimes \I_{\V}(X\times Y).
\end{split}
\end{equation*}

This proves that $-\otimes \I_{\mathcal{V}}$ is a monoidal $\infty$-functor.

In conclusion, we have an adjunction

\begin{equation}
\label{eqPrecatUnder}
(-\otimes \I_{\mathcal{V}})_!\dashv Hom_{\mathcal{V}}(\I_{\mathcal{V}},-)_{!}:PreC(\mathcal{V})\to PreC(\mathcal{S}).
\end{equation}
\end{example}

\begin{remark}    
In \cite[Example 4.3.20]{GepHauEnriched}, the authors define $-\otimes \I_{\mathcal{V}}$ with a different technique, and in their case, $-\otimes \I_{\mathcal{V}}$ is monoidal by construction.
\end{remark}

Roughly speaking, an associative algebra in $Quiv_{X}(\mathcal{V})$ $\mathcal{C}$ is a pair \[(\mathcal{C}:X\times X^{op}\to \mathcal{V},m:\mathcal{C}\otimes_{a}\mathcal{C}\to \mathcal{C})\in Fun(X\times X^{op},\mathcal{V})\times Fun(X\times X^{op},\mathcal{V})^{\mathbb{2}}, \] whose the first component is an $\infty$-functor which is the hom-$\mathcal{V}$-object, the second component is a natural transformation between $\infty$-functors which contains all the information about composition. 

Indeed, for each triple of objects $x,y,z   \in X$ in the space of objects of $\mathcal{C}$, we can find the composition morphism in $\mathcal{V}$ as the following morphism in $\mathcal{V}$:

\[ \mathcal{C}(x,y)\otimes_{\mathcal{V}}\mathcal{C}(y,z)\xrightarrow{\lambda_{(y,y)} } \oint^{w\in X} \mathcal{C}(x,w)\otimes_{\mathcal{V}}\mathcal{C}(w,z)= \mathcal{C}\otimes_{a}\mathcal{C} (x,y)\xrightarrow{m_{x,y}}\mathcal{C}(x,y),     \]
here, $\lambda_{(y,y)}$ is the $(y,y)$-component of the colimit cone of $\oint^{w\in X} \mathcal{C}(x,w)\otimes_{\mathcal{V}}\mathcal{C}(w,z)$ (every coend is a colimit).

Moreover, the axioms that define an enriched category are all a consequence of the dinaturality of the coend. 

Actually, every $\mathcal{V}$-enriched precategory has everything we expect from a $\mathcal{V}$-enriched $\infty$-category, but we already mentioned this is not the correct concept.

The problem is that the $\infty$-category $PreC(\mathcal{V})$ does not have as equivalence the correct $\infty$-categorical analogous of the Dwyer-Kan equivalence. We will localize $PreC(\mathcal{V})$ to obtain the correct equivalences.

Haugseng, Gepner, and Hinich prove that this problem is not only similar to that of Segal Space, as mentioned in the introduction, but it is the same.

Indeed, there exists an equivalence of $\infty$-categories, see \cite[Theorem 4.4.7]{GepHauEnriched} or \cite[\S 5.6]{HinichYon}, 
\begin{equation}
\label{eqPrecSegSpace}     
PreC(\mathcal{S}^\times)\simeq Seg(\mathcal{S})
 \end{equation}
where $Seg(\mathcal{S})$ is the $\infty$-category underlying Rezk's model category of Segal Spaces (see \cite[Definition 1.3.1]{HinichDK},\cite[Definition 3.4]{DoniPhDThesis} and \cite[\S 7.1]{Rezk}) and $\mathcal{S}^{\times}$ is the monoidal category of spaces whose tensor is the product and whose unit is the terminal space $*$.
\begin{notation}
The $\infty$-functor $Hom_{\mathcal{V}}(\I_{\mathcal{V}},-)_{!}$ in \eqref{eqPrecatUnder} sends every $\mathcal{V}$-enriched precategory $\A$ to a Segal Space, which we call the \textit{underlying Segal space of $\A$}.
\end{notation}

Finally, we can define what a $\mathcal{V}$-enriched $\infty$-category is.

\begin{definition}
    \label{defHigherenrichedCategory}
Let $X$ be a space, $(\mathcal{V},\otimes_{\mathcal{V}},\I_{\mathcal{V}})$ be a presentably monoidal $\infty$-category. An $\mathcal{V}$\textit{-enriched category $\A$ with space of objects $X$} is a $\mathcal{V}$-enriched precategory such that its underlying Segal space is complete. 

Moreover, we define the \textit{$\infty$-category of $\mathcal{V}$-enriched $\infty$-categories}, denoted $\Cat^{\mathcal{V}}$, as the full $\infty$-subcategory of $PreC(\mathcal{V})$ spanned by $\mathcal{V}$-enriched categories.
 
\end{definition}

\begin{notation}
Let $\mathcal{A}$ be a $\mathcal{V}$-enriched $\infty$-category. We will use the notation $\mathcal{A}_{ob}$ to denote the space of elements of $\mathcal{A}$. 
\end{notation}

\begin{remark}
Using \Cref{defHigherenrichedCategory} with $\mathcal{V}=(\mathcal{S},\times, *)$ be the $\infty$-category of spaces with the cartesian structure, and \eqref{eqPrecSegSpace}, we obtain tautologically the equivalence \[\Cat^{\mathcal{S}}\simeq\Cat.\]

From now on, we use both notations $\Cat^{\mathcal{S}}$ and $\Cat$ for the $\infty$-category of small $\infty$-categories.
\end{remark}

In \cite[Proposition 5.4.3]{GepHauEnriched}, the authors prove that, for each $\V$ presentably monoidal $\infty$-category, there is a localization

% https://q.uiver.app/#q=WzAsMixbMCwwLCJQcmVDKFxcVikiXSxbMiwwLCIkQ2F0XntcXFZ9JCJdLFswLDEsIiIsMCx7Im9mZnNldCI6Mn1dLFsxLDAsImkiLDIseyJvZmZzZXQiOjMsInN0eWxlIjp7InRhaWwiOnsibmFtZSI6Imhvb2siLCJzaWRlIjoidG9wIn19fV1d
\begin{equation}
\label{eqLocPreInEnr}
\begin{tikzcd}
	{PreC(\V)} && {\Cat^{\V},}
	\arrow["L"',shift right=2, from=1-1, to=1-3]
	\arrow["i"', shift right=3, hook, from=1-3, to=1-1]
\end{tikzcd}
\end{equation}

which exhibits $\Cat^{\V}$ as the localization of $PreC(\V)$ with respect the essentially surjective fully faithful $\V$-enriched $\infty$-functor, see \Cref{defFullyFaith} and \Cref{defEssSurj}. In particular $\Cat^{\V}$ is a presentable $\infty$-category.

The same constructions about enriched precategory hold for enriched $\infty$-category. It has the same functoriality and the following $\infty$-functors are well-defined:

\begin{equation}
\label{eqVaccaLogiaCAt}    
Cat_{\infty,(-)}(-):=Alg_{Ass}(Quiv): \mathcal{S}^{op} \times PMon^{Lax}\to \Cat; 
\end{equation}

\begin{equation}
\label{eqVaccaLogia1CAt}    
Cat_{\infty,(-)}(-):=Alg_{Ass}(Quiv): \mathcal{S}^{op} \times Mon^{Pr}\to {Pr^{L}};  
\end{equation}

\begin{equation}
\label{eqFunCatErFanculo}
\Cat^{(-)}: PMon^{Lax}\to  \Cat;  
\end{equation}

and

\begin{equation}
\label{eqFunCatEr}
\Cat^{(-)}: Mon^{Pr}\to  Pr^{L}.  
\end{equation}
In \cite[\S 5.7]{GepHauEnriched}, the authors prove that for $\Cat^{(-)}$ has the same properties as $PreC(-)$. In fact, we can rewrite \Cref{notPreC} using $Cat^{(-)}$ instead of $PreC(-)$.

\begin{notation}
\label{notprec3}
    In this paper, we use the same notation for both the $\infty$-functors \eqref{eqVaccaLogiaCAt} and \eqref{eqVaccaLogia1CAt}, and for \eqref{eqFunCatErFanculo}, \eqref{eqFunCatEr}.
    Let $f:\A^{\otimes}\to \B^{\otimes}$ be an arrow of $PMon^{Lax}$ or $Mon^{pr}$, we denote by $f_{!}$ the $\infty$-functor 
\[f_{!}:=\Cat^{f}:\Cat^{\A}\to \Cat^{\B}. \]  
and we call $f_{!}$ \textit{the change of base via $f$}.
\end{notation}

\begin{notation}[{\cite[Proposition 5.7.17]{GepHauEnriched}}]
\label{notCatPuntoEsc}
 Let $f:\A^{\otimes}\to \B^{\otimes}$ be a monoidal $\infty$-functor between presentably monoidal $\infty$-categories $\A^{\otimes}$ and $\B^{\otimes}$, such that its underlying $\infty$-functor $f:\A\to \B$ preserves small colimits and has a right adjoint $f\dashv g:\B\to \A$ (this implies that $g$ is lax monoidal).
Then there exists an adjunction 
\[f_{!} \dashv g_{!}: \Cat^{\B} \to \Cat^{\A}, \] in the slice $\infty$-cosmos $\mathbf{QCat}_{/\mathcal{S}}$. Note that the right adjoint $\infty$-functor is not necessarily an arrow in $Pr^{L}$.  

In particular, for each space $X$, we have the specific case: 
\[f_! := \Cat^{(-)}(id_{X}, f): Cat^{\A}_{\infty,X} \to Cat^{\B}_{\infty,X},\]
and an adjunction
\[f_{!} \dashv g_{!}: Cat^{\B}_{\infty,X} \to Cat^{\A}_{\infty,X} \] in $\Cat$.
Again, $f_{!}$ belongs to $Pr^{L}$ but $g_{!}$ does non-necessarily need to be.

We use the same notation for both specific and general cases.
    
\end{notation}

\begin{remark}
In \cite{GepHauEnriched}, the authors use the notation $f_*$ instead of our notation $f_{!}$.
\end{remark}

\begin{example}[Free $\V$-enriched $\infty$-category and underlying $\infty$-category]
\label{RmkFreeUnderlyingEnrichedCatAdjunction}

Let $\V$ be a presentably monoidal $\infty$-category. In \Cref{exFreePrecUnderSeg} we found an adjunction \[-\otimes\I_{\V}\dashv Hom_{S}(\I_{\V}, -):\V\to \mathcal{S},\] that satisfies the condition in \Cref{notCatPuntoEsc}.

Hence there exists an adjunction
\begin{equation}
(-\otimes \I_{\V})_!\dashv \Hom_{\V}(\I_{\V},-)_{!}:\Cat^{\V}\to \Cat^{\mathcal{S}}.
\end{equation}

We call the \textit{underlying $\infty$-category $\infty$-functor} the following composition of $\infty$-functors:
\begin{equation}
% \label{eqHomotopyHigherFunctorCat}  
(-)_{o}:\Cat^{\V}\xrightarrow{\Hom(\I_{\V},-)_{!}}\Cat^{\Top}\simeq \Cat:\;
 \A\mapsto \A_{o}.
\end{equation}
Let $\A$ be an enriched $\infty$-category we call $\A_o$ its underlying $\infty$-category. In particular, for each $x,y$ elements of $\A$, the space $\A_{o}(x,y)$ is the mapping space $\Hom_{\A_{o}}(x,y)$ of the $\infty$-category $\A_o$.  
%Moreover, we denote by $hA$ its homotopy category $h\A_{o}$ and call it \textit{the homotopy category of $\A$}. 
%Let $f:\A\to\B$ be an object in $\Hom_{\Cat^{\V}}(\A,\B)$, we call $ f_o:\A_{o}\to\B_{o}$ the \textit{underlying $\infty$-functor of $f$}.

%Let $\A$ be a $\V$-enriched $\infty$-category; we call $Hom_{\V}(\I_{\V},-)_{!}(\A)$ the \textit{underlying $\infty$-category of $\A$} and denote it by $A_o$.
 
 Let $f:\A\to\B$ be an object in $\Hom_{\Cat^{\V}}(\A,\B)$; we call $f_o:\A_{o}\to \B_{o}$ the \textit{underlying $\infty$-functor of $f$}.% and we denote it by $ f_o:\A_{o}\to\B_{o}$.

Let $\C$ be an $\infty$-category; we call $(-\otimes \I_{\V})_!(\C)$ the \textit{associated $\V$-enriched $\infty$-category of $\C$} and denote it by the same $\C$.% This example has a categorical analog, see \cite[Example 3.3.4, Definition 3.4.5]{riehl2014categorical}.
\end{example}

\begin{notation}
\label{notHomotopyCategoryOfA}
Let $\A$ be an $\V$-enriched $\infty$-category. We call \textit{homotopy category of $\A$} and denote it by $h\A$ the homotopy category of the underlying $\infty$-category of $\A$, i.e. $h\A:= h\A_{o}$.
Furthermore, we define the \textit{homotopy $\infty$-functor} as the following $\infty$-functor:
\begin{equation}
% \label{eqHomotopyHigherFunctorCat}  
h:\Cat^{\V}\xrightarrow{Hom(\I_{\V},-)_{!}}\Cat^{\Top}\simeq\Cat\xrightarrow{h} Cat. 
\end{equation}
\end{notation}

Let us finish this section with an example of monoidal and lax $\infty$-functors. Before, we fix a notation.

\begin{example}
\label{exAdjunctionEnrichedVLMODRV}
Let $\V$ be an $\mathbb{E}_2$-monoidal $\infty$-category and let $R$ be an $\mathbb{E}_1$-ring of $\V$.
Then the following adjunction exists %\textcolor{red}{Qua devo citare qualcosa per dire che esiste la seguente aggiunzione}

\begin{equation}
\label{eqMannaggiaStronzaGenerale1}
-\otimes_{\V} R : \V\simeq \mathrm{LMod}_{\mathbb{\mathbb{I}_{\V}}}(\V)\rightleftarrows \mathrm{LMod}_{R}(\V):U. 
\end{equation}

In particular, this implies that the right adjoint $U_{!}$ is lax monoidal.
We are in the situation of \Cref{notCatPuntoEsc}, so \eqref{eqMannaggiaStronzaGenerale1} induces the following adjunction between enriched $\infty$-categories
\begin{equation}
\label{eqMannaggiaStronzaGenerale}
-\otimes_{\V}R_{!}:\Cat^{\V}\rightleftarrows \Cat^{\mathrm{LMod}_{R}(\V)}: U_{!}. 
\end{equation}
   
\end{example}

\begin{example}
\label{exAdjunctionEnrichedSpLMod}
In \Cref{rmkAdjunctionHigherCatextentionrestrictioncorestriction}, we recall that there is the extension-restriction adjunction
\begin{equation}
\label{eqMannaggiaStronza}
-\otimes_{\Sp} \Hk : \mathrm{LMod}_{\mathbb{S}}\rightleftarrows \LModk:U. 
\end{equation}

Now we prove that the left adjoint $-\otimes_{\Sp} \Hk$ is monoidal. 

Recall that the smash product of $\Sp$ is the same as the tensor product $_{\mathbb{S}}\otimes_{\mathbb{S}}$ in $\mathrm{LMod}_{\mathbb{S}}$, see \Cref{exsmashproduct}.

For each $X,Y$ elements of $\Sp$, the following list of equivalences holds:
\begin{equation*}    
\begin{split}
\otimes\Hk(X)\; \otimes_{\Hk} \otimes\Hk(Y) \simeq &  
(X\otimes_{\Sp}\Hk)\; _{\Hk}\otimes_{\Hk} (Y\otimes_{\Sp}\Hk)  \\ \simeq & (X\; _{\mathbb{S}}\otimes_{\mathbb{S}}\Hk)\; _{\Hk}\otimes_{\Hk} (Y\; _{\mathbb{S}}\otimes_{\mathbb{S}}\Hk) \\ \simeq & X\; _{\mathbb{S}}\otimes_{\mathbb{S}}(\Hk\; _{\Hk}\otimes_{\Hk} \Hk)\; _{\mathbb{S}}\otimes_{\mathbb{S}}Y\\ \simeq & X\; _{\mathbb{S}}\otimes_{\mathbb{S}} \Hk\; _{\mathbb{S}}\otimes_{\mathbb{S}}Y \\ \simeq & X\; _{\mathbb{S}}\otimes_{\mathbb{S}} Y\; _{\mathbb{S}}\otimes_{\mathbb{S}}\Hk \\ \simeq & \otimes\Hk(X\otimes_{\Sp} Y).
\end{split}
\end{equation*}
The above equivalences hold because $\Sp$ and $\LModk$ are symmetric monoidal $\infty$-categories, the relative tensor product is associative (\cite[Proposition 4.4.3.14]{HA}), and because there is a canonical equivalence  \[\Hk\;_{\Hk}\otimes_{\Hk}\Hk=Bar_{\Hk}(\Hk,\Hk)\simeq \Hk .\]

In particular, this implies that the right adjoint $U_{!}$ is lax monoidal.
We are in the situation of \Cref{notCatPuntoEsc}, so \eqref{eqMannaggiaStronza} induces the following adjunction between enriched $\infty$-categories
\begin{equation}
-\otimes_{\Sp}\Hk_{!}:\Cat^{\Sp}\rightleftarrows \Cat^{\LModk}: U_{!}. 
\end{equation}
   
\end{example}

\subsection{Tensored Enriched categories}
\label{subsectensoredEnrichedCat}

In this paper, we are interested in the tensored enriched $\infty$-category, and in this subsection, we will explain what this means for us. This definition is expected, but it needs some clarification, which we will include in \Cref{rmkDefEnrichedCat}, which will allow us to give a second definition, \Cref{defTensoredEnrichedCategory1}.
This second definition allows us to define an enriched version of the definition.

Let $\V$ be a presentably monidal $\infty$-category and 
let $\A$ be an $\V$-enriched $\infty$-category.
We start with the following canonical $\infty$-functor:
\begin{equation}
    \label{eqHom(A())}
    \Hom_{\V}(-,\A(-,-)):\V^{op}\times\A_o^{op}\times \A_o\xrightarrow{id_{\V}\times \A(-,-)}\V^{op}\times\V\xrightarrow{\Hom_{\V}(-,-)} \Top,
\end{equation}
and then we apply the Grothendieck construction (for a good presentation of the Grothendieck construction, see \cite[\S 2.2.6]{HinichYon}) and we obtain a module \[\int\Hom_{\V}(-,\A(-,-)) \to (\V\times\A_o)\times \A_o.\]

\begin{warning}
Two clarifications regarding the previous construction
\begin{itemize}
    \item $\int\Hom_{\V}(-,\A(-,-))$ is a module but not necessarily a comma $\infty$-category; 
    \item a bifibration according to Hinich is different from a bifibration as defined by Verity-Riehl. Indeed in \cite{RiehlElements}, the authors call modules what Hinich calls bifibration. In this work, we follow the Verity-Riehl convention.
\end{itemize}
\end{warning}

\begin{remark}[{Technical point}]
\label{rmkDefEnrichedCat}
In our definition of $\V$-enriched $\infty$-category, a $\V$-enriched $\infty$-category $\A$ is an $\infty$-functor  \[\A(-,-):\A_{ob}^{op}\times\A_{ob}\to \V \] with other external properties. That is a problem because this implies that the $\infty$-functor $\Hom_{\A_o}(-,\A(-,-))$ in \eqref{eqHom(A())} is not well defined:
it is defined on the product of the underlying $\infty$-category of $\A$, i.e.  $\A_o^{op}\times \A_o$, while $\A(-,-)$ is defined only on the product of the space of objects of $\A$, i.e. $\A_{ob}^{op}\times\A_{ob}$.
This problem is linked to our choice of definition of $\V$-enriched $\infty$-category. Indeed if we had chosen Lurie's definition, an $\V$-enriched $\infty$-category $\A$ would have been an $\infty$-functor $\A(-,-):\A^{op}\times\A\to \V$ with extra properties, see  \cite[Remark 4.2.1.31]{HA}.   
Since H.Heine in \cite{Heine} has proved that the two notions of $\V$-enriched $\infty$-category are equivalent, the $\infty$-functor in (\ref{eqHom(A())}) is well defined.
\end{remark}

\begin{definition}
\label{defTensoredEnrichedCategory}
Let $\A$ be an $\V$-enriched $\infty$-category, we say that it is \textit{tensored} if the module $\int\Hom_{\V}(-,\A(-,-))$ is left representable: there exists an $\infty$-functor $-\otimes -: \V\times\A_o\to \A_o$ and an equivalence of modules  \[\int\Hom_{\V}(-,\A(-,-))\simeq_{\V\times\A_o\times \A_o} \Hom_{\A_o}(-\otimes-,\A_o).\]

\end{definition}

Since we work in the $\infty$-cosmos of $(\infty,1)$-categories, equivalences of modules and fiberwise equivalences are the same concept, \Cref{rmkfiberequivalencesFiberwise}, we can redefine the above definition as follows.

\begin{definition}
\label{defTensoredEnrichedCategory1}
Let $\A$ be a $\V$-enriched $\infty$-category. It is \textit{tensored} if there exists an $\infty$-functor $-\otimes -: \V\times\A_o\to\A_o:(v,a)\mapsto v\otimes a$ such that for each  element $(v,n,m):\mathbb{1}\to \V\times\A_{o}\times\A_{o}$ of $\V\times\A_{o}\times\A_{o}$, there exists an element $v\otimes m:\mathbb{1}\to \A$ in the image of $-\otimes-$ together with an isomorphism of spaces
\begin{equation}
    \label{eqTensor}
    \Hom_{\V}(v,\A(m,n))\simeq\Hom_{\A_o}(v\otimes m,n). 
\end{equation}
\end{definition}

\begin{remark}    
\label{remarkLurieLeftTensorAndOur}
In particular, we have recovered the definition of tensored $\V$-enriched $\infty$-category of Lurie. However, Lurie calls left-tensored what we call tensored, see \cite[Definitions 4.2.1.19, 4.2.1.28]{HA} or our \Cref{exVLeftTensoredCategory}.
Roughly speaking, the difference is that Lurie starts with a left-tensored $\infty$-categories and then he defines the enrichment. Instead, we start from the enrichment and then we define the tensor.
We decided to use the word \textit{tensored} instead of \textit{left-tensored} because we want to be coherent with \cite[\S 3.7]{riehl2014categorical} and \cite[Definition 6.9]{heine2023equivalence}.
\end{remark}

\begin{example}
\label{exPresentablyMonoidalistensoredOnItself}
Every presentably monoidal $\infty$-category $\V$ may be seen as a tensored $\V$-enriched $\infty$-category. 

Indeed, it is $\V$-left tensored $\infty$-category with a tensor product that preserves componentwise the colimits, so the Adjoint Functor Theorem implies that it is a $\V$-enriched $\infty$-category. We conclude with \Cref{remarkLurieLeftTensorAndOur}. 
\end{example}

\begin{definition}
\label{defTensoredinEnrichedWay}
Let $\A$ be a $\V$-enriched $\infty$-category, it is \textit{tensored in $\V$-enriched way} if there exists an $\infty$-functor $-\otimes -: \V\times\A_o\to\A_o:(v,a)\mapsto v\otimes a$ such that for each  element $(v,n,m):\mathbb{1}\to \V\times\A_{o}\times\A_{o}$ of $\V\times\A_{o}\times\A_{o}$, there exists an element $v\otimes m:\mathbb{1}\to \A$ in the image of $-\otimes-$ together with an isomorphism in $\V$
\begin{equation}
    \label{eqTensorEnrichedWay}
    \V(v,\A(m,n))\simeq \A(v\otimes m,n). 
\end{equation}
\end{definition}

\begin{remark}
Applying the $\infty$-functor $\Hom_{\V}(\I_\V,-): \V\to \Top$ to the isomorphism in (\ref{eqTensorEnrichedWay}), we obtain the isomorphism in (\ref{eqTensor}), hence every $\infty$-category tensored in $\V$-enriched way is tensored.
\end{remark}

\begin{remark}
\Cref{defTensoredinEnrichedWay} is the enriched $\infty$-categorical analogous of \cite[Definition 3.7.2]{RiehlElements}.
We expect that the $\V$-isomorphisms (\ref{eqTensorEnrichedWay}) implies that, for each element $m:\mathbb{1}\to \A$, the adjunction $-\otimes m\dashv \A(m,-)$ is actually a $\V$-adjunction. In the ordinary case, it is easy to prove this fact using some properties of the tensor, see \cite[Proposition 3.7.10]{riehl2014categorical}. But we will not analyze this situation.
\end{remark}

\begin{proposition}
Let $\V$ be a presentably monoidal $\infty$-category.

Let $\A$ be a $\V$-enriched $\infty$-category, it is tensored in $\V$-enriched way if and only if it is tensored.
\end{proposition}

\begin{proof}
We need to prove that for each $v:\mathbb{1}\to \V$ element of $\V$ and each pair of elements of $\A$ $n:\mathbb{1}\to \A$, $m:\mathbb{1}\to \A $, the elements of $\V$ $\V(v,\A(m,n))$ and  $\A(v\otimes m,n)$ are equivalent. To prove that we show that they represent the same $\infty$-functor $\V\to \Top$.
Let $w:\mathbb{1}\to \V$ be an element of $\V$, the following equivalence in $\V$ holds: 

\begin{equation}
\begin{split}
\Hom_{\V}(w,\A(v\otimes m,n)) & \simeq \Hom_{\A_{o}}(w\otimes (v\otimes m),n) \\ 
& \simeq \Hom_{\A_o}((w\otimes_{\V}v)\otimes m,n) \\ 
& \simeq \Hom_{\V}(w\otimes v, \A(m,n)) \\ 
& \simeq \Hom_{\V}(w,\V(\A(v\otimes m))). 
\end{split}    
\end{equation}

The first three equivalences hold because $\A$ is tensored.
In \Cref{exPresentablyMonoidalistensoredOnItself} we showed that a $\V$ is a tensored $\V$-enriched $\infty$-category, so the last equivalence holds.
\end{proof}

\begin{remark}
\label{rmkTensoredcatForMeAndHeine}
The last proof is important because it implies that when $\V$ is a presentably monoidal $\infty$-category our definition of tensored $\V$-enriched $\infty$-categories and Heine's definition of $\V$-tensored $\infty$-categories, \cite[Definition 6.9]{heine2023equivalence}, coincide.

Moreover, in this paper, we work only with $\infty$-categories enriched in presentably monoidal $\infty$-categories, so the three definitions in this subsection are always equivalent to us.
\end{remark}

Similarly, we can define as above when a $\V$-enriched $\infty$-category is cotensored.

\begin{definition}
Let $\A$ be an $\V$-enriched $\infty$-category, we say that it is \textit{cotensored} if the module $\int\Hom_{\V}(-,\A(-,-))$ is right representable:
 if there exists an $\infty$-functor \[(-)^ {(-)}: \V^{op}\times\A_o\to\A_o:(v,a)\mapsto v\otimes a\] such that for each  element $(v,n,m):\mathbb{1}\to \V\times\A_{o}\times\A_{o}$ of $\V\times\A_{o}\times\A_{o}$, there exists an element $ n^{v}$ in the image of $(-)^{(-)}$  together with an isomorphism of spaces
\begin{equation}
    \label{eqCoTensor}
    \Hom_{\V}(v,\A(m,n))\simeq\Hom_{\A_o}(m,n^v). 
\end{equation}

Moreover, we say that $\A$ is \textit{cotensored in $\V$-enriched way} if there exists an $\infty$-functor $(-)^ {(-)}: \V^{op}\times\A_o\to\A_o:(v,a)\mapsto v\otimes a$ such that for each  element $(v,n,m):\mathbb{1}\to \V\times\A_{o}\times\A_{o}$ of $\V\times\A_{o}\times\A_{o}$, there exists an element $ n^{v}$ in the image of $(-)^{(-)}$  together with an isomorphism of spaces together with an isomorphism in $\V$
\begin{equation}
    \label{eqcoTensorEnrichedWay}
    \V(v,\A(m,n))\simeq \A(m,n^v). 
\end{equation}
\end{definition}

%%%%%%%%%%%%%%%%%%%%%%%%%%%%%%%%%%%%%%%%%%%%%%%%%%%%%%%%%%%%%%%%%%%%%%%%%%%%%%%%%%%%%%%%%%%%%%%%%%%%%%%%%%%%%%%%%%%%%%%%%%%%%%%%%%%%%%%%%%%%%%%%%%%%%%%%%%%%%%%%%%%%%%%%%%%%%%%%%%%%%%%%%%%%%%%%%%%%%%

%%%%%%%%%%%%%%%%%%%%%%%%%%%%%%%%%%%%%%%%%%%%%%%%%%%%%%%%%%%%%%%%%%%%%%%%%%%%%%%%%%%%%%%%%%%%%%%%%%%%%%%%%%%%%%%%%%%%%%%%%%%%%%%%%%%%%%%%%%%%%%%%%%%%%%%%%%%%%%%%%%%%%%%%%%%%%%%%%%%%%%%%%%%%%%%%%%%%%%%%%%%%%%%%%%%%%%%%%%%%%%%%%%%%%%%%%%%%%%%%%%%%%%%%%%%%%%%%%

%%%%%%%%%%%%%%%%%%%%%%%%%%%%%%%%%%%%%%%%%%%%%%%%%%%%%%%%%%%%%%%%%%%%%%%%%%%%%%%%%%%%%%%%%%%%%%%%%%%%%%%%%%%%%%%%%%%%%%%%%%%%%%%%%%%%%%%%%%%%%%%%%%%%%%%%%%%%%%%%%%%%%%%%%%%%%%%%%%%%%%%%%%%%%%%%%%%%%%%%%%%%%%%%%%%%%%%%%%%%%%%%%%%%%%%%%%%%%%%%%%%%%%%%%%%%%%%%%

%%%%%%%%%%%%%%%%%%%%%%%%%%%%%%%%%%%%%%%%%%%%%%%%%%%%%%%%%%%%%%%%%%%%%%%%%%%%%%%%%%%%%%%%%%%%%%%%%%%%%%%%%%%%%%%%%%%%%%%%%%%%%%%%%%%%%%%%%%%%%%%%%%%%%%%%%%%%%%%%%%%%%%%%%%%%%%%%%%%%%%%%%%%%%%%

\subsection{Monoidal structure in $\Cat^{\V}$}
\label{sussecMonoidalEnrichedCat}

 For generality, we write this susection with $(\V, \otimes_{\V}, \I_{\V})$ a presentably monoidal $\infty$-category and $R$ a commutative algebra of $\V$. However, if the reader is only interested in our results, he can replace $\V$ and $R$ with $\Sp$ and $\Hk$ respectively.

This section only provides an overview of how the tensor between enriched $\infty$-categories works and is far from being technical. For a detailed presentation, see \cite{haugseng2023tensor}, \cite{HinichYon}.

Firstly, we need to define an $\infty$-functor.

We begin with the canonical cocartesian fibration $\Cat^{\V} \to \Top$ and take the pullback with respect to the terminal object $\{*\}: \mathbb{1} \to \Top$:

% https://q.uiver.app/#q=WzAsNCxbMCwwLCJDYXRfe1xcaW5mdHksKn1ee1xcVn0iXSxbMSwwLCJcXENhdF57XFxWfSJdLFswLDEsIlxceypcXH0iXSxbMSwxLCJcXFRvcCJdLFswLDEsInAiXSxbMCwyXSxbMiwzXSxbMSwzXSxbMCwzLCIiLDEseyJzdHlsZSI6eyJuYW1lIjoiY29ybmVyIn19XV0=
\begin{equation}
\label{eqMikiMouse}    
\begin{tikzcd}[ampersand replacement=\&]
	{Cat_{\infty,*}^{\V}} \& {\Cat^{\V}} \\
	{\mathbb{1}} \& {\Top,}
	\arrow["p", from=1-1, to=1-2]
	\arrow[from=1-1, to=2-1]
	\arrow["\{*\}",from=2-1, to=2-2]
	\arrow[from=1-2, to=2-2]
	\arrow["\lrcorner"{anchor=center, pos=0.125}, draw=none, from=1-1, to=2-2]
\end{tikzcd}
\end{equation}
the apex of the pullback is the $\infty$-category of $\V$-enriched $\infty$-categories with a contractible space of objects $Cat_{\infty,*}^{\V}$. Tautologically, a $\V$-enriched $\infty$-category with a contractible space of objects is an associative algebra in $\V$, so there is a canonical equivalence $\phi: Alg(\V) \simeq Cat_{\infty,*}^{\V}$.   

\begin{definition}
\label{defunderlineHigher}    
We call the \textit{underline $\infty$-functor} the composition of the $\infty$-functor $\phi: Alg(\V) \simeq Cat_{\infty,*}^{\V}$ and the projection $p: Cat^{\V}_{*} \to \Cat^{\V}$ in (\ref{eqMikiMouse}):
\begin{equation}
\label{eqUnderlineFunctorGeneral}
\underline{(-)}: Alg_{Ass}(\V) \simeq Cat^{\V}_{\infty,*} \to \Cat^{\V}: \A \to \underline{\A}.
\end{equation}
\end{definition}

The above $\infty$-functor sends an algebra object in $\V$ to the $\V$-enriched $\infty$-category with spaces of objects being the terminal space $*$.

In \cite{haugseng2023tensor}, the authors prove that when $(\V, \otimes_{\V}, \I_{\V})$ is a presentably monoidal $\infty$-category, the $\infty$-category $\Cat^{\V}$ has a monoidal structure $(\Cat^{\V}, \otimes^{\V}, \I)$ whose unit is $\I := \underline{\I_{\V}}$ and $\Cat^{\V}$ is presentably monoidal.

\begin{notation}
\label{notServePerIntro}
We denote the tensor between enriched $\infty$-categories by  
\[
-\otimes^{\V}-: \Cat^{\V} \times \Cat^{\V} \to \Cat^{\V}: (\A, \B) \mapsto A\otimes^{\V}\B.
\]
\end{notation}

We can say more about this. Indeed, this monoidal structure behaves exactly as in the categorical case (see \cite[\S 7.3]{riehl2014categorical} for the ordinary case).

In \Cref{exMonoidalStructureAlgSpCAlgSp}, we recalled that both $CAlg(\V)$ and $Alg(\V)$ are presentably monoidal $\infty$-category, and the forgetful $\infty$-functors: 
\[Alg(\V)\to \V,\]
and
\[CAlg(\V)\to \V,\]
are monoidal.

The tensor of two $\V$-enriched $\infty$-categories both with trivial spaces of objects is still a $\V$-enriched $\infty$-category with trivial space of objects and the unit $\un{\I}_{\V}$ of the monoidal structure belong to $Cat^{\V}_{\infty,*}$, so the full $Cat^{\V}_{\infty,*}$ $\infty$-subcategory $\Cat^{\V}$ ineriths a structure of monoidal $\infty$-category (\cite[Remark 2.2.1.2]{HA}).  

So,  the equivalence $Alg(\V)\simeq Cat_{\infty,*}^{\V}$ is actually a monoidal equivalence.

We deduce that the underline $\infty$-functor $\underline{(-)}$ is monoidal, and then the $\infty$-categories $\Cat^{\V}$ is a left $Alg(\V)$-module object of $Pr^{L}$.

\begin{notation}
We denote by $^\V\otimes^{\V}$ the left action of $Alg(\V)$ on $\Cat^{\V}$
\begin{equation}
\label{eqAlgVtensoredCAtCatV}
-^{\V}\otimes^{\V} -: Alg(\V) \times \Cat^{\V} \xrightarrow{\underline{(-)}\times id} \Cat^{\V}\times\Cat^{\V} \xrightarrow{-\otimes^{\V} -} \Cat^{\V}.
\end{equation}

\end{notation}

Now, we want to define the $\infty$-categorical version of the double-underline functor founded in \cite{DoniCategorical}.

\begin{remark}
We tried to write the entire paper using only associative (or $\mathbb{E}_{1}$-algebras ) and commutative algebras (or $\E$-algebras ) of $\infty$-operad, but in the next part, it was not possible to avoid the $\mathbb{E}_2$-algebras. For an introduction to this topic, see \cite[\S 5]{HA}. 
The idea is that while in ordinary category theory, there is only the concept of associativity and commutativity monoids, in higher category there are \enquote{an infinite number of degrees of associativity}: associativity is the lowest degree and being commutative is the highest degree. 
The existence of infinite levels of associativity is the reason why we cannot define \Cref{defDoubleFunctorHigherCateg} without the use of $\mathbb{E}_{2}$: in particular that is why our fixed $k$ must be a commutative ring and not only a ring. Otherwise, we do not have the description of a $k$-linearization as the $\infty$-category of $\Cat^{\Sp}$-enriched, see \Cref{thAction}, because $\2un{\Hk}$ does not exist. 
In one other part of the paper, we shall use $\mathbb{E}_{2}$, \Cref{rmkHKAlgebraDiAlgebraDgCat}.
\end{remark}

Since the underline $\infty$-functor is monoidal, it induces an $\infty$-functor between algebras:
\begin{equation}
\label{eqAggiuntaPERChiarezza}
Alg(\underline{(-)}): Alg_{\mathbb{E}_{2}}(\V) \simeq Alg(Alg(\V)) \to Alg(\Cat^{\V}),
\end{equation}

where the isomorphism holds by Dunn Additivity Theorem, \cite[Theorem 5.1.2.2]{HA}.

There is an inclusion $j:\mathbb{E}_{2}\to Comm\simeq\E$. The precomposition with this inclusion defines an $\infty$-functor:
\begin{equation}
 \label{eqAggiuntaperComodità}   
 j^*: CAlg(\V)=Alg_{\E}(\V) \to Alg_{\mathbb{E}_2}(\V). 
\end{equation}

By composing $j^*$, \eqref{eqAggiuntaPERChiarezza} and the underline $\infty$-functor $\underline{(-)}:Alg(\Cat^{\V})\to\Cat^{\Cat^{\V}}$, we obtain the following definition:

\begin{definition}
\label{defDoubleFunctorHigherCateg}
We call the \textit{double-underline $\infty$-functor} the $\infty$-functor
\begin{equation}
\label{eqDoubleUnderlineGeneral}
\2un{(-)}:CAlg(\V) \xrightarrow{\un{(-)}\circ j^{*}} \Cat^{\Cat^{\V}}\simeq \mathcal{C}at_{(\infty,2)}^{\V}.
\end{equation}

\end{definition}

In \cite[\S 6.3]{GepHauEnriched}, the authors build a more general construction of the underline and double underline $\infty$-functors.

Let $n\in\N_{\geq 0}\cup \{\infty\}$ and let $\V$ be a $\mathbb{E}_{n+1}$-monoidal $\infty$-category, there exists a fully faithful and monoidal $\infty$-functors 
\begin{equation}
\label{eqBfunctor}
\mathcal{B}^{n}(-): Alg_{\mathbb{E}_n}(\V)\to \mathcal{C}at_{(\infty,n)}^{\V}: R\mapsto B^{n}
\end{equation}

The case $n=1$ and $n=2$ are equivalent to the $\infty$-functors \eqref{eqUnderlineFunctorGeneral} and \eqref{eqDoubleUnderlineGeneral}, the latter after precomposing them with $j^{*}$. 
To have notations consistent with the \cite{DoniklinearMorita} we retain the notation  $\underline{(-)}$ and $\underline{\underline{(-)}}$
for the special case just mentioned.

%%%%%%%%%%%%%%%%%%%%%%%%%%%%%%%%%%%%%%%%%%%%%%%%%%%%%%%%%%%%%%%%%%%%%%%%%%%%%%%%%%%%%%%%%%%%%%%%%%%%%%%%%%%%%%%%%%%%%%%%%%%%%%%%%%%%%%%%%%%%%%%%%%%%%%%%%%%%%%%%%%%%%%%%%%%%%%%%%%%%%%%%%%%%%%%%%%%%%%%%%%%%%%%%%%%%%%%%%%%%%%%%%%%%%%%%%%%%%%%%%%%%%%%%%%%%%%%%%%%%%%%%%%%%%%%%%%%%%%%%%%%%%%%%%%%%%%%%%%%%%%%%%%%%%%%%%%%%%%%%%%%%%%%%%%%%%%%%%%%%%%%%%%%%%%%%%%%%%%%%%%%%%%%%%%%%%%%%%%%%%%%%%%%%%%%%%%%%%%%%%%%%%%%%

%%%%%%%%%%%%%%%%%%%%%%%%%%%%%%%%%%%%%%%%%%%%%%%%%%%%%%%%%%%%%%%%%%%%%%%%%%%%%%%%%%%%%%%%%%%%%%%%%%%%%%%%%%%%%%%%%%%%%%%%%%%%%%%%%%%%%%%%%%%%%%%%%%%%%%%%%%%%%%%%%%%%%%%%%%%%%%%%%%%%%%%%%%%%%%%%%%%%%%%%%%%%%%%%%%%%%%%%%%%%%%%%%%%%%%%%%%%%%%%%%%%%%%%%%%%%%%%%%%%%%%%%%%%%%%%%%%%%%%%%%%%%%%%%%%%%%%%%%%%%%%%%%%%%%%%%%%%%%%%%%%%%%%%%%%%%%%%%%%%%%%%%%%%%%%%%%%%%%%%%%%%%%%%%%%%%%%%%%%%%%%%%%%%%%%%%%
\subsection{The $\infty$-category of $\V$-enriched $\infty$-functors }
\label{sebInfiniteEnrichedFunctor}

This subsection is far from technical, and for a detailed presentation, see \cite{HinichYon} and \cite{GepHauEnriched}.

In this section, we fix a presentably monoidal $\infty$-category $\V$.

This section can be divided into two parts. In the first one, we define what a $\V$-enriched $\infty$-functor is between two $\V$-enriched $\infty$-categories. Furthermore, we characterize when a $\V$-enriched $\infty$-functor is an equivalence; this part is taken from \cite{GepHauEnriched}.

In the second part, we will define the $\infty$-category of $\V$-enriched $\infty$-functors. To obtain this, we study the $\V$-enriched $\infty$-functors with a source a $\V$-enriched $\infty$-category and with a target a $\V$-enriched presentable $\infty$-category tensored in $\V$. This second part is taken from \cite{HinichYon}.

\begin{definition}
\label{defHigherEnrichedFunctor}
Let $\A$ and $\C$ be two enriched $\infty$-categories. \textit{A $\V$-enriched $\infty$-functor $F:\A\to \C$ with source $\A$ and target $\C$} is an object in the spaces $\Hom_{\Cat^{\V}}(\A,\C)$.

Moreover, we say that a $\V$-enriched $\infty$-functor $F$ is an \textit{equivalence} if it is an isomorphism as arrows of $\Cat^{\V}$.
\end{definition}

Our definition is categorical, it does not make it clear what a $\V$-enriched $\infty$-functor is. The equivalent Haugseng-Gepner's definition, see \cite[Definition 2.4.9]{GepHauEnriched}, is combinatorial and makes it possible to find analogies between a $\V$-enriched $\infty$-functor and an ordinal-enriched functor (\cite[Definition 2.3]{DoniCategorical}).

Indeed, starting from a $\V$-enriched $\infty$-functor $F$, it is possible to derive a morphism between the spaces of objects:
\begin{equation}
\label{eqMapBetweenObject}
    F_{ob}:\A_{ob}\to \C_{ob}:x\mapsto Fx ;
\end{equation}
and, for each pair of objects $x,y$ of $\A$, a morphism in $\V$:
\begin{equation}
\label{eqmapBetweenMorphismObject}    
\gamma_{x,y}:\A(x,y)\to \C(Fx,Fy).
\end{equation}

The next definition is very natural.

\begin{definition}[{\cite[Definition 5.3.1]{GepHauEnriched}}]
\label{defFullyFaith}
 Let $\V$ be a presentably monoidal $\infty$-category. A $\V$-enriched $\infty$-functor $F : \A \to \C$ is \textit{fully faithful} if the maps
\[\gamma_{x,y}:\A(x,y) \to \C(Fx, Fy)\] 
are equivalences in $\V$ for all pair $x,y$ of objects in $\A$.
\end{definition}

In \Cref{notHomotopyCategoryOfA} we defined the homotopy category $\infty$-functor $h:\Cat^{\V}\to Cat$, using this $\infty$-functor the next definition is natural. 

\begin{definition}[{\cite[Definition 5.3.3 and Lemma 5.3.4]{GepHauEnriched}}]
    \label{defEssSurj}
 Let $\V$ be a presentably monoidal $\infty$-category. A $\V$-enriched $\infty$-functor $F : \A \to \C$ is \textit{essentially surjective} if the induced map between categories \[hF : h\A \to
h\C\] is essentially surjective.
\end{definition}

In \cite{GepHauEnriched}, the authors prove that for a $\V$-enriched $\infty$-functor, being an equivalence or being essentially surjective and fully faithful is the same. So, we give a second definition of equivalence.

\begin{definition}
\label{defEquivalenceVHigherFunctor}
A $\V$-enriched $\infty$-functor $F$ is an equivalence if it is essentially surjective and fully faithful.
\end{definition}

Moreover, in \cite{GepHauEnriched} or \cite{heine2023equivalence}, the authors define the $\infty$-category of $\V$-enriched $\infty$-functor between two $\V$-enriched $\infty$-categories $\A$ and $\B$, we denote it by $Fun^{\V,GH}(\A,\B)$.

From now on, we analyze the point of view of Hinich. We start with his definitions of $\V$-enriched $\infty$-functor and of $\infty$-category of $\V$-enriched $\infty$-functors, and then we explain why his definitions are equivalent to \Cref{defHigherEnrichedFunctor} in some situations, those situations that interest us.

Let $\V$ be a presentably monoidal $\infty$-category; in particular, $\V\in Alg(Pr^{L})$. We know that $Pr^{L}$ has a monoidal structure, \eqref{eqTensorPrl}, and it is canonical to see $Pr^{L}$ as a $Pr^{L}$-left tensored $\infty$-categories. 
Let $\B\in \mathrm{LMod}_{\V}(Pr^{L})$ be a left $\V$-module object of $Pr^{L}$.

Let $\A$ be a $\V$-enriched $\infty$-category. A $\V$-enriched $\infty$-category is an associative algebra in $Quiv_{\Aob}(\V)$ and, from \Cref{rmkHigherAlgebraicStructureInQuivBM}, we know that $Quiv_{\Aob}^{LM}(\B)$ is a $Quiv_{\Aob}(\V)$-left $\infty$-category, so the next definition makes sense.

\begin{definition}[{\cite[\S 6.1.3]{HinichYon}}]
\label{defHigherCategoryofEnrichedFunctor}
 Let $\V$ be a presentably monoidal $\infty$-category. Let $\A\in Cat_{\infty,\Aob}^{\V}$ be an $\V$-enriched $\infty$-category with spaces of $\Aob$ and let $\B$ be a left $\V$-module object in $Pr^{L}$. An $\V$-enriched $\infty$-functor $F : \A\to \B$ is a left $\A$-module object of $Fun(\Aob,B)$.
We call $\mathrm{LMod}_{\A}(Fun(\Aob,\B))$ the \textit{$\infty$-category of $\V$-enriched $\infty$-functors from $\A$ to $\B$} and we denote it by $Fun^{\V,H}(\A,\B)$.
\end{definition}

Now we study a subcase of the above definition to make clear why \Cref{defHigherCategoryofEnrichedFunctor} is a correct definition. This subcase is the case in which we are interested.

Let $\V$ be a presentably monoidal $\infty$-category and let $\B$ be a tensored $\V$-enriched $\infty$-category with a presentable underlying $\infty$-category $\B_o$. In particular, this category is a left $\V$-module object in $Pr^{L}$, so the $\infty$-category of $\V$-enriched $\infty$-functors is well-defined. Let $F\in \mathrm{LMod}_{\A}(Fun(\Aob,\B))$ be a $\V$-enriched $\infty$-functor; it is composed of two parts: an object $F_{ob}$ of $Fun(\Aob,\B)$, which is a morphism between the spaces of objects:
\begin{equation}
\label{eqObjetRitrovato}   
F_{ob}:\Aob\to \Bob,
\end{equation}
and an arrow in $Fun(\Aob,\B)$ (i.e. a natural trasformation).

\[m_*: \oint^{x\in \Aob} A(x,*)\otimes F(x)\Rightarrow F(*). \]

In particular, for each object $y\in \Aob$, there is a morphism in $\V$:

\[m: \oint^{x\in \Aob} A(x,y)\otimes F(x)\to F(y). \]

Every coend is a colimit, so we can precompose the above morphism with the $(x,x)$-component of the colimit cone $\lambda_{(x,x)}$, and we obtain: 

\[A(x,y)\;^{\V}\otimes_{\B} F(y)\xrightarrow{\lambda_{(x,x)}} \oint^{y\in \Aob} A(x,y)\otimes F(y)\xrightarrow{m} F(y). \]

Finally, since the enriched $\infty$-category $\B$ is tensored, we can take the corresponding adjoint morphism in $\V$:

\begin{equation}
\label{eqMorphismRitrovato}
   \gamma_{x,y}: \A(x,y)\to \B(F(x),F(y)). 
\end{equation}

Summing up, in (\ref{eqMorphismRitrovato}) and (\ref{eqObjetRitrovato}), we found (\ref{eqmapBetweenMorphismObject}) and (\ref{eqMapBetweenObject}). We conclude that if $\B$ is a tensored $\V$-enriched $\infty$-category with a presentable underlying $\infty$-category $\B_o$, \Cref{defHigherCategoryofEnrichedFunctor} and \Cref{defHigherEnrichedFunctor} coincide. After the above, the next result can be expected.

\begin{proposition}
Let $\V$ be a presentably monoidal $\infty$-category, and let $\A$ be a $\V$-enriched $\infty$-category and $B$ is a tensored $\V$-enriched $\infty$-category with a presentable underlying $\infty$-category $\B_o$. Then there is an equivalence of $\infty$-categories
\[ Fun^{\V,H}(\A,\B)\simeq Fun^{\V,GH}(\A,\B) ,\]
\end{proposition}
$\\$
For the proof of the above proposition, see \cite{heine2023equivalenceTwo}.
$\\$
As a consequence of the last result, we set a common notation.
\begin{notation}
We denote by $Fun^{\V}(\A,\B)$ both the $\infty$-category $Fun^{\V,GH}(\A,\B)$ than $Fun^{\V,H}(\A,\B)$, where $\A$ and $\B$ are two suitable $\V$-enriched $\infty$-categories. 
\end{notation}
$\\$
In this series of articles, except in \cite[ Proposition 6.7]{DoniklinearMorita} with $\infty$-category of enriched $\infty$-functors we intend the Hinich's $\infty$-category.

%\begin{proposition}
%Let $\V$ be a presentably monoidal $\infty$-category, and let $\A$ be a $\V$-enriched $\infty$-category and $B$ is a tensored $\V$-enriched $\infty$-category with a presentable underlying $\infty$-category $\B_o$. Then there is an equivalence of $\infty$-categories
%\[ Fun^{\V,GH}(A,B)\simeq Fun^{\V,GH}(A,B) ,\]
%\end{proposition}

%For a proof of the above proposition, see \cite{heine2023equivalence}.
%As a consequence of the last result, we set a common notation.
%\begin{notation}
%Let $\A$ and $\B$ be two suitable $\V$-enriched $\infty$-categories. 
%We denote by $Fun^{\V}(\A,\B)$ both the $\infty$-categories $Fun^{\V,GH}(\A,\B)$ and $Fun^{\V,H}(\A,\B)$.
%\end{notation}
%In the last part of this section, we say something about the $\infty$-category of $\V$-enriched presheaves.

%First of all, we need to explain what the opposite of a $\V$-enriched $\infty$-category is. Since, for us, it is enough, we add the symmetric property to $\V$, and so we fix a symmetric monoidal $\infty$-category $\V$.

%%%%%%%%%%%%%%%%%%%%%%%%%%%%%%%%%%%%%%%%%%%%%%%%%%%%%%%%%%%%%%%%%%%%%%%%%%%%%%%%%%%%%%%%%%%%%%%%%%%%%%%%%%%%%%%%%%%%%%%%%%%%%%%%%%%%%%%%%%%%%%%%%%%%%%%%%%%%%%%%%%%%%%%%%%%%%%%%%%%%%%%%%%%%%%%%%%%%%%%%%%%%%%%%%%%%%%%%%%%%%%%%%%%%%%%%%%%%%%%%%%%%%%%%%%%%%%%%%%%%%%%%%%%%%%%%%%%%%%%%%%%%%%%%%%%%%%%%%%%%%%%%%%%%%%%%%%%%%%%%%%%%%%%%%%%%%%%%%%%%%%%%%%%%%%%%%%%%%%%%%%%%%%%%

%%%%%%%%%%%%%%%%%%%%%%%%%%%%%%%%%%%%%%%%%%%%%%%%%%%%%%%%%%%%%%%%%%%%%%%%%%%%%%%%%%%%%%%%%%%%%%%%%%%%%%%%%%%%%%%%%%%%%%%%%%%%%%%%%%%%%%%%%%%%%%%%%%%%%%%%%%%%%%%%%%%%%%%%%%%%%%

In the last part of this subsection, we reformulate \Cref{LemmaHighLeftRightAdjointforgMod} and \Cref{LemmaRightaadjointforgetful}, in particular, we start to consider module objects and enriched $\infty$-functors as something strongly related. 

\begin{remark}
\label{rmkAdjunctionHigherCatextentionrestrictioncorestriction}
Let $R\in Alg$ be an algebra object of $\Sp$. By the definition of the $\infty$-category of $\Sp$-enriched $\infty$-functors, \Cref{defHigherEnrichedFunctor}, the following equivalences hold \[\mathrm{LMod}_{R}\simeq \mathrm{LMod}_{\un{R}}(Fun(*,Sp))\simeq Fun^{\Sp}(\un{R},\un{Sp}),\] and so, \Cref{LemmaHighLeftRightAdjointforgMod} and \Cref{LemmaRightaadjointforgetful} induces the following adjunction:

% https://q.uiver.app/?q=WzAsMyxbMCwwLCJcXHVue0N9XlxcdW57MX0iXSxbMiwwLCJcXHVue0N9XlxcdW57Un0iXSxbNCwwLCJcXHVue0N9XlxcdW57MX0iXSxbMCwxLCItXFxvdGltZXMgUiIsMSx7ImN1cnZlIjotMn1dLFsyLDEsIlJhbl97cn0tIiwxLHsiY3VydmUiOi0yfV0sWzEsMiwiVSIsMSx7ImN1cnZlIjotMn1dLFsxLDAsIlUiLDEseyJjdXJ2ZSI6LTJ9XSxbMyw2LCIiLDAseyJsZXZlbCI6MSwic3R5bGUiOnsibmFtZSI6ImFkanVuY3Rpb24ifX1dLFs1LDQsIiIsMCx7ImxldmVsIjoxLCJzdHlsZSI6eyJuYW1lIjoiYWRqdW5jdGlvbiJ9fV1d
\[\begin{tikzcd}
	{\Sp\simeq Fun^{\Sp}(\un{\mathbb{S}},\Sp)} && {Fun^{\Sp}(\un{R},\Sp)} && { Fun^{\Sp}(\un{\mathbb{S}},\Sp)\simeq \Sp.}
	\arrow[""{name=0, anchor=center, inner sep=0}, "{-\otimes R}"{description}, bend left, from=1-1, to=1-3]
	\arrow[""{name=1, anchor=center, inner sep=0}{description}, bend left, from=1-5, to=1-3]
	\arrow[""{name=2, anchor=center, inner sep=0}, "U"{description}, bend left, from=1-3, to=1-5]
	\arrow[""{name=3, anchor=center, inner sep=0}, "U"{description}, bend left, from=1-3, to=1-1]
	\arrow["\dashv"{anchor=center, rotate=-90}, draw=none, from=0, to=3]
	\arrow["\dashv"{anchor=center, rotate=-90}, draw=none, from=2, to=1]
\end{tikzcd}\]

\end{remark}

It is possible to prove that $U$ is the restriction $\infty$-functor via the unit morphism $\I_{R}:\mathbb{S}\to R$. The above construction can be generalized to every morphism of algebras $f: E\to R$.

\begin{corollary}
\label{corExtentionRestriction}
Let $f: E\to R$ be an arrow of algebras objects of $\Sp$. The restriction (or forgetful) $\infty$-functor $U: \mathrm{Mod}_{R}\to \mathrm{Mod}_{E}$ has both a left adjoint, which is called \textit{the extension of scalar}, and the right adjoint, which is called \textit{the coextension of scalar}. 
\end{corollary}

% https://q.uiver.app/?q=WzAsMyxbMCwwLCJcXHVue0N9XlxcdW57MX0iXSxbMiwwLCJcXHVue0N9XlxcdW57Un0iXSxbNCwwLCJcXHVue0N9XlxcdW57MX0iXSxbMCwxLCItXFxvdGltZXMgUiIsMSx7ImN1cnZlIjotMn1dLFsyLDEsIlJhbl97cn0tIiwxLHsiY3VydmUiOi0yfV0sWzEsMiwiVSIsMSx7ImN1cnZlIjotMn1dLFsxLDAsIlUiLDEseyJjdXJ2ZSI6LTJ9XSxbMyw2LCIiLDAseyJsZXZlbCI6MSwic3R5bGUiOnsibmFtZSI6ImFkanVuY3Rpb24ifX1dLFs1LDQsIiIsMCx7ImxldmVsIjoxLCJzdHlsZSI6eyJuYW1lIjoiYWRqdW5jdGlvbiJ9fV1d
\begin{equation}
    \label{eqSixFnctormod}
\begin{tikzcd}
	{Fun^{\Sp}(\un{E},\Sp)} && {Fun^{\Sp}(\un{R},\Sp)} && {Fun^{\Sp}(\un{E},\Sp).}
	\arrow[""{name=0, anchor=center, inner sep=0}, "{f_{!}}"{description}, bend left, from=1-1, to=1-3]
	\arrow[""{name=1, anchor=center, inner sep=0},"f_*"{description}, bend left, from=1-5, to=1-3]
	\arrow[""{name=2, anchor=center, inner sep=0}, "f^*"{description}, bend left, from=1-3, to=1-5]
	\arrow[""{name=3, anchor=center, inner sep=0}, "f^*"{description}, bend left, from=1-3, to=1-1]
	\arrow["\dashv"{anchor=center, rotate=-90}, draw=none, from=0, to=3]
	\arrow["\dashv"{anchor=center, rotate=-90}, draw=none, from=2, to=1]
\end{tikzcd}
\end{equation}

%% file: Chapters/HigherCategoricalArticle.tex
\section{Results}
\label{secLModEnrichedareModuleObjects}

%%%%%%%%%%%%%%%%%%%%%%%%%%%%%%%%%%%%%%%%%%%%%%%%%%%%%%%%%%%%%%%%%%%%%%%%%%%%%%%%%%%%%%%%%%%%%%%%%%%%%%%%%%%%%%%%%%%%%%%%%%%%%%%%%%%%%%%%%%%%%%%%%%%%%%%%%%%%%%%%%%%%%%%%%%%%%%%%%%%%%%%%%%%%%%%%%%%%%%%%%%%%%%%%%%%%%%%%%%%%%%%%%%%%%%%%%%%%%%%%%%%%%%%%%%%%%%%%%%%%%%%%%%%%%%%%%%%%%%%%%%%%%%%%%%%%%%%%%%%%%%%%%%%%%%%%%%%%%%%%%%%%%%%%%%%%%%%%%%%%%%%%%%%%%%%%%%%%%%%%%%%%%%%%%%%%%%%%%%%%%%%%%%%%%%%%%%%

%%%%%%%%%%%%FINO A qUA

\begin{theorem}
\label{thCatEnInLmodAreLmod}
Let $(\V,-\otimes_{\V},\mathbb{I}_{\V}),$ be a presentable $\mathrm{E}_{2}$-monoidal $\infty$-category such that $\otimes_{\V}$ preserves small colimits componentwise. Let $R$ be an $\mathbb{E}_2$-algebra of $\V$,
then there exists an equivalence of $\infty$-categories 
\begin{equation}
\label{eqCatEnInLmodAreLmod}
\Cat^{\mathrm{LMod}_{R}(\V)}\simeq \mathrm{LMod}_{R}(\Cat^{\V}).
\end{equation}
\end{theorem}

\begin{proof}
We aim to apply the Barr-Beck-Lurie theorem \cite[Theorem 4.7.3.5]{HA} to the adjunction
\[-\otimes R_{!}:\Cat^{\V}\rightleftarrows \Cat^{\mathrm{LMod}_{R}(\V)}:U_{!},\]
which we found in \Cref{exAdjunctionEnrichedVLMODRV}.

We claim that the adjunction $\otimes R_{!}\dashv U_{!}$ satisfies the hypotheses of the Barr-Beck-Lurie theorem, i.e. the functor $U_{!}$ is conservative and it creates colimits of $U_{!}$-split simplicial objects.

Assuming the claim to be true, then we obtain an equivalence:
\begin{equation}
\label{eqAdriano}
\Cat^{\mathrm{LMod}_{R}(\V)}\simeq \mathrm{LMod}_{T}(\Cat^{\V})
\end{equation}
where $T$ is the monad $T=U_{!}\circ (-)\otimes{R}_{!}: \Cat^{\V}\to \Cat^{\V}$.

In particular, since the $\infty$-functor $\Cat^{(-)}$ is functorial, the following chain of equivalences holds:
\[T=U_{!}\circ (-)\otimes{R}_{!}\simeq (U\circ (-)\otimes R)_{!} \simeq (-)\otimes_{\V} R_{!}.\]

The adjunction forgetful-free 
\[ F:\Cat^{\V}\leftrightarrows \mathrm{LMod}_{R}(\Cat^{\V}):G, \]
is monadic, and its associated monad is $T_{1}=R ^{\V}\otimes^{\V}-: \Cat^{\V}\to \Cat^{\V}$, see \cite[Example 4.7.3.9]{HA}.
From the definitions of the various tensors, the following equivalences of endofunctors hold:
\[  
T_{1}=R ^{\V}\otimes^{\V}(-)\simeq \un{R}\otimes^{\V}(-)\simeq (-)\otimes^{\V}\un{R}\simeq  (U\circ (-)\otimes{R})_{!}= T:\Cat^{\V}\to \Cat^{V}.
\]
Thus, we obtain the following chain of equivalences:
\begin{equation}
\label{eqObaObaMartins}    
 \mathrm{LMod}_{T}(\Cat^{\V})\simeq \mathrm{LMod}_{T'}(\Cat^{\V})\simeq \mathrm{LMod}_{R}(\Cat^{\V}).
\end{equation}
Combining \eqref{eqAdriano} and \eqref{eqObaObaMartins}, we obtain the thesis.

Next, we prove the claim.
We start by proving that $U_{!}$ is conservative.
Let $H:\A\to \B$ be a morphism in $\Cat^{\LModk}$ such that its image through $U_{!}$, \[U_{!}H:U_{!}\A\to U_{!}\B,\] is an equivalence in $\Cat^{\V}$. This means that $U_{!}H$ has two properties:
\begin{itemize}
    \item It is fully faithful, i.e., for each pair of elements $m,n$ of $U_{!}\A$, the morphism \begin{equation}
    \label{eqU!H}    
    U_{!}H_{m,n}:U_{!}\A(m,n)=U\A(m,n)\to U_{!}\B(U_{!}H(m),U_{!}H(n))=U\B(U_{!}H(m),U_{!}H(n))
    \end{equation}
    (where the equals hold by definition of $\Cat^{(-)}$) is an equivalence in $\V$;
    \item It is essentially surjective, i.e., the induced functor between the homotopy categories \[ hU_{!}H: hU_{!}\A\to hU_{!}\B\] is an (ordinary) equivalence.
 \end{itemize}
Now we prove that $H$ is essentially surjective and also fully faithful.

Since $hH$ and $hU_{!}H$ are the same (ordinary) functors, $U_{!}H$ is essentially surjective. 
In particular, this means that for each $m\in\A$ the two objects $U_{!}H(m)$ and $H(m)$ are equivalent in $\B$ (i.e. they are equivalent in $\B_o$) and, without loosing generalizaion, we can considered $U_{!}H(-)_{ob}$ and $H(-)_{ob}$ the same $\infty$-functor between $\infty$-groupoids. 
With an abuse of notation, for each $m\in\A$, we denote by $H(m)$ the object $U_{!}H(m)$.   

On the hom-objects the $\infty$-functor $U_{!}$ sends for each pair of elements $m,n\in \A$ sends a $\mathrm{LMod}_{R}(\V)$-enriched $\infty$-functor \[H_{m,n}:\A(m,n)\to \B(H(m),H(n))\] into a $\infty$-functor in $\V$ \[UH_{m,n}:U\A(m,n)\to U\B(H(m),H(n))\]
where is the monoidal $\infty$-forgetful $U:\mathrm{LMod}_{R}(\V)\to \V$. In particular, $U$ is conservative so if $UH_{m,n}$ is an equivalence then $H_{m,n}$ is. This proves that $U_{!}$ is conservative.

Now, we will show that $U_{!}$ creates $U_{!}$-split simplicial objects.

By \cite[Corollary 4.2.3.6]{HA} and \eqref{eqFunCatEr}, $\Cat^{\mathrm{LMod_{R}(\V)}}$ admits all small colimits and then every $U_{!}$-split simplicial object has colimits.

So it only remains to prove that $U_{!}$ preserves colimits of $U_{!}$-split simplicial objects.  

Since $\Delta^{op}$ is a sifted $\infty$-category, \cite[Lemma 5.5.8.4]{HTT}, every $U_{!}$-split simplicial object is a sifted diagram, so if we prove that $U_{!}$ preserves sifted colimits then it preserves colimits of $U_{!}$-split simplicial objects.

Furthermore, using \cite[Lemma 5.7.7]{GepHauEnriched}, it suffices to prove that the $\infty$-functor between the $\infty$-categories of enriched precategories 
\[ U_{!}: PreC(\mathrm{LMod}_{R}(\V))\to PreC(\V)  \]
preserves the sifted colimits. 

Now, we will use the strategy that Gepner and Haugseng used in the first part of the proof of \cite[Proposition A.5.10]{GepHauEnriched}.  

The strategy is: $U_{!}$ preserves sifted colimits if $U_{!}^{triv}$ preserves them. In our case, $U_{!}^{triv}:Quiv(\mathrm{LMod}_{R}(\V))\to Quiv(\V)$ is a morphism in the $\infty$-category of cocartesian fibration over $\mathcal{S}$, it is denoted by $CoCart(\mathcal{S})$, which we will explicitly describe as a natural transformation between $\infty$-functors with target $\Top$ and target $\Cat$ using the Grothendieck construction.  
Thus, it is suffices to show that the morphism $U_{!}^{triv}:Quiv(\mathrm{LMod_{R}(\V)})\to Quiv(\V)$ in $CoCart(\mathcal{S})$ preserves sifted colimits.

Using the Grothendieck construction, we obtain that the above $Quiv(\mathrm{LMod}_R(\V))$ is the $\infty$-functor \[Quiv(\mathrm{LMod}_R(\V)):\mathcal{S}\to \Cat:X\mapsto Quiv_{X}(\mathrm{LMod}_R(\V))= Fun(X^{op}\times X,\mathrm{LMod}_R(\V)),\]
$Quiv(\V)$ is the $\infty$-functor \[Quiv(\V):\mathcal{S}\to \Cat:X\mapsto Quiv_{X}(\V)=Fun(X^{op}\times X,\V)\] and $U_{!}^{triv}:Quiv(\mathrm{LMod}_R(\V))\to Quiv(\V)$ is a natural transformation such that, for each $X\in \mathcal{S}$, the $X$-component is the post-composition with $U$ \[U^{triv}_{!,X}:Fun(X^{op}\times X,\mathrm{LMod}_R(\V))\xrightarrow{U\circ (-)}Fun(X^{op}\times X,\V).\]
Since colimits in an $\infty$-functors $\infty$-category are computed in the target, a natural transformation preserves colimits if and only if any component preserves them.
In our case, for any $X\in\mathcal{S}$, $U^{triv}_{!,X}$ preserves colimits  because it is a postcomposition with an $\infty$-functor that preserves colimits, see \Cref{LemmaRightaadjointforgetful}.  
In particular, this means that $U_{!}^{triv}$ preserves sifted colimits and hence $U_{!}: PreCat(\mathrm{LMod}_R(\V))\to PreCat(\V)$ preserves them too.

In conclusion, we proved the claim and also the thesis.
\end{proof}

The next corollary was suggested to me by Haugseng, and for that I thank him.
We start with a remark.

\begin{remark}
\label{rmkMonoidalStructureOfLModR}%\textcolor{red}{to clarify}
Let $R$ and $\V$ as in \Cref{thCatEnInLmodAreLmod}.
Using \eqref{eqCatEnInLmodAreLmod}, we can suit $\mathrm{LMod}_{R}(\V)$ with a canonical structure of monoidal $\infty$-category such that the equivalence \eqref{eqCatEnInLmodAreLmod} is monoidal.
Is not clear to me if when $R$ is an $\E$-ring, e.g. $R=\Hk$, this induced monoidal structure on $\mathrm{LMod}_{R}(\V)$ is equivalent to relative tensor product monoidal structure, see \Cref{notRelativeTensorProduct}.
    
\end{remark}

\begin{corollary}
\label{corIteratedResult}
Let $\V$ be a presentable $\mathbb{E}_{n+1}$-monoidal $\infty$-category with $n\in\N_{\geq 2}\cup \{\infty\}$ and let $R$ be an $\mathbb{E}_{n+1}$-algebra of $\V$.
Then there exists an equivalence of $\infty$-category
\[   \mathcal{C}at_{(\infty,n)}^{\mathrm{LMod}_{R}(\V)}\simeq \mathrm{LMod}_{\mathcal{B}^{n}R}(\mathcal{C}at_{(\infty,n)}^{\V}).\]
    
\end{corollary}
\begin{proof}

We proceed by induction. 
The base case is $n=2$.
The following equivalences hold and are well defined:

\begin{equation}
\begin{split}
\mathcal{C}at_{(\infty,2)}^{\mathrm{LMod}_{R}(\V)} & \simeq \mathcal{C}at_{(\infty,1)}^{Cat_{(\infty,1)}^{\mathrm{LMod_{R}(\V)}}} \simeq \mathcal{C}at_{(\infty,1)}^{\mathrm{LMod}_{\mathcal{B}R}(Cat_{(\infty,1)}^{\V})}
\\ & \simeq  \mathrm{Mod}_{\mathcal{B}^{2}R}(\mathcal{C}at_{(\infty,2)}^{\V}).
\end{split}
\end{equation}
All the $\infty$-categories are well define because \cite[Corollary 5.7.12, Remark 5.7.13]{GepHauEnriched}.

The first equation holds by definition, while the second and the last ones follow from \Cref{thCatEnInLmodAreLmod}.

We assume the theorem holds for $n-1$ and proceed to prove the induction step. In the following chain of equivalences, all the $\infty$-categories are well-defined due to \cite[Corollary 5.7.12, Remark 5.7.13]{GepHauEnriched}. Moreover, the second equivalence holds by \Cref{thCatEnInLmodAreLmod}, the first and the last ones by definition, and the third by induction hypothesis.

\begin{equation}
\begin{split}
\mathcal{C}at_{(\infty,n)}^{\mathrm{LMod}_{R}(\V)} & \simeq \mathcal{C}at_{(\infty,n-1)}^{Cat_{(\infty,1)}^{\mathrm{LMod_{R}(\V)}}} \simeq \mathcal{C}at_{(\infty,n-1)}^{\mathrm{LMod}_{\mathcal{B}R}(Cat_{(\infty,1)}^{\V})}
\\ & \simeq  \mathrm{Mod}_{\mathcal{B}^{n-1}\mathcal{B}R}(\mathcal{C}at_{(\infty,n)}^{\V})\simeq  \mathrm{Mod}_{\mathcal{B}^{n}R}(\mathcal{C}at_{(\infty,n)}^{\V}).
\end{split}
\end{equation}

\end{proof}

Above, we have established the presentable $\infty$-categorical version of \cite[Theorem 3.1]{DoniCategorical}: we have demonstrated that even in the $\infty$-categorical context, the theory of $\mathrm{LMod}_{R}(\V)$-enriched $\infty$-categories can be comprehensively studied within the framework of $\V$-enriched $\infty$-categories.

We are interested in further describing $\mathrm{LMod}_{R}(\V)$-enriched $\infty$-categories as $\Cat^{\V}$-enriched $\infty$-functors with source $\2un{\V}$. Since $\Cat^{\V}$ is presentable and its tensor preserves colimits componentwise (\cite[Corollary 9.4.]{haugseng2023tensor}), it is a left $\Cat^{\V}$-module in $Pr^{L}$.
So, $Fun^{\Cat^{\V}}(\underline{\underline{R}},\Cat^{\V})$ is well-defined in Hinich's setting, see \Cref{defHigherCategoryofEnrichedFunctor}.

Now, we can state the second result of this article.
\begin{theorem}
\label{thAction}
Let $\V$ be a presentable $\mathbb{E}_3$-monoidal $\infty$-category and let $R$ be an $\mathbb{E}_2$-algebra of $\V$.

There exists a chain of equivalences of $\infty$-categories:
\begin{equation}
\label{eqHkaction}
\mathrm{LMod}_{R}(\Cat^{\V}) \simeq \mathrm{LMod}_{\underline{R}}(\Cat^{\V})\simeq  Fun^{\Cat^{\V}}(\underline{\underline{R}}, \Cat^{\V}).
\end{equation}
\end{theorem}
\begin{proof}
The following chain of equivalences holds:
\begin{equation}
    \begin{split}
    Fun^{\Cat^{\V}}(\underline{\underline{R}}, \Cat^{\V}) \simeq & \mathrm{LMod}_{\underline{\underline{\Hk}}}(Fun(*,\Cat^{\V})) \\ \simeq & \mathrm{LMod}_{\underline{R}}(\Cat^{\V}) \\ \simeq & \mathrm{LMod}_{R}(\Cat^{\V}) . 
    \end{split}
\end{equation}

Since the spaces of objects of $\underline{\underline{R}}$ is the final spaces, $\underline{\underline{R}}_{ob}\simeq *$, the first equivalence is valid by definition of the $\infty$-category of $\Cat^{\V}$-enriched $\infty$-functors with target a left $\Cat^{\V}$-module in $Pr^{L}$, \Cref{defHigherCategoryofEnrichedFunctor}. The second follows from the well-known equivalence $ev:Fun(*,\Cat^{\V})\simeq \Cat^{\V}$. The latter follows from the definition of the left $Alg$-module object of $Pr^{L}$ structure of $\Cat^{\V}$, see \eqref{eqAlgVtensoredCAtCatV}.
\end{proof}

Let $\V$ an presentable $\mathbb{E}_{2}$ monoidal $\infty$-category, the following underlying $\infty$-functor is well-defined:
\begin{equation}
    \underline{(-)}:Alg(Alg(\V))\to \Cat^{Alg(\V)}.
\end{equation}
We know that $\Cat^{\V}$ is a left $Alg(\V)$-module object in $Pr^{L}$, see \eqref{eqCalgTensoredCatsp}, so it is well-defined the $\infty$-category of $Alg(\V)$-enriched $\infty$-functors $Fun^{Alg(\V)}(\underline{R},\Cat^{\V})$ in Hinich's setting, see \Cref{defHigherCategoryofEnrichedFunctor}. 

Now, we can state the last result of this paper.

\begin{proposition}
\label{eqCatLModkdescritionCalgtensored}
There exists a chain of equivalences of $\infty$-categories:
\begin{equation}
\Cat^{\mathrm{LMod}_{R}(\V)}\simeq \mathrm{LMod}_{R}(\Cat^{\V})\simeq Fun^{Alg(\V)}(\underline{R},\Cat^{\V}).
\end{equation}
\end{proposition}
\begin{proof}
The first equivalence follows from \Cref{thAction} and \eqref{eqCatEnInLmodAreLmod}.
Using the definition of $\infty$-category of $Alg(\V)$-enriched $\infty$-functors $\infty$-category, see \Cref{defHigherCategoryofEnrichedFunctor}, and the fact that $\underline{\V}$ is a $Alg(\V)$-enriched $\infty$-category with a trivial space of objects, i.e., $\underline{R}_{ob}\simeq * $, we obtain the following chain of equivalences:
\begin{equation}
\label{eqFareChiarezza}
\mathrm{LMod}_{R}(\Cat^{\V})\simeq \mathrm{LMod}_{\underline{R}}(Fun(*,\Cat^{\V}))\simeq Fun^{Alg(\V)}(\underline{R},\Cat^{\V}).
\end{equation}
The latter equivalence in \eqref{eqFareChiarezza} follows from the well-known equivalence $ev_*:Fun(*,\V)\simeq \V$.  
\end{proof}

%%%%%%%%%%%%%%%%%%%%%%%%%%%%%%%%%%%%%%%%%%%%%%%%%%%%%%%%%%%%%%%%%%%%%%%%%%%%%%%%%%%%%%%%%%%%%%%%%%%%%%%%%%%%%%%%%%%%%%%%%%%%%%%%%%%%%%%%%%%%%%%%%%%%%%%%%%%%

\subsection{Specific case: $\V=\Sp$ and $R=\Hk$}
\label{subsecSpecificCase}

In this subsection, we specialize the results of this paper to an interesting case for derived algebraic geometry: $\V=\Sp$ and $R=\Hk$.
 
At the end, an $\LModk$-enriched $\infty$-categories will mean five things to us:
\begin{itemize}
    \item[(1)] a $\LModk$-enriched $\infty$-category;
    
    \item[(2)] or, a $\Cat^{\Sp}$-enriched $\infty$-functor from the $\Cat^{\Sp}$-enriched $\infty$-category $\2un{\Hk}$ to $\Cat^{\Sp}$;
    
    \item[(3)] or, a left $\un{\Hk}$-module object of $\Cat^{\Sp}$;

    \item[(4)] or, a $\Cat^{\Sp}$-enriched $\infty$-functor from the $Alg$-enriched $\infty$-category $\un{\Hk}$ to the $Alg$-left tensored $\infty$-category $\Cat^{\Sp}$;

    \item[(5)] or, a left $\Hk$-module object of $\Cat^{\Sp}$. 

\end{itemize}

This result is crucial because it will allow us to define a $k$-linearization in different settings. Indeed, the quintuple vision above is the starting point for a general notion of $k$-linearization, and it will allow us to define a comparable Morita theory in the pretriangulated $dg$-categories and stable $\infty$-categories settings. However, this is a topic for another article, \cite{DoniklinearMorita}, that makes up the series. 

\begin{notation}
We refer to objects of $\Cat^{\Sp}$ as \textit{spectral $\infty$-categories}. Sometimes, we call this as the \textit{non-geometric case}.
Instead, we call $\Cat^{\LModk}$ the \textit{$k$-linear case}. 
\end{notation}

\begin{notation}
\label{notunderlyingspectralcategory}
    Let $\A\in\Cat^{\LModk}$ be a $\LModk$-enriched $\infty$-category, we call \textit{the underlying spectral $\infty$-category of $\A$} the spectral $\infty$-category $U_{!}(\A)$ where $U_{!}$ is as in \Cref{exAdjunctionEnrichedSpLMod}. We denote it by $\A_{sp}$. 
\end{notation}

For this part, we fix a commutative unitary ring $k$.
 
Let us specialize some results of \Cref{sussecMonoidalEnrichedCat} to the cases $\V=\Sp$ or $\V=\LModk$; both are presentably symmetric monoidal $\infty$-categories.
In \Cref{sussecMonoidalEnrichedCat}, we defined the underline $\infty$-functor:  
\[\underline{(-)}:Alg\simeq Cat^{\Sp}_{\infty,*}\to\Cat^{\Sp}: \A\to \underline{\A}.\]
It maps an algebra object of $\Sp$ to the $\Sp$-enriched $\infty$-category whose space of objects is the terminal space $*$, with only non-trivial hom-$\Sp$-object $\A\simeq \underline{\A}(*,*)$.
Moreover, we defined a presentably symmetric monoidal $\infty$-category $(\Cat^{\Sp},\otimes^{\Sp}, \underline{\mathbb{S}})$, and using $\underline{(-)}$, we equipped $\Cat^{\Sp}$ with a structure of left $Alg$-module objects of $Pr^{L}$, see (\ref{eqAlgVtensoredCAtCatV}):

\begin{equation}
    \label{eqCalgTensoredCatsp}
       -^{\Sp}\otimes^{\Sp} -:Alg \times \Cat^{\Sp}\xrightarrow{\underline{(-)}\times id}  \Cat^{\Sp}\times\Cat^{\Sp}\xrightarrow{-\otimes^{\Sp} -} \Cat^{\Sp}.
\end{equation}

\begin{remark}
\label{rmkHKAlgebraDiAlgebraDgCat}
Since $\Hk$ is an object of $CAlg$, it is in particular an algebra in the symmetric monoidal $\infty$-category $Alg$ because there is the canonical forgetful $\infty$-functor: 
\begin{equation}
\label{eqHKAlgebraDiAlgebra2}
U:CAlg= Alg_{\E}(\Sp)\to Alg_{\mathbb{E}_{2}}(\Sp)\simeq Alg(Alg).
\end{equation}
As a consequence, the $\infty$-category $\mathrm{LMod}_{\Hk}(\Cat^{\Sp})$ is well defined, and so the statement in \Cref{thCatEnInLmodAreLmod} and \Cref{thAction} will be well placed.
This observation clarifies why it is important that $\Hk$ is a commutative ring (or an $\E$-ring).
We suspect that it is sufficient that it is an $\mathbb{E}_3$-ring, but in this paper, we are not interested in this generalization.
\end{remark}

Furthermore, we defined the double-underline $\infty$-functor
\[\2un{(-)}:CAlg\to\Cat^{\Cat^{\Sp}}: \A\to \2un{\A}.\]
It sends a commutative algebra object in $\Sp$ to the $\Cat^{\Sp}$-enriched $\infty$-category whose space of objects is the terminal space $*$, with only non-trivial hom-$\Cat^{\Sp}$-object $\underline{A}\simeq \2un{\A}(*,*)$.

We also defined a symmetric monoidal $\infty$-category $(\Cat^{\LModk},\otimes^{\LModk}, \underline{\Hk})$, using the same strategy.

\begin{corollary}
\label{crCatEnInLmodAreLmodHKSp}
There exists an equivalence of $\infty$-categories 
\begin{equation}
\label{eqCatEnInLmodAreLmod}
\Cat^{\LModk}\simeq \mathrm{LMod}_{\Hk}(\Cat^{\Sp}).
\end{equation}
\end{corollary}

\begin{proof}
 This result is \Cref{thCatEnInLmodAreLmod} with $\V=\Sp $ and $R=\Hk$.   
\end{proof}

Above, we have established the $\infty$-categorical version of \cite[Theorem 3.1]{DoniCategorical}: we have demonstrated that even in the $\infty$-categorical context, the theory of $\LModk$-enriched $\infty$-categories can be comprehensively studied within the framework of $\Sp$-enriched $\infty$-categories.

We are interested in further describing $\LModk$-enriched $\infty$-categories as $\Cat^{\Sp}$-enriched $\infty$-functors with source $\2un{\Hk}$. Since $\Cat^{\Sp}$ is presentable and its tensor preserves colimits componentwise (\cite[Corollary 9.4.]{haugseng2023tensor}), it is a left $\Cat^{\Sp}$-module in $Pr^{L}$.
So, $Fun^{\Cat^{\Sp}}(\underline{\un{\Hk}},\Cat^{\Sp})$ is well-defined in Hinich's setting, see \Cref{defHigherCategoryofEnrichedFunctor}.

\begin{corollary}
\label{crAction}
There exists a chain of equivalences of $\infty$-categories:
\begin{equation}
\label{eqHkactionHKSp}
\mathrm{LMod}_{\Hk}(\Cat^{\Sp}) \simeq \mathrm{LMod}_{\un{\Hk}}(\Cat^{\Sp})\simeq  Fun^{\Cat^{\Sp}}(\underline{\un{\Hk}}, \Cat^{Sp}).
\end{equation}
\end{corollary}

\begin{proof}
    This result is \Cref{thAction} with $\V=\Sp$ and $R=\Hk$. 
\end{proof}

\begin{remark}
\Cref{thAction} can be seen as the categorical version of the equivalent condition of $R$-action in \cite[\S D.1]{SAG} without $(d)$. We hope that $(d)$ will become a topic for future work because it is linked to the theory of additive $dg$-categories. 
\end{remark}

Since $\Hk$ is an object of $CAlg$, in particular, it is an algebra in a symmetric monoidal $\infty$-category $Alg$ because there is the canonical forgetful $\infty$-functor:
\begin{equation}
\label{eqHKAlgebraDiAlgebra1}
U:CAlg\simeq Alg_{\E}(\Sp)\to Alg_{\mathbb{E}_{2}}(\Sp)\simeq Alg(Alg).
\end{equation}
Furthermore, we have the underlying $\infty$-functor of this situation:
\begin{equation}
    \underline{(-)}:Alg(Alg)\to \Cat^{Alg}.
\end{equation}
We know that $\Cat^{\Sp}$ is a left $Alg$-module object in $Pr^{L}$, see (\ref{eqCalgTensoredCatsp}), so it is well-defined the $\infty$-category of $Alg(Sp)$-enriched $\infty$-functors $Fun^{Alg}(\un{\Hk},\Cat^{\Sp})$ in Hinich's setting, see \Cref{defHigherCategoryofEnrichedFunctor}. 

Now, we can state the last result of this paper.

\begin{corollary}
\label{crCatLModkdescritionCalgtensored}
There exists a chain of equivalences of $\infty$-categories:
\begin{equation}
\Cat^{\LModk}\simeq \mathrm{LMod}_{\Hk}(\Cat^{\Sp})\simeq Fun^{Alg}(\un{\Hk},\Cat^{\Sp}).
\end{equation}
\end{corollary}
\begin{proof}
    This result is \Cref{eqCatLModkdescritionCalgtensored} with $\V=\Sp$ and $R=\Hk$.
\end{proof}

%\begin{corollary}
%\label{eqCatLModkdescritionCalgtensored}
%There exists a chain of equivalences of $\infty$-categories:
%\begin{equation}
%\Cat^{\LModk}\simeq \mathrm{LMod}_{\Hk}(\Cat^{\Sp})\simeq Fun^{Alg}(\un{\Hk},\Cat^{\Sp}).
%\end{equation}
%\end{corollary}

%This last description is important because it allows us to define a $k$-linearization without considering the entire $\infty$-category of $\Cat^{\Sp}$-enriched $\infty$-functors but only its full $\infty$-subcategory of $Alg$-enriched $\infty$-functors. \Cref{eqCatLModkdescritionCalgtensored} will be a crucial fact when, in \cite{DoniklinearMorita}, we will compare the theories of Morita for pretriangulated $dg$-categories and stable $\infty$-categories.  

%%%%%%%%%%%%%%%%%%%%%%%%%%%%%%%%%%%%%%%%%%%%%%%%%%%%%%%%%%%%%%%%%%%%%%%%%%%%%%%%%%%%%%%%%%%%%%%%%%%%%%%%%%%%%%%%%%%%%%%%%%%%%%%%%%%%%%%%%%%%%%%%%%%%%%%%%%%%%%%%%%%%%%%%%%%%%%%%%%%%%%%%%%%%%

%%%%%%%%%%%%%%%%%%%%%%%%%%%%%%%%%%%%%%%%%%%%%%%%%%%%%%%%%%%%%%%%%%%%%%%%%%%%%%%%%%%%%%%%%%%%%%%%%%%%%%%%%%%%%%%%%%%%%%%%%%%%%%%%%%%%%%%%%%%%%%%%%%%%%%%%%%%%%%%%%%%%%%%%%%%%%%%%%%%%%%%%%%%%%